\documentclass[msom,nonblindrev]{informs3}

\def\anonymize{0}

\usepackage{amsmath, amssymb}
\usepackage{graphicx}
\usepackage{booktabs}
\OneAndAHalfSpacedXI
\usepackage{bm, bbm}
\usepackage{amsfonts}
\usepackage{xcolor}
\usepackage[margin=1in]{geometry}
\usepackage[hidelinks]{hyperref}
\usepackage{subcaption}
\usepackage[export]{adjustbox}
\usepackage{placeins}
\usepackage[nice]{nicefrac}

\newcommand{\msom}{\color{black}}
\newcommand{\latest}{\color{black}}

\usepackage{natbib}
 \bibpunct[, ]{(}{)}{,}{a}{}{,}%
 \def\bibfont{\small}%
 \def\bibsep{\smallskipamount}%
\TheoremsNumberedThrough     

\EquationsNumberedThrough    

\begin{document}


\if\anonymize0
\RUNAUTHOR{Bertsimas et al.}
\fi

\RUNTITLE{Decarbonizing OCP}

\TITLE{Decarbonizing OCP}

\if\anonymize0
\ARTICLEAUTHORS{
\AUTHOR{Dimitris Bertsimas}
\AFF{Sloan School of Management, Massachusetts Institute of Technology, Cambridge, MA, USA.\\
ORCID: \href{https://orcid.org/0000-0002-1985-1003}{$0000$-$0002$-$1985$-$1003$}\\ \EMAIL{dbertsim@mit.edu} \url{dbertsim.mit.edu}}
\AUTHOR{Ryan Cory-Wright}
\AFF{Department of Analytics, Marketing \& Operations, Imperial College Business School, Imperial College, SW7 2AZ, UK.\\
IBM Thomas J. Watson Research Center, USA\\
ORCID: \href{https://orcid.org/0000-0002-4485-0619}{$0000$-$0002$-$4485$-$0619$}\\
\EMAIL{r.cory-wright@imperial.ac.uk}, \url{ryancorywright.github.io}}
\AUTHOR{Vassilis Digalakis Jr}
\AFF{Operations Research Center, Massachusetts Institute of Technology, Cambridge, MA, USA.\\
\EMAIL{vvdig@mit.edu}, \url{vvdigalakis.github.io}}
} 
\fi

\ABSTRACT{%
\textbf{Problem definition:} We present our collaboration with the OCP Group, one of the world's largest producers of phosphate and phosphate-based products, in support of a green initiative designed to reduce OCP's carbon emissions significantly. We study the problem of decarbonizing OCP's electricity supply by installing a mixture of solar panels and batteries to minimize its time-discounted investment cost plus the cost of satisfying its remaining demand via the Moroccan national grid. {\msom OCP is currently 
designing its renewable investment strategy, using insights gleaned from our optimization model, and has pledged to invest $\$130$ billion MAD (approximately $\$13$ billion USD) in a green initiative by $2027$, a subset of which involves decarbonization.}

\textbf{Methodology/Results:} We immunize our model against deviations between forecast and realized solar generation output via a combination of robust and distributionally robust optimization. To account for variability in daily solar generation, we propose a data-driven robust optimization approach that prevents excessive conservatism by averaging across uncertainty sets. To protect against variability in seasonal weather patterns induced by climate change, we invoke distributionally robust optimization techniques. {\msom Under a ten billion MAD (approx. one billion USD) investment by OCP, the proposed methodology reduces the carbon emissions which arise from OCP's energy needs by more than $70\%$ while generating a net present value (NPV) of five billion MAD over a twenty-year planning horizon. Moreover, {\latest a twenty billion MAD investment induces a $95\%$ reduction in carbon emissions} and generates an NPV of around two billion MAD.} 

\textbf{Managerial Implications:} To fulfill the Paris climate agreement, rapidly decarbonizing the global economy in a financially sustainable fashion is imperative. Accordingly, this work develops a robust optimization methodology that enables OCP to decarbonize at a profit by purchasing solar panels and batteries. Moreover, the methodology could be applied to decarbonize other industrial consumers. Indeed, our approach suggests that decarbonization's profitability depends on solar capacity factors, energy prices, and borrowing costs.

}%



\maketitle

\section{Introduction}
On December $12$, $2015$, $195$ countries met at COP-$21$ in Paris and signed a climate agreement that aims to limit average global temperature increases to within $2$ degrees Celsius—and preferably within $1.5$ degrees Celsius—of preindustrial levels by the end of the $21$st century \citep{davenport2015inside}. As $1.1$ degrees of global temperature rises have already occurred, this pledge will be implemented by rapidly decarbonizing industrialized nations and subsequently decarbonizing the rest of the world. To minimize hardship during this transition, all participants in the global economy should decarbonize wherever doing so is profitable as soon as possible. The Kingdom of Morocco is a signatory to the Paris agreement, and in $2021$ it proposed to reduce its carbon emissions by $45.5\%$ by $2030$ \citep{du2021contribution}.

{\msom 
The OCP group (formerly Office Cherifien des Phosphates) constitutes a significant $5.6\%$ of Morocco's GDP \citep{geissler2018striving} and is responsible for managing Morocco's $70\%$ share of the world's phosphate rock reserves to manufacture fertilizer \citep{summaries2020us}. To help implement Morocco's green energy pledge, in $2022$ OCP announced a green initiative to fully decarbonize its production process by investing in solar panels and batteries which will fully power its production process by $2027$, and be fully carbon neutral by $2040$ \citep{hammond2022ocp}. This commitment is significant for two reasons. First, as of $2020$, {$57\%$ of OCP's energy needs} are satisfied via non-renewable resources, and OCP is interested in dramatically reducing this portion. Second, fertilizer is responsible for around $30\%$-$50\%$ of the world's food production \citep{stewart2005contribution}, and with a growing global population and increasingly protein-rich dietary habits, it is desirable to produce food in as sustainable a fashion as possible.
}

In this paper, we describe a robust optimization (RO) methodology that OCP is currently using to size its investment and partially implement its decarbonization pledge by installing a judicious mix of solar panels and batteries throughout its production system. We also propose techniques which allow OCP to successfully operate its system in the presence of uncertainty induced by intermittent solar generation, particularly concerning deciding when to release or store energy in batteries. By co-optimizing the cost of installing renewable energy and the cost of procuring electricity from the Moroccan national grid, we drastically curtail the carbon emissions which arise due to OCP's energy needs. We remark that there are other sources of emissions in OCP's production process that we do not address, and thus decarbonizing OCP's electricity supply does not fully decarbonize OCP; see \citet{becker2022toward} for an overview of OCP's production process.

Our model forecasts that decarbonizing OCP's electricity supply is profitable for OCP due to the abundance of solar energy in Morocco, relatively high energy prices, and relatively low real interest rates. This forecast has significant managerial implications; while companies and nation-states recognize decarbonization as a laudable goal—indeed a profitable one for the global economy \citep{adrian2022great}—they sometimes view it as an expensive luxury that members of the developing world cannot afford while industrializing their economies \citep{schmall2022india}. Indeed, the Paris climate agreement reflects this belief, by requiring that wealthier industrialized nations decarbonize sooner since they can afford to do so. However, we demonstrate by example that decarbonization can (at least sometimes) be profitable; therefore, access to global capital markets and low regulatory barriers to installing renewables may sometimes be sufficient to decarbonize. 

Ultimately, whether decarbonization is profitable for a large industrial consumer depends on the amount of solar capacity available, the price of solar panels and batteries, local energy prices, and interest rates. Therefore, we hope that our methodology can be applied elsewhere to ascertain whether decarbonization via solar panels (or other technologies) is profitable for a given consumer.

\subsection{An Overview of OCP's Production Process}\label{ssec:overview-ocp}
\begin{figure}[h!] 
    \centering
\includegraphics[width=0.8\columnwidth]{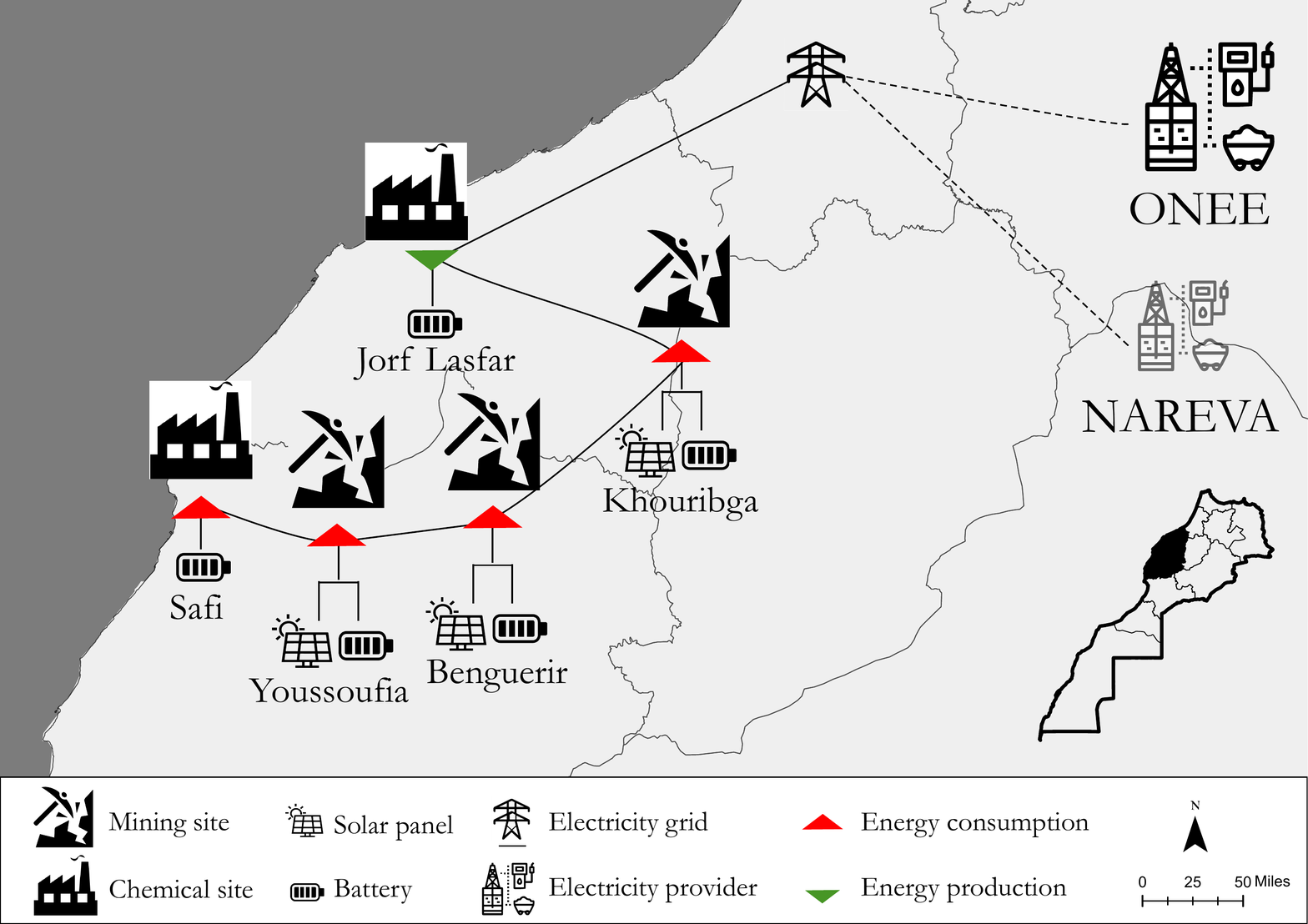}
\caption{OCP's manufacturing and phosphate rock mining sites across Morocco, and its energy suppliers. {\msom Energy production (consumption) refers to a site owned by OCP that produces (consumes) net energy. }}
\label{fig:map}
\end{figure}

In this section, we provide an overview of OCP’s current manufacturing process and its energy consumption behavior \citep[see][for further details]{ocp}, with a view to model and eventually decarbonize the process. Figure \ref{fig:map} summarizes OCP's current process visual-spatially.

\paragraph{Current Production Process (as of 2020).} Each year, OCP extracts around $40$ million tons of raw phosphate rock from eight mining sites. This extracted rock is first enriched at washing facilities, and then either exported via a port or transported to one of two chemical sites to undergo further processing. A portion of the phosphate rock that arrives at the two processing platforms is combined with sulfuric acid to produce about seven million tons of phosphoric acid, which is then directly exported to be used by a variety of economic sectors, including the food and pharmaceutical industries. The remainder is processed with ammonia to produce about $12$ million tons of fertilizers. All in all, OCP supplies more than half of all phosphate-based products sold in Africa.

\paragraph{Energy Consumption Behavior.} OCP’s energy needs are currently satisfied via carbon-emitting sources, with a small minority of its needs met via wind energy. Specifically, its demand is currently met by a combination of: wind energy ($6\%$), which OCP procures via power purchase agreements, cogeneration ($37\%$), which is generated by recovering waste heat released during the sulfuric acid production within OCP’s processing sites,  industrial fuel ($35\%$), diesel ($9\%$), natural gas ($4\%$), and electricity purchased from the grid ($9\%$). Excluding cogeneration, around $89.5\%$ of OCP's demand is currently satisfied via carbon-emitting sources.

\paragraph{Gentailers.} Excluding cogeneration, OCP purchases more than $60\%$ of the electricity which it uses. OCP’s electricity providers are the Office National de l’Electricité et de l’Eau Potable (ONEE), the national electricity company and leading operator in the field of electricity in Morocco, and Nareva, a private energy company. 
ONEE has a global installed net generation capacity of $11$ GW and generates around $34\%$ of its energy via renewable sources. 
Nareva has an installed capacity of about $2$ GW and generates around $37\%$ of its energy via renewable sources.

In this paper, we develop a methodology which reduces the amount of carbon-based electricity OCP procures from the grid by installing solar panels and batteries across OCP's system. Our methodology explicitly models (a) a complex of four mining sites in Khouribga, (b) the mining sites in Benguerir and (c) Youssoufia, as well as OCP’s two chemical sites, located in (d) Jorf and (e) Safi. Due to physical constraints, solar panels can be installed only at mining sites, whereas batteries can be installed at all sites. All sites are connected to the grid; however, the OCP-Nareva contract allows purchasing electricity from Nareva only at mining sites. We now review the relevant literature on this topic (Section \ref{ssec:related-work}), before stating our full contributions explicitly (Section \ref{ssec:contributions}).

\subsection{Related Work}\label{ssec:related-work}
Our work arises at the intersection of two related areas of the Analytics and Operations Management literature: (a) techniques for optimally expanding a production system's capacity and (b) data-driven methods for solving multi-stage RO problems. We review both areas.

\paragraph{Generator Capacity Expansion Literature.} The problem of expanding a generator's capacity is one of the most frequently studied in the power systems literature. Classically, generators rank their options according to their long-run marginal costs using a screening curve \citep[see, e.g.,][]{masse1957application, stoft2002power}, which, as discussed by \citet[]{ferris2021dynamic} corresponds to solving a linear optimization problem to determine which generators should expand at each location. 

Unfortunately, applying this approach out of the box can be inaccurate, since it assumes future supply and demand are deterministic. In reality, wind and solar generation are intermittent and could be considered as part of a future generation mix even when not currently present. Moreover, future consumer demand is uncertain since it depends on local weather conditions and population growth among other factors. As reflected in stochastic \citep{birge2011introduction, shapiro2021lectures} and robust \citep[]{ben2009robust, bertsimas2020robust} optimization textbooks, screening curves are therefore highly suboptimal in the presence of uncertainty.

To address capacity expansion {\msom more rigorously}, several authors have proposed models that explicitly account for uncertainty via stochastic or robust optimization. 
Among others, \cite{ahmed2003approximation} model capacity expansion problems as multi-stage stochastic integer programs and approximately solve them via their linear relaxations; \citet{ahmed2003multi} solve capacity expansion problems to global optimality via branch-and-bound; \citet{singh2009dantzig} solve them in a more scalable manner via Dantzig-Wolfe decomposition
; while \citet{zou2018partially} apply a very efficient adaptive optimization scheme; see \citet{gacitua2018comprehensive, guo2022generation} for reviews.

More recently, capacity expansion models have also been developed for regulators aiming to address sustainability concerns by decarbonizing electricity markets. Among others, \citet{boffino2019two} proposes a generation and transmission capacity expansion model for decarbonizing the ERCOT (Texas) electricity market by installing additional wind and solar capacity. In a related direction, \citet{ferrisa100} recently developed a model for decarbonizing the New Zealand Electricity Market while accounting for its hydro-dominated nature. Both models either constrain the amount of non-renewable energy consumed explicitly, or equivalently (from a Lagrangian perspective) impose a carbon price to reduce carbon emissions.


In this work, we take a different perspective on capacity expansion. Rather than considering the perspective of an individual generator that sells its production to consumers or a regulator aiming to reduce emissions across an entire electricity system, we take the perspective of a vertically integrated generator-consumer pair, or prosumer, which produces electricity to satisfy its own demand. Our model is also designed in collaboration with the prosumer, thus making it arguably more physically realistic than existing models in the literature. 


\paragraph{Data-driven methods for multi-stage RO.} The idea of planning over time by partitioning decision variables across multiple stages has been popular since it was originally proposed by George Dantzig and Evelyn Beale in $1955$, eventually giving rise to the paradigm of multistage stochastic optimization \citep{dantzig1955linear, beale1955minimizing}. In this paradigm, at each stage, a decision maker optimizes a subset of the decision variables while knowing the past deterministically and the future stochastically. In response, nature selects some of the uncertain parameters from a (known) probability distribution. In practice, inferring nature's joint probability distribution requires a prohibitive amount of historical data and therefore optimizers replace the true underlying distribution with its empirical distribution. This technique is known as sample average approximation (SAA) \citep[c.f.][Chapter 5]{shapiro2021lectures}. In the single-stage case, given a sample $\mathcal{S}=\{\xi_1, \xi_2, \ldots, \xi_N\}$ of a random variable $\xi(\omega)$ defined on a probability space $(\Omega, \mathcal{F}, \mathbb{P})$, an optimizer approximately minimizes the expected value of a cost function $c(\bm{x}, \bm{\xi})$, $\mathbb{E}[c(\bm{x}, \bm{\xi})]$ over a convex set $\mathcal{X}$ by solving
$
    \min_{\bm{x} \in \mathcal{X}} \frac{1}{N}\sum_{i=1}^n c(\bm{x}, \bm{\xi}_i).
$
Owing to the success of techniques like SAA and related multi-stage algorithms like Stochastic Dual Dynamic Programming \citep{pereira1991multi}, stochastic optimization is frequently and sometimes successfully used to address important problems in energy, financial and logistics.

Unfortunately, stochastic optimization has three main drawbacks. First, without making potentially heroic assumptions on the structure of a problem, it is computationally intractable. Indeed, \citet{hanasusanto2016comment} 
have shown that solving a two-stage stochastic optimization problem exactly is $\#P$-hard, i.e., as hard as counting the number of optimal solutions to a binary quadratic problem, while sample-based methods require a number of scenarios exponential in the number of stages \citep{shapiro2006complexity}. Second, stochastic optimization assumes that uncertain parameters can be drawn from a (known) probability distribution, when optimizers usually only have access to a limited amount of historical data, and estimating a potentially high-dimensional distribution from this data can introduce significant estimation error. Third, SAA tends to generate policies {\msom that} overfit the empirical distribution 
and 
perform poorly out-of-sample \citep{van2021data}.

To avoid the computational intractability of stochastic optimization, several authors \citep{soyster1973convex, ben1999robust, bertsimas2004price} designed an alternative modeling paradigm called RO, which plans under uncertainty by designing policies {\msom that} perform best under the worst-case parameter realization contained within an uncertainty set. Crucially, if uncertainty sets are designed appropriately, RO problems can be reformulated as deterministic equivalents of comparable size and similar complexity as the original nominal problem \citep[]{bertsimas2006tractable}. Moreover, as discussed in \cite{bandi2012tractable}, one can use uncertainty sets that are derived from asymptotic implications of probability theory
, and therefore obtain tangible probabilistic guarantees on the feasibility of the solution to the RO problem. 
A frequent critique of RO is that a suboptimal solution to a stochastic optimization problem may in fact perform better in practice than an optimal solution to a robust one, because RO aims to generate solutions with good worst-case performance and thus is too conservative in circumstances where we are interested in optimizing average-case performance and uncertainty can be amortized over time. 

To reduce the conservatism of RO while retaining its tractability, several authors \citep{delage2010distributionally, wiesemann2014distributionally} propose a related modeling paradigm called distributionally robust optimization (DRO), which unifies stochastic optimization with RO by optimizing for the worst-case measure over an ambiguity set of probability measures. Surprisingly, 
unlike in the stochastic case (which typically overfits), or the robust case (which typically underfits), one can easily generalize DRO to account for the fact that we work with historical data rather than probability distributions without sacrificing much out-of-sample performance. Indeed, as established by \cite{van2021data}, a variant of SAA where we pessimize over all distributions within a given distance from the empirical one, i.e., solve $\min_{\bm{x} \in \mathcal{X}} \sup_{\mathbb{Q} \in \mathcal{P}_{\delta}} \sum_{i=1}^N \mathbb{Q}(\omega^i) c(\bm{x}, \bm{\xi}^i),$
where $\mathcal{P}_{\delta}$ is the set of all measures close to the empirical measure and with the same support, guarantees that our out-of-sample disappointment will never be too high. 

\subsection{Contributions and Structure}\label{ssec:contributions}
The key contributions of the paper are twofold. First, we present an optimization \emph{methodology} that guides OCP's investment in renewable energy and {\msom supports} its operational strategy over the next $20$ years. Via RO, our proposed model guards against uncertainty associated with solar generation on a day-to-day basis. Moreover, using DRO, the model guards against changes in seasonal generation output induced by shifting weather patterns as the Moroccan climate changes with increasing atmospheric $\text{CO}_{2}$ levels. Furthermore, the model is data-driven: it combines optimization with machine learning (ML) to learn from simulated data while accounting for prediction inaccuracies. 
Second, \emph{from an implementation perspective}, our approach {\msom provides insights that facilitate} a significant (and profitable) investment in renewable energy. {\msom As of our most recent submission, OCP is using our model to size its investment in solar panels and batteries, which it expects to be at least {\msom twenty billion MAD (approx. two billion USD)}, as part of its overall $130$ billion MAD green initiative. We remark that our methodology supports a long-term investment, and therefore the full implications of our implementation will emerge over time as OCP physically installs solar panels and batteries. In terms of immediate impact, our model has already demonstrated its value as {a powerful decision support tool that OCP has used to evaluate a suite of investment options, including deciding how many solar panels and batteries should be installed in each year of a planning horizon and at each site of its production process under a given investment budget}. We detail the full impact of our collaboration on OCP's operations in Section \ref{ssec:exper-impact}.}

\paragraph{Technical Contributions.} From a technical standpoint, {(i)} we develop a multi-period SAA formulation to co-optimize OCP's investment and operational decisions in a tractable fashion. To tackle the large dimensionality of the problem, we propose a novel, ML-based scenario reduction technique, which allows us to plan for a reduced set of “typical days” rather than all days in the planning horizon. {(ii)} We use RO to account for variability in daily solar capacities (uncertainty in the weather). We propose a novel, data-driven RO approach that constructs uncertainty sets  around the amount of solar power available in {\msom each} day to reflect the idea that at the start of a day{\msom ,} a weather forecast is not perfect. This gives rise to a notion of averaging over uncertainty sets and helps prevent conservatism. {(iii)} We use DRO to protect against variability in seasonal weather patterns induced by climate change (uncertainty in the climate).

As a philosophical remark on the use of RO vs DRO: to our knowledge, there has not been work that uses both RO and DRO in the same problem. Our modeling choice of RO vs DRO follows from the amount of data we have. We use considerable amounts of historical data to inform the uncertainty in solar availability that RO guards against, but only a small amount of data to model ambiguity in seasonal weather patterns induced by climate change.

\paragraph{Managerial Implications.} Our results provide valuable insights for large electricity consumers seeking to reduce their carbon emissions or 
energy-related expenses. Our methodology demonstrates that whether decarbonizing via installing solar panels and batteries is profitable is primarily a function of local solar capacity factors, the cost of borrowing capital, and current energy prices. Therefore, although we focus on decarbonizing a fertilizer manufacturer's electricity supply, our approach could be useful in other contexts, such as manufacturing steel, cement, or paper.

\paragraph{Outline.} The rest of this paper is laid out as follows. In Section \ref{sec:SAA}, we lay out OCP's full multi-period SAA problem which allows OCP to co-optimize its investment and operational decisions. In Section \ref{sec:RO}, we develop a data-driven methodology to robustify our model against uncertainty in both day-to-day weather and seasonal climate patterns, via a combination of RO and DRO. In Section \ref{sec:exper}, we study the impact of our model and demonstrate that installing solar panels and batteries is profitable if energy costs are sufficiently high.

\section{A Sample Average Approximation Model}\label{sec:SAA}
In this section, we lay out a mathematical optimization model which optimizes OCP's green energy strategy, by making optimal \emph{strategic} decisions regarding where it should build new solar capacity and batteries, and optimal \emph{operational} decisions, which allow it to take advantage of its new assets. The key assumption in this section is that the amount of solar energy at a given site varies over time in a deterministically knowable way; we relax this assumption in Section~\ref{sec:RO}.

We first provide an overview of (Section \ref{ssec:SAA-model-description}) and formulate (Section \ref{ssec:SAA-model-formulation}) our sample-average model, which optimizes both strategic and operational decisions at the start of the planning horizon given an empirical distribution of solar capacity factors. Next, we develop a highly scalable real-time model which allows OCP to make real-time operating decisions for the next $24$ hours after receiving a solar capacity factor forecast (Section \ref{ssec:SAA-operationalizing}). Finally, we {\msom describe the scenario reduction procedure we perform to obtain a small number of highly representative scenarios for our model (Section \ref{ssec:scenario-reduction})}. All-in-all, we provide a highly accurate yet tractable approximation of OCP's full multi-period planning problem, which eventually guides its investment decisions.

\subsection{Model Description}\label{ssec:SAA-model-description}
Our overall multi-period model takes as inputs (i) the costs of procuring OCP's strategic assets, e.g., the cost of procuring and installing one kW of solar panels or one kWh of battery storage, OCP's internal cost of capital, and the salvage value of one kWh of batteries or one kW of solar panels at the end of the planning horizon, (ii) the time-dependent costs of operating OCP's system, and (iii) other information required to operate OCP's system on a day-to-day basis. It then provides as outputs (i) optimal \emph{strategic} decisions regarding where OCP should build new solar capacity and batteries at the start of each year, and (ii) an optimal \emph{policy} which allows OCP to operationalize its assets in each {\msom hour-long} period of the planning horizon. Our assumption that strategic assets are only installable at the start of each year is for computational tractability: we could consider installing assets at, e.g., the start of each month, although this would not significantly change the policies computed by the model while reducing its tractability.

{\msom  Observe that we optimize over two different time scales simultaneously: we make strategic decisions over a time horizon of years, and operational decisions over a time horizon of hours. This two-scale approach is necessary to decarbonize OCP economically: the optimal investment policy involves purchasing batteries, and the value of these batteries depends on the hour-by-hour dynamics of OCP's production process. Moreover, the capacities of batteries and solar panels degrade over a life cycle of years, which implies a trade-off between building all generation assets at the beginning of the time horizon and periodically installing assets to maintain generation capacity throughout the time horizon. 

}

\subsection{Model Formulation}\label{ssec:SAA-model-formulation}
We now formulate the capacity expansion model as a multi-period linear optimization problem, which we refer to as the sample average approximation (SAA) model. Table \ref{tab:multiperioddata} lists the parameters used throughout the model, which we set in collaboration with the OCP team{\msom , by iteratively running the models described in the paper with values prescribed by OCP, jointly examining the results, and improving the problem data and constraints until both ourselves and OCP were satisfied with the results}. {\msom We do not disclose values for all of the problem data discussed in Table \ref{tab:multiperioddata}, due to its commercially sensitive nature and to preserve OCP's privacy. }Indices that correspond to time periods are written as superscripts, and indices that correspond to locations or different electricity providers are written as subscripts. 
\paragraph{Sets and Indices.}
In the capacity expansion planning model, we are given a set of nodes $\mathcal{N}$ corresponding to mines and chemical factories which OCP operates, {a set of {\msom undirected} arcs $\mathcal{A}$ corresponding to pairs of nodes which are connected via transmission lines}, a set of hours $\mathcal{H}$, days $\mathcal{D}$, and months $\mathcal{M}$ which OCP operates its system over, and a set of years $\mathcal{Y}$ which OCP optimizes investment decisions over. Given these sets, OCP decides which strategic assets to procure, when it should procure them, and how it will operationalize them over different hours, months, and years.

{Owing to the nature of the Moroccan national grid, OCP has the option of both purchasing electricity at a given node $n$ or purchasing it at a different node $n'$ and then ``renting'' a line (or set of lines) connecting $n$ and $n'$. In this case, OCP pays the marginal price of procuring electricity at node $n'$, plus a transmission price for transmitting energy from node $n'$ to node $n$. 
}

\paragraph{\msom Scenario reduction:} To reduce the dimensionality of the problem while {\msom }treating it in a data-driven fashion, we 
{\msom  apply a clustering approach described in detail in Section \ref{ssec:scenario-reduction} to reduce the $365$ days in each year of the planning horizon to a set of scenarios $\mathcal{D}$ with a smaller cardinality (typically $|\mathcal{D}|=10$) and the same essential characteristics. 
Further, we use the same set of scenarios $\mathcal{D}$ in each month $m$, with a different weight $P^{d,m,y}$ to account for the fact that sunrise and sunset times vary predictably by month, and there are different weather patterns in different months. This allows us to represent each month of the planning horizon with the same set of scenarios.
}



\begin{table}[h!]
\centering
\caption{Summary of notation. Calligraphic letters refer to sets, Roman/Greek letters refer to problem data. {\msom To preserve OCP's privacy, we do not disclose the data values not explicitly stated in this table}.}
{\footnotesize
\begin{tabular}{l l}
    \toprule
     \textbf{Symbol} & \textbf{Description}\\\midrule
     \multicolumn{2}{l}{\textbf{General Setting}}\\
    $\mathcal{H}$ & Set of hours in each day, $\{1, \ldots, 24\}$\\
    $\mathcal{D}$ & Set of reduced scenarios\\
    $\mathcal{M}$ & Set of months in a calendar year, $\{1, \ldots, 12\}$\\
     $\mathcal{Y}$ & Set of years in OCP's planning horizon; i.e., $\{1, \ldots, 20\}$\\
     $\mathcal{N}$ & Set of nodes in the network, i.e., $\{$Jorf, Safi, Benguerir, Youssoufia, Khouribga$\}$\\
    $\mathcal{A}$ & Set of all arcs in the network\\
    \msom  $\text{D}^{m,y}$ & Number of days in month $m \in \mathcal{M}$ in year $y \in \mathcal{Y}$\\
    \msom  $\text{P}^{d,m,y}$ & Probability of reduced scenarios of type $d \in \mathcal{D}$ in month $m \in \mathcal{M}$ in year $y \in \mathcal{Y}$\\
    \multicolumn{2}{l}{\textbf{Investment Decisions}}\\
    $B$ & Investment budget in MAD (Moroccan dirham)\\
    $\rho$ & Discount factor, i.e., {\msom $0.96$}\\
    $c_b^y$ & Cost of purchasing and installing one kWh of batteries in year $y \in \mathcal{Y}$\\
    $c_s^y$ & Cost of purchasing and installing one kW DC of solar panels in year $y \in \mathcal{Y}$\\

     \multicolumn{2}{l}{\textbf{OCP Operations}}\\
    \multicolumn{2}{l}{\textit{Operational Data}}\\
    $R$ & Constant which converts {energy} released from batteries into a rate, i.e., $1$\\
    $\xi$  & Annual rate of solar generation capacity degradation, i.e., $0.995$\\
    $\nu$ & Annual rate of battery storage degradation, i.e., $0.96$\\
    $\psi$ & Proportion of {energy} stored in a battery available an hour later, i.e., {\msom $0.997$}\\
    $\beta$ & Fraction of daily amount of {power} produced by solar panels that may be sold, i.e., 0.2\\
    $\mathcal{I}(n)$ & Set of arcs $a=(i,n)$ flowing into node $n \in \mathcal{N}$\\
    $\mathcal{O}(n)$ & Set of arcs $a=(n,i)$ flowing out of node $n \in \mathcal{N}$\\[1ex]
    $K_a$ & Capacity limit {in kW} on the flow through arc $a \in \mathcal{A}$\\[1ex]
    $\eta_{a}$ & The transmission efficiency coefficient for arc $a \in \mathcal{A}$, i.e., $\eta=0.99$\\
    \multicolumn{2}{l}{\textit{Time-Dependent Data}}\\
    $G_o^{h}$ & ONEE generation capacity {in kW} at node $n \in \mathcal{N}$ in hour $h \in \mathcal{H}$\\
    $G_n^{h}$ & NAREVA generation capacity {in kW} at node $n \in \mathcal{N}$ in hour $h \in \mathcal{H}$\\
    $v^{h,d}$ & Capacity factor for a solar panel in hour $h \in \mathcal{H}$ of scenario $d \in \mathcal{D}$\\
    $p_{O,n}^{h,m}$ & Marginal cost of energy {in MAD/kWh} from ONEE at node $n \in \mathcal{N}$ at time $h \in \mathcal{H}$, $m \in \mathcal{M}$\\[1ex]
    $p_{N,n}^{h,m}$ & Marginal cost of energy {in MAD/kWh} from NAREVA at node $n \in \mathcal{N}$ at time $h \in \mathcal{H}$, $m \in \mathcal{M}$\\[1ex]
    $p_{w,n}^{h,m}$ & Marginal feed-in price {in MAD/kWh} for selling electricity at node $n \in \mathcal{N}$ at time $h \in \mathcal{H}$, $m \in \mathcal{M}$\\[1ex]
    $c_{r,a}^{h,m}$ & Marginal cost {in MAD/kWh} of renting line $a \in \mathcal{A}$ at time $h \in \mathcal{H}$, $m \in \mathcal{M}$\\[1ex]
    $d^{h,m,y}_{n}$ & Aggregate demand {in kWh} at node $n \in \mathcal{N}$ at time $h \in \mathcal{H}$, $m \in \mathcal{M}$, $y\in \mathcal{Y}$\\[1ex]

     \bottomrule  
\end{tabular}
}
\label{tab:multiperioddata}
\end{table}

\paragraph{Decision Variables.} {We now lay out the strategic decision variables optimized by the model, followed by the operational decisions, which are made on an hourly basis.}

\paragraph{Strategic Decision Variables.} The key strategic decision variables are:
\begin{subequations}
\begin{align}
    b^y_{n} \in \mathbb{R}_{+}: \ & \text{The number of kWh of battery storage built at node} \ n \in \mathcal{N} \ \text{in year}\ y \in \mathcal{Y};\\[-2pt]
    z^y_{n} \in \mathbb{R}_{+}: \ & \text{The number of kW of solar panels built at node} \ n \in \mathcal{N} \ \text{in year}\ y \in \mathcal{Y}.
\end{align}
\end{subequations}
{\msom 
Note that, as discussed in Section \ref{ssec:overview-ocp}, we cannot build solar panels at the Jorf and Safi factory sites, as these sites are in populated areas. Accordingly, we require that $z^y_{n}=0$ for $n \in \{\text{Jorf}, \text{Safi}\}$.
}

\paragraph{Operational Decision Variables.}
The key operational decision variables we use are:
\begin{subequations}
\begin{align}
     f_{a}^{h,d,m,y}: \ & \text{The DC load-flow (in kW) through line} \\[-2pt]
     &  a \in \mathcal{A} \ \text{in} \ \text{hour} \ h \in \mathcal{H} \ \text{of} \ \text{day} \ d \in \mathcal{D} \ \text{of month} \ m \in \mathcal{M} \ \text{of year} \ y \in \mathcal{Y};\nonumber\\[-2pt]
     s_n^{h,d,m,y}: \ & \text{Energy stored (in kWh) in batteries at node} \\ 
     &  n \in \mathcal{N} \ \text{in} \ \text{hour} \ h \in \mathcal{H} \ \text{of} \ \text{day} \ d \in \mathcal{D} \ \text{of month} \ m \in \mathcal{M} \ \text{of year} \ y \in \mathcal{Y};\nonumber\\
     r_{n}^{h,d,m,y}: \ & \text{{Power} (in kW) discharged at node} \\ 
     &  n \in \mathcal{N} \ \text{in} \ \text{hour} \ h \in \mathcal{H} \ \text{of} \ \text{day} \ d \in \mathcal{D} \ \text{month} \ m \in \mathcal{M} \ \text{year} \ y \in \mathcal{Y};\nonumber\\[-2pt]
     x_{O,n}^{h,d,m,y}: \ & \text{{Power} procured (in kW) from grid from ONEE at node} \\ 
     &  n \in \mathcal{N} \ \text{in} \ \text{hour} \ h \in \mathcal{H} \ \text{of} \ \text{day} \ d \in \mathcal{D} \ \text{month} \ m \in \mathcal{M} \ \text{year} \ y \in \mathcal{Y};\nonumber\\[-2pt]
     x_{N,n}^{h,d,m,y}: \ & \text{{Power} procured (in kW) from grid from NAREVA at node} \\ 
     &  n \in \mathcal{N} \ \text{in} \ \text{hour} \ h \in \mathcal{H} \ \text{of} \ \text{day} \ d \in \mathcal{D} \ \text{month} \ m \in \mathcal{M} \ \text{year} \ y \in \mathcal{Y};\nonumber\\[-2pt]
     w_{n}^{h,d,m,y}: \ & \text{{Power} sold (in kW) to grid at node} \\ 
     &  n \in \mathcal{N} \ \text{in} \ \text{hour} \ h \in \mathcal{H} \ \text{of} \ \text{day} \ d \in \mathcal{D} \ \text{month} \ m \in \mathcal{M} \ \text{year} \ y \in \mathcal{Y};\nonumber
\end{align}
\end{subequations}

Optimizing the above decision variables using the problem data in Table \ref{tab:multiperioddata} yields the problem:
\begin{subequations}\label{prob:detmultiperiod}
\begin{align}
    \min \quad & 
    \sum_{y \in \mathcal{Y}}
    \left[\underbrace{\sum_{n \in \mathcal{N}}c_b^y({\msom (\rho)^y}-{\msom (\rho)^{|\mathcal{Y}|}}) b_n^y}_{\text{cost of batteries}}
    +\underbrace{\sum_{n \in \mathcal{N}}c_s^y({\msom (\rho)^y}-{\msom (\rho)^{|\mathcal{Y}|}}) z_n^y}_{\text{cost of solar}}
    +\underbrace{\sum_{a,m,d,h} {\msom (\rho)^y} {\msom \text{D}^{m,y} \text{P}^{d,m,y} }c_{r,a}^{h,m} \vert f_{a}^{h,d,m,y}\vert 
    }_{\text{cost to rent lines}}
    \right. \nonumber\\
    & \left.
    +\underbrace{\sum_{n, m, d, h}{\msom (\rho)^y} {\msom \text{D}^{m,y} \text{P}^{d,m,y}} \left(
     p_{O,n}^{h,m} x_{O,n}^{h,d,m,y} 
     + p_{N,n}^{h,m} x_{N,n}^{h,d,m,y}
     - p_{w,n}^{h,m} w_{n}^{h,d,m,y}
    \right)
    }_{\text{cost to procure and sell energy}}
    \right] \label{eqn:obj}
\intertext{which is to be minimized subject to the following constraints:
}
 \text{s.t.} \quad & \msom  \sum_{n,y}
    c_b^y {\msom (\rho)^y} b_n^y + c_s^y {\msom (\rho)^y} z_n^y \leq B, & \label{eqn:budget} \\
  & \sum_{a \in \mathcal{I}(n)} \tau_a(f_{a}^{h,d,m,y})+ \sum_{a \in \mathcal{O}(n)} \tau_a(-f_{a}^{h,d,m,y})+ R\cdot r_{n}^{h,d,m,y}+ x_{O,n}^{h,d,m,y}+x_{N,n}^{h,d,m,y} & \nonumber\\
&  \geq d_{n}^{h,d,m,y} {\msom +} w_{n}^{h,d,m,y}-v^{h,d} \left( \sum_{y'=1}^y \xi^{y-y'} z_n^{y'} \right) , \label{eqn:kirchhoff}\\ 
 &\sum_{h} w_{n}^{h,d,m,y} \leq \beta  \sum_{h} \left[ v^{h,d} \left( \sum_{y'=1}^y \xi^{y-y'} z_n^{y'} \right) + \max\left\{ 0, -d_n^{h,m,y} \right\} \right] , \label{eqn:sell} \\ 
& s^{h+1,d,m,y}_n=\psi s_n^{h,d,m,y}-r_{n}^{h,d,m,y} ,  \ s^{1,d,m,y}_n=\psi s_n^{24,d,m,y}-r_{n}^{24,d,m,y} ,  \nonumber\\ 
& s^{h,d,m,y}_n \leq \sum_{y'=1}^y \nu^{y-y'} b^{y'}_n ,\label{eqn:batteries}\\ 
& x_{O,n}^{h,d,m,y} \leq G_o^h, \ x_{N,n}^{h,d,m,y} \leq G_n^h ,\label{eqn:gridcapacity} \\ 
& \vert f_{a}^{h,d,m,y}\vert \leq K_a ,\label{eqn:flow}\\[-2pt] 
& s_{n}^{h,d,m,y}, x_{O,n}^{h,d,m,y},
     x_{N,n}^{h,d,m,y},
     w_{n}^{h,d,m,y}, b_n^y, z_n^y \geq 0 , \nonumber\\ & z^y_{n}=0 \ \text{for} \ n \in \{\text{Jorf}, \text{Safi}\}.\nonumber
\end{align}
\end{subequations} 
{where $\tau_a(f):=\eta_a f_{+}-f_{-}: f=f_{+}-f_{-}, f_{+}, f_{-}\geq 0$ models the net flow through arc $a$ after transmission losses, and all constraints are taken over the indices $n, a, h, d, m, y$ whenever these indices are present in a constraint and not prescribed by a sum.}

\paragraph{Objective.} Equation \eqref{eqn:obj} minimizes the time-discounted cost of OCP's investment plus the time-discounted cost of procuring electricity from the Moroccan national grid. Correspondingly, we minimize the sum of four terms which each correspond to annual strategic or operational decisions and are discounted by year {\latest via a discount factor $\rho=0.96$ for a discount rate of $4\%$}. {\msom Note that $(\rho)^y$ denotes the discount factor raised to a power $y$ and all non-bracketed superscripts denote indices.} 

The first term models the cost of installing batteries, minus {their end-of-horizon salvage value}. Similarly, the second term represents the cost of installing solar panels {minus their salvage value}. The third term corresponds to the operational cost of renting power lines. We take the planning horizon to be $T=20$ years and assume that batteries and solar panels degrade over time. Finally, the fourth term captures the operational cost of procuring energy from the Moroccan power grid (from either provider), minus the profit from selling energy to the grid.

\paragraph{Constraints.} To ensure that the model is physically realistic, it contains a number of constraints.
Constraint \eqref{eqn:budget} imposes that the total {\msom time-discounted} investment cannot exceed the allocated budget. 
Constraint \eqref{eqn:kirchhoff} imposes Kirchoff's current law at each node $n \in \mathcal{N}$, after accounting for transmission losses, batteries, and solar generation. Note that we treat power generation output at each factory as negative demand, and impose an inequality constraint to allow factories/solar panels to shed load, as this will rarely be profitable at optimality. {\msom Also note that the solar capacity factors $v^{h,d}$ do not explicitly depend on the month of the year; this dependence is captured implicitly, by placing a higher weight $P^{d,m,y}$ on scenarios $d$ in month $m$ which more closely align with the historical solar capacity factors in month $m$, as described in detail in Section \ref{ssec:scenario-reduction}.}
Constraint \eqref{eqn:sell} imposes that, daily, we cannot sell more than a pre-specified fraction (here $20\%$) of the energy we locally produce at a given node (either via solar panels or via negative demand). This constraint is imposed both for legal purposes and to prevent energy providers from viewing OCP as a competitor in the energy generation business. 
Constraint \eqref{eqn:batteries} ensures that, for each day, the amount of energy stored in batteries at hour $h+1$ is equal to that at hour $h$, plus any extra energy injected into the battery, minus any injected into the grid, after accounting for efficiency losses. {Moreover, the amount of energy stored in batteries at any time cannot exceed capacity or be negative. For} the sake of tractability, we assume that energy stored at hour $24$ can be used at hour $1$. Constraint \eqref{eqn:gridcapacity} ensures that the electricity procured from the grid cannot exceed its capacity. 
Finally, Constraint \eqref{eqn:flow} ensures that net flow through a line does not exceed the line's capacity.

{\msom We remark that our model implicitly assumes that the amount of energy stored in batteries at the end of hour $24$ of each scenario is the amount of energy available in hour $1$ of the same scenario, and thus load cannot be shifted between different days in OCP's planning horizon. It could be argued that this assumption overly simplifies OCP's investment problem. Fortunately, this is not the case, because it is not profitable to shift a significant amount of load between days. Indeed, Morocco's weather is very consistent on a month-to-month basis, and energy stored in batteries depletes over time, which strongly suggests that it is more profitable to use energy the day it is collected than to store it in batteries for multiple days. {\latest In Section \ref{appx:exper-sensitivity}, we verify this claim empirically, by investigating the impact of partially relaxing this assumption on the optimal investment policy. Namely, we take the outer product of our set of scenarios with itself to generate $48$ hour-long scenarios where it is possible to shift load between the first and second day in each scenario. Following the same logic, we also generate $72$ hour-long scenarios. We find that using $48$ or $72$, rather than $24$ hour-long scenarios has a negligible impact on the optimal investment. Rather, the benefit of installing batteries is that they allow OCP to store energy when the sun shines, and use this energy at night time.}}
 
 


\subsection{Operationalizing the SAA Model} \label{ssec:SAA-operationalizing}
In this section, we convert the SAA model we derived in the previous section into one which OCP can use to make real-time decisions. This is an important practical issue. Indeed, the SAA model can be run at the beginning of each year of the planning horizon to make strategic decisions regarding how many solar panels and batteries OCP should install that year and where it should install them, but does not actually provide OCP with any guidance on how it should operate its system on an hourly basis. In particular, the SAA model proposes a daily policy for each reduced scenario in our set of scenarios, but, since the amount of sunlight is continuous and uncertain, real days will almost surely differ from these scenarios.

To resolve this issue, {\msom in Section \ref{appx:SAA-operationalizing}, we develop a real-time linear optimization model which, for a given time $t$ in a given day $d$ of year $y$, and a given—highly accurate—solar capacity forecast for the rest of the day, prescribes an optimal policy for how to operate OCP's system across the rest of the day. 
Notice that, at time $t$, decision variables that correspond to times $i \leq t-1$ have already been decided (we use tilde to indicate such quantities). Fixing this data in our SAA model gives the real-time model {\msom laid out in} Problem \eqref{eqn:operational}; a linear optimization problem with $O\big(T(|\mathcal{N}|+|\mathcal{A}|)\big)$ variables which can be solved at scale in real-time.}

OCP operates its system in time period $t$ by acting according to an optimal value of the decision variables $(f_{a}^{t}, s_{n}^{t}, r_{n}^{t}, x_{O,n}^{t}, x_{N,n}^{t}, w_{n}^{t})$ in Problem \eqref{eqn:operational}. Respectively, these variables correspond to the amount of power we transmit through each line, the amount of energy we store at each node's batteries, the amount of energy we charge or discharge from each node's batteries, the amount of power we procure from each provider at each node, and the amount of power we sell back to the grid in time period $t$.
We remark that if the real-time solar capacity factors in time period $t+1$ deviate from their forecast value in time period $t$, then we can rerun the real-time model in period $t+1$ and proceed accordingly.

{We close this section by noting that in the same way as solar forecasts deviates over a period of days, their distribution deviates over a period of years, due to, e.g., climate change. To address this, OCP should also rerun the SAA model at the start of each year (with the number of solar panels and batteries previously purchased at each site fixed as sunk costs). This approach is potentially useful for a second reason: the analysis in this paper is built upon one year of hourly solar capacity data, since this is all OCP currently has access to. Rerunning the SAA model at the start of each year allows OCP to iteratively incorporate more data in its forecasts, thus driving the implemented solutions closer to that we would obtain with full knowledge of the underlying stochastic process.}


\subsection{Solar Capacity Factor Scenario Reduction} \label{ssec:scenario-reduction}
{\msom In this section, we outline our scenario reduction approach, which leverages clustering techniques to reduce the number of days in each year of our planning horizon. The primary motivation for our approach is that optimizing over a small representative set of $|\mathcal{D}|$ days instead of all $365$ days in each year of the horizon reduces the dimensionality of the problem and our policies' sensitivity to outliers while obtaining qualitatively similar policies in practice. Indeed, accounting for the robustification we perform in the next section, there are strong theoretical guarantees that our clustering approach does not significantly alter the optimal investment compared to the full sample-average-approximation \citep[see][for a general theory]{wang2022mean}. Moreover, as we show in Section \ref{appx:exper-sensitivity}, our investment policy does not change significantly as we vary the number of clusters.

\paragraph{Simulated Solar Capacity Factors.} We are given a year of simulated solar capacity factors for each OCP mining site. The simulated solar capacity factors are obtained using \texttt{PVsyst 7.2}, a software package for the study, sizing, and data analysis of complete PV systems. OCP engineers performed the corresponding simulations using several years of solar generation data from all OCP sites. To generate solar capacity factors over a longer time horizon, we combine the simulated solar capacity factors with forecasts of solar generation in North Africa over the next 20 years generated by \citet[][]{jerez2015impact,bichet2019potential}. In particular, we take a conservative approach and assume solar generation will degrade by a worst-case $1.5 \%$ in total over the next 20 years. We let $i \in \{1, \ldots, 20 \cdot 365\}$ index our data points, which are $72$-dimensional, and represent the hourly solar capacity factor on a given day for all three mining sites. Correspondingly, the solar capacity factor at mining site $n$ in hour $h \in \mathcal{H}$ of data point $i \in [N]$ at all three mining sites is given by ${v}_n^{h,i}$.


\paragraph{Clustering-based Scenario Reduction With Wasserstein Distance Metric.} We reduce the $20\cdot 365$ days of capacity factor data into a smaller set of scenarios $\mathcal{D}$ according to the following methodology: First, we perform $k$-means clustering and take the centroids of each cluster as the solar capacity factors in our scenarios, denoting by $\bar{\bm{v}}^{d}$ the solar capacity factor of each reduced scenario $d \in \mathcal{D}$. 
We remark that rather than $k$-means clustering, one could alternatively use another clustering technique to perform scenario reduction, e.g., $k$-medoids or hierarchical clustering; see \citet[][]{teichgraeber2019clustering} for a review of clustering techniques for scenario reduction in energy systems. However, in the absence of the distributional robustness techniques explored in the next section, \cite{teichgraeber2019clustering} found that $k$-means clustering tends to perform more predictably than other clustering approaches reviewed therein, which arguably justifies our approach. 

Second, to account for monthly deviations in solar capacities induced by seasonal weather patterns and different sunrise and sunset times, we assign a different mass $\text{P}^{d,m}$ on scenario $d$ in each month $m$. We note that our technique can be extended to assign different weights $\text{P}^{d,m,y}$ in each year; however, we avoid this as one year of simulated capacity factors is not enough to accurately estimate annual deviations in solar capacity, beyond the $1.5\%$ degradation assumed over the time horizon. 
To compute the mass $\text{P}^{d,m}$ assigned to each scenario, we follow the standard Wasserstein mass transportation approach \citep[see, e.g.,][]{rujeerapaiboon2022scenario}. Namely, we set $\textbf{P}^{m}$ in each month by minimizing the type-$1$ Wasserstein distance between the empirical distribution of capacity factors in that month, and the set of all discrete distributions supported on our reduced scenarios $\mathcal{D}$. For a set of days $j \in \mathcal{J}_m$ in month $m$, this corresponds to solving the linear problem:
\begin{equation} \label{eqn:scenario-reduction}
\begin{aligned} 
\min \quad & \sum_m \sum_{i \in \mathcal{J}_m} \sum_{d \in \mathcal{D}} W^{d, i, m} \left\| \bar{\bm{v}}^d - \bm{v}^{i} \right\|^2 & \\
\text {s.t.} \quad & \sum_{d \in \mathcal{D}} W^{d, i, m}=1, \quad & \forall m,\ i \in \mathcal{J}_m \\
& \text{P}^{d,m} = \frac{1}{|\mathcal{J}_m|} \cdot \sum_{i \in \mathcal{J}_m} W^{d, i, m}, \quad & \forall m,\ d \in \mathcal{D} \\
& \bar{\bm{v}} \geq 0,\ \textbf{P} \geq 0,\ \bm{W} \geq 0. & 
\end{aligned}    
\end{equation}


{\latest
In Section \ref{appx:scenarios-visualization}, we provide a visualization of the proposed scenario reduction methodology.
Moreover, in the next section, we explore techniques for robustifying our point estimates $\bm{v}^d$ to account for uncertainty in renewable generation output, as clustering-based scenario reduction techniques may perform arbitrarily poorly in the absence of robustification techniques \citep[c.f.][]{rujeerapaiboon2022scenario}, but perform well with robustification \citep{wang2022mean}.
}}

\section{Immunizing Against Uncertainty in Solar Generation} \label{sec:RO}

{\latest When we introduced our SAA model in Section \ref{sec:SAA}, we assumed that the amount of sunlight in {\msom each} day is drawn from a probability distribution at the start of the day, and is subsequently deterministic.} Unfortunately, this is not true in practice: when OCP plans its operations for a given day, it has access to a highly accurate forecast of the amount of sunlight available, but this forecast is not deterministically accurate. To address this, OCP solves a real-time problem in each hour of the day, as suggested in Section \ref{ssec:SAA-operationalizing}. However, this strategy is suboptimal for two reasons. First, the SAA problem does not account for deviations between the SAA and real-time solutions. Second, the real-time problem in each time $t$ does not account for further deviations {\msom in} the weather.

{\latest From a RO perspective, the SAA model also suffers from a further defect: the amount of solar generation in each day of the planning horizon is drawn from an underlying probability distribution; we only have access to a sample from this distribution, but not the distribution itself.} Thus, the SAA model overfits the simulated data in a finite-sample setting and underperforms out-of-sample.


Accordingly, in this section, we robustify the model in two ways. First, we address the deviation between the SAA and real-time problems: we outline how we model uncertainty in solar generation in Section \ref{ssec:RO-uncertainty-defn}, perform an exploratory analysis of OCP's solar capacity factor data in Section \ref{ssec:RO-uncertainty-analysis}, use our analysis to propose a data-driven uncertainty set in Section \ref{ssec:RO-uncertainty-set}, and robustify the SAA model against deviations from forecasts in Section \ref{ssec:RO-constraint-aggregation}. Second, we immunize the model against overfitting to the simulated data by using recent advances in DRO in Section \ref{ssec:RO-distributionally}. Interestingly, this robustification does not impact the SAA problem's tractability, since the robust model can be reformulated as a tractable linear optimization model with exponential cone constraints, featuring a number of decision variables and constraints proportional to the SAA model. All in all, our approach immunizes the SAA model against uncertainty in a tractable and unified fashion.

{\msom As a philosophical remark, we note that decarbonization must occur in many different contexts to fulfill the Paris climate treaty, and the nature of uncertainty might be quite different in some of these other contexts. For instance, industries near large amounts of hydroelectric power might install hydroelectric dams to decarbonize, in which case the nature of the uncertainty may be different \citep[c.f.][]{duque2020distributionally}. In this case, it may be appropriate to use different uncertainty sets than the sets designed here in conjunction with our overall approach, rather than applying the uncertainty sets designed here directly out-of-the-box.}

\subsection{Uncertainty in Solar Generation} \label{ssec:RO-uncertainty-defn}

In the reduced SAA formulation \eqref{prob:detmultiperiod} laid out in the previous section, we assumed that the solar capacity factor $v^{h,d}$ for hour $h \in \mathcal{H}$ of scenario $d\in \mathcal{D}$ is unknown at the start of each year of the planning horizon, and takes the value $v^{h,d} := \bar{v}^{h,d}$ before OCP plans its operations in scenario $\mathcal{D}$. This model of uncertainty is analogous to OCP being given a probability distribution of capacity factors when making its asset purchases at the start of each year, and subsequently receiving a perfectly accurate forecast at the beginning of each day $d$, for the entirety of day $d$, which allows OCP to plan its operations for that day without accounting for uncertainty.

In practice, the true solar capacity factors at each node $n$ and time $(h,d,m,y)$ are highly correlated to the forecasts OCP receives at time $(1, d, m, y)$, but not exactly equal to them. To account for this discrepancy, we now introduce a vector of uncertain parameters $\boldsymbol{u}$ and require that our OCP's operational policy is feasible for all capacity factors $\boldsymbol{v}$ in the set
\begin{equation}
    \mathcal{U} = \left\{ \boldsymbol{v} \in \mathbb{R}_+^{|\mathcal{H}| \times |\mathcal{D}| \times |\mathcal{M}| \times |\mathcal{Y}| \times |\mathcal{N}|}: \quad v^{h,d,m,y}_n = \bar{v}^{h,d}+u^{h,d,m,y}_n \quad \forall h,d,m,y,n \quad \forall \boldsymbol{u} \in \mathcal{U}_u \right\}.
\end{equation}

{\color{red}


}

In the rest of this section, we describe how we model the uncertainty in solar generation through a tailored, data-driven uncertainty set $\mathcal{U}_u$, which accounts for adversarial perturbations in the nominal solar generation, albeit in a controlled and structured way.




\subsection{Exploring the Uncertainty}\label{ssec:RO-uncertainty-analysis}
We now study the predictability of the solar capacity factors to, in the next section, design a data-driven uncertainty set that accurately reflects uncertainty in renewable generation output without being excessively conservative. Correspondingly, the main quantity which we are interested in is the within-reduced-scenario distance between an average hourly solar capacity factor and a given simulated capacity factor. Therefore, we explore how this quantity evolves across time in order to propose an uncertainty set which accurately models deviations in solar generation capacity. For the analysis that follows, we take $|\mathcal{D}|=10.$

\begin{figure}[!ht]
\centering
\begin{subfigure}{.45\textwidth}
  \centering
  \includegraphics[width=\linewidth]{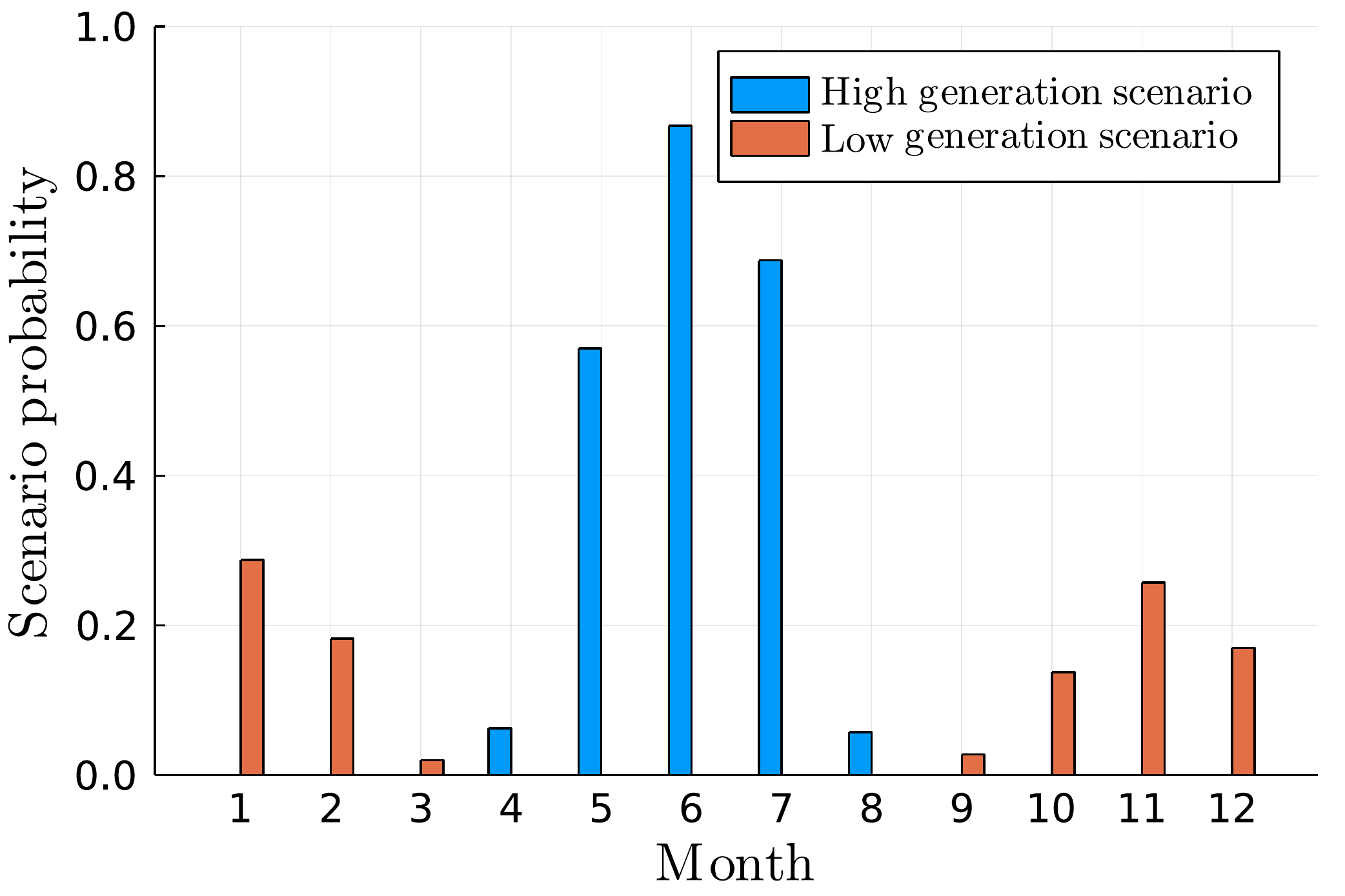}
\end{subfigure}%
\hfill
\begin{subfigure}{.45\textwidth}
  \centering
  \includegraphics[width=\linewidth]{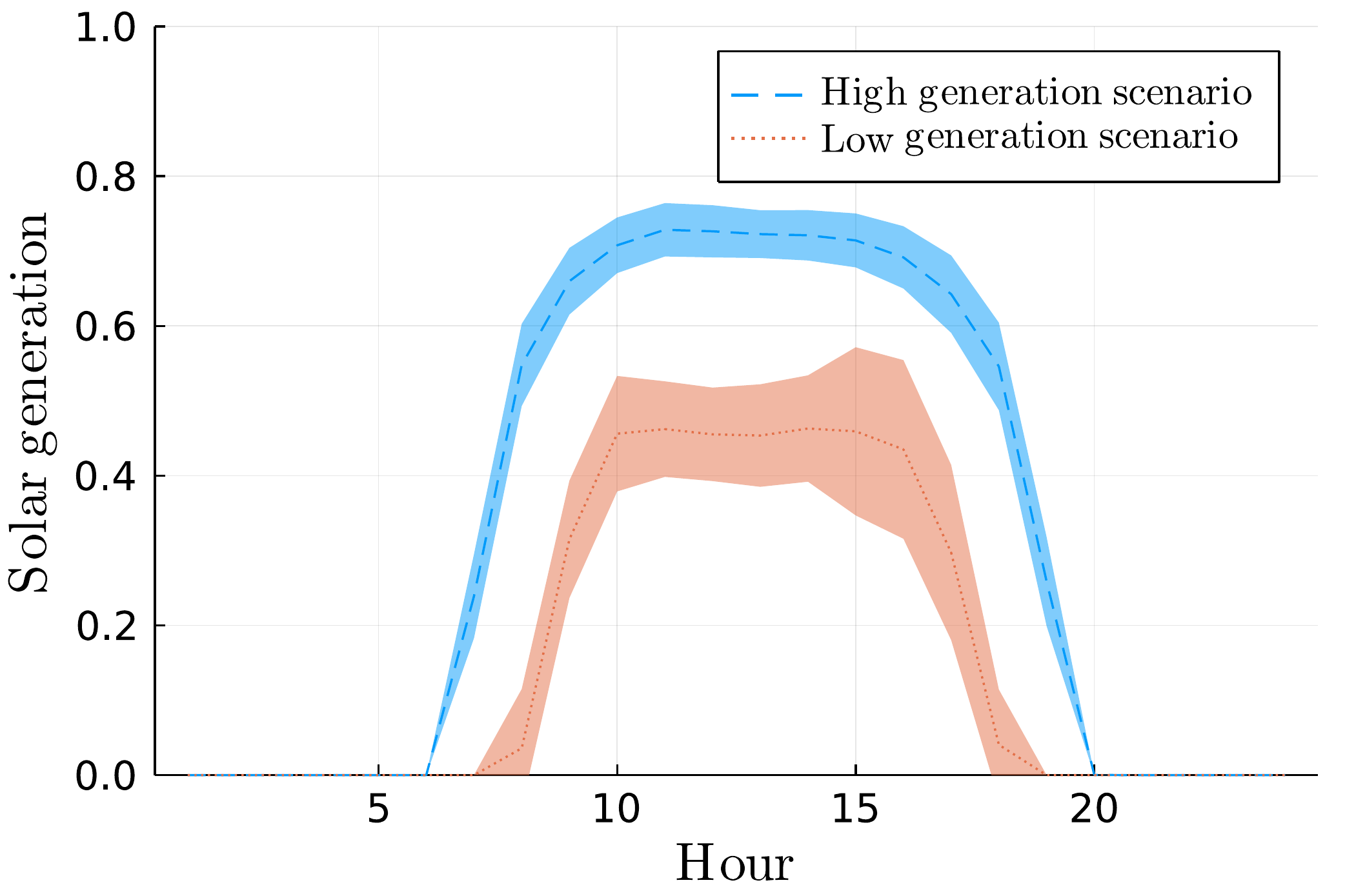}
\end{subfigure}%
\vspace{10pt}
\caption{Probability of occurrence during each month of selected reduced scenarios (left) and mean and standard deviation of solar generation for selected reduced scenarios (right).}
\label{fig:unc1}
\end{figure}

{\latest We first study the two most common reduced scenarios obtained by clustering our solar capacity factors according to the methodology described in Section \ref{ssec:scenario-reduction}. As shown in Figure \ref{fig:unc1} (left), the first reduced scenario (RS1) is a high-generation scenario that occurs mostly during the spring and summer months, while the second reduced scenario (RS2) is a low-generation scenario that occurs mostly in winter and fall months.} 
Figure \ref{fig:unc1} (right) presents, for both RS1 and RS2, the mean and standard deviation in the hourly solar generation capacity factors. The mean corresponds to the nominal solar generation $\bar{v}^{h,d}$, while the standard deviation varies depending on the hour $h$ of the day and is typically zero during nighttime. Unsurprisingly, the mean solar generation for RS1 is higher. More interestingly however, the variance for RS2 is higher than RS1. This is perhaps best explained by the more variable cloudy weather conditions that RS2 captures and suggests that the structure of uncertainty varies between scenarios and should be modeled as such.

\begin{figure}[h!]
\centering
\begin{subfigure}{.45\textwidth}
  \centering
  \includegraphics[width=\linewidth]{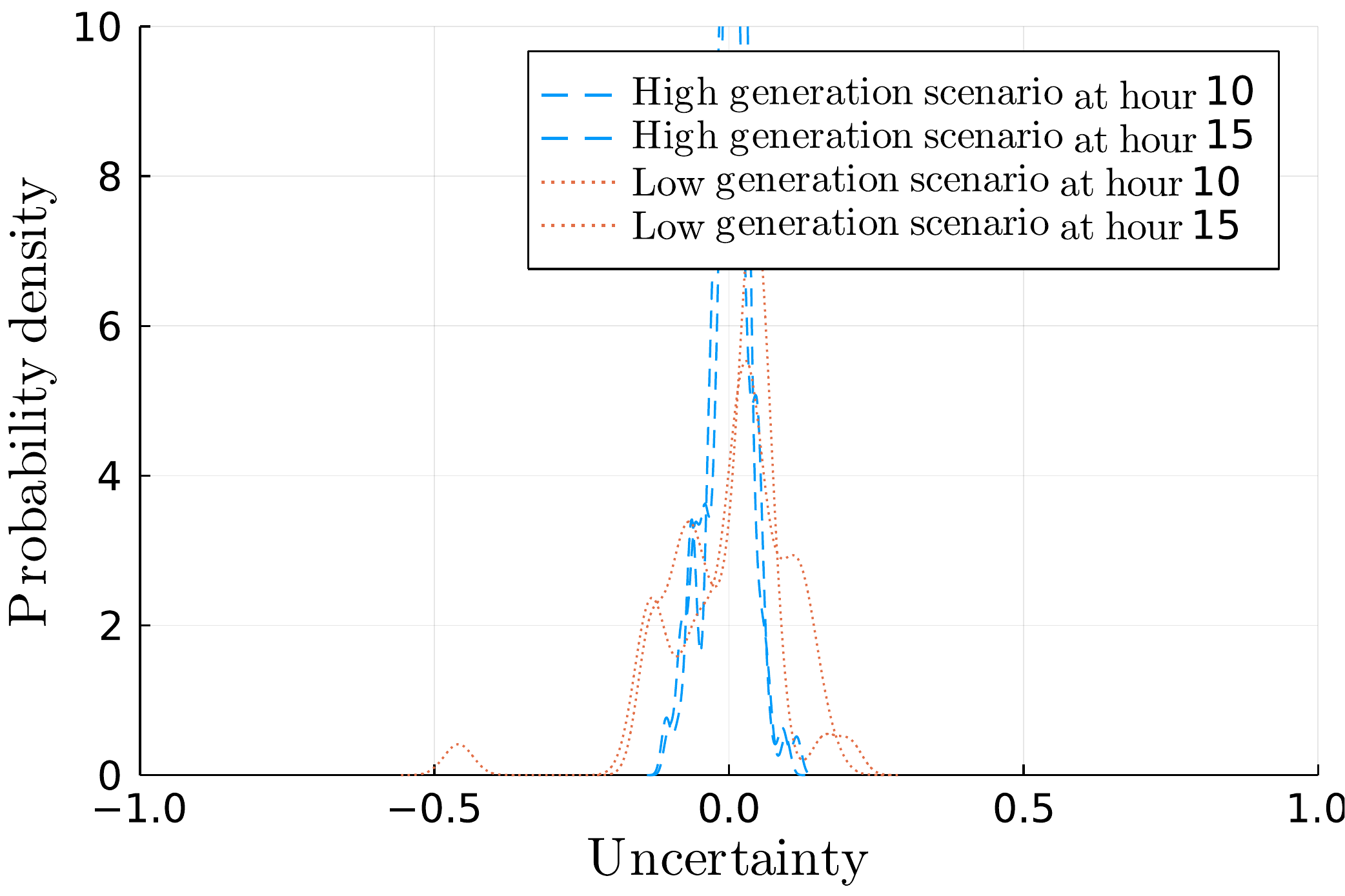}
  \label{fig:plot_unc_clt}
\end{subfigure}%
\hfill
\begin{subfigure}{.45\textwidth}
  \centering
  \includegraphics[width=\linewidth]{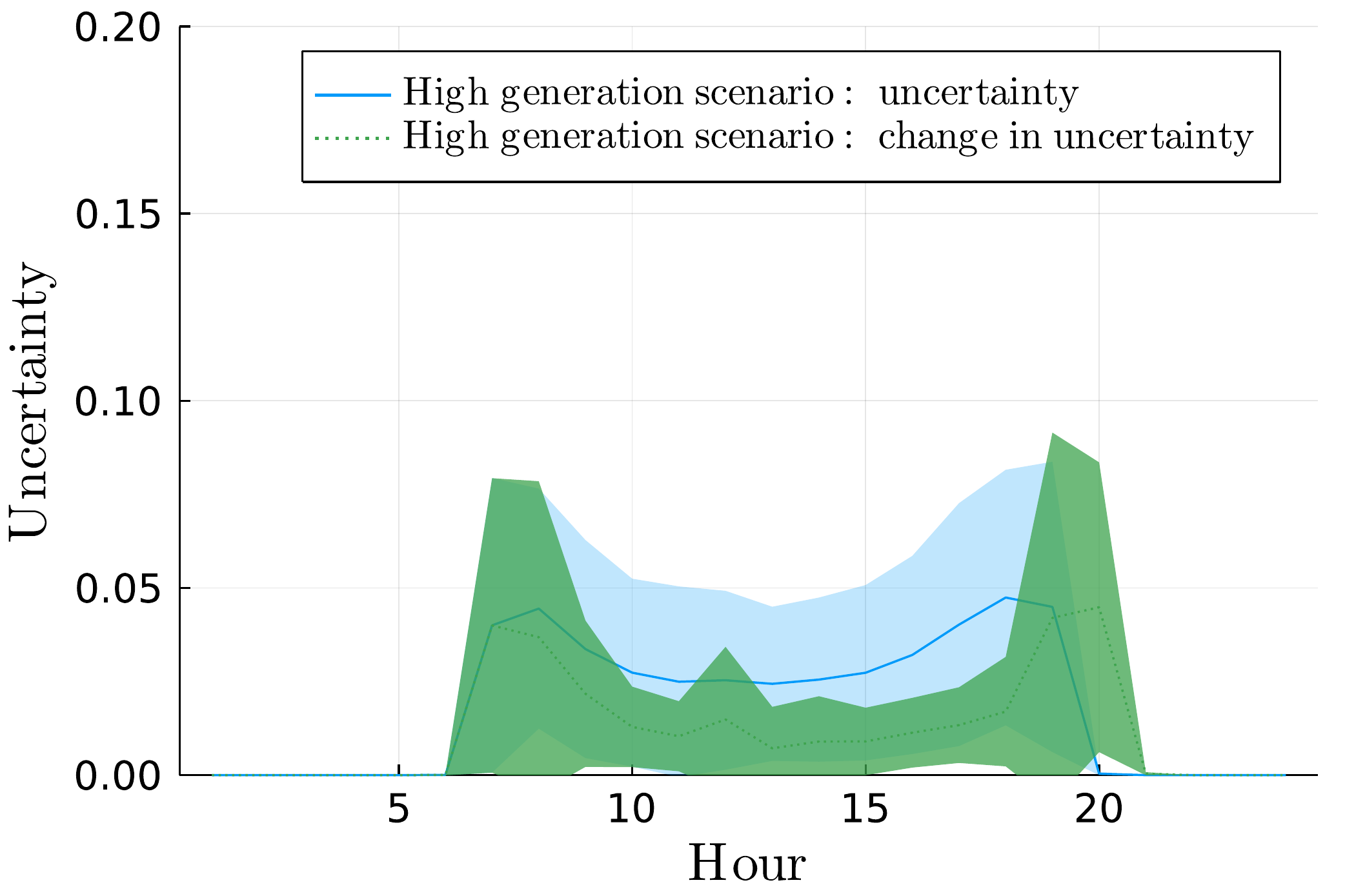}
\end{subfigure}
\vspace{10pt}
\caption{The empirical distribution of the uncertainty per hour and of the daily sum are both centered and concentrated around zero (left). Moreover, mean absolute uncertainty/change in uncertainty during consecutive hours is small and is maximized during sunrise/sunset (right). All in all, mean uncertainty/change in uncertainty is only a small fraction of mean solar generation.}
\label{fig:unc3}
\end{figure}

Figure \ref{fig:unc3} investigates the empirical distribution of deviations between solar capacity factors and their forecasts, and the change in uncertainty between consecutive hours. 
Figure \ref{fig:unc3} (left) investigates the empirical distribution of uncertainty, $u_n^{h,d,m,y}$: across all reduced scenarios, the distribution is centered and 
concentrated around zero. 
Figure \ref{fig:unc3} (right) shows the mean and standard deviation of the absolute uncertainty, $|u_n^{h,d,m,y}|$ and the absolute change in uncertainty during consecutive hours, $|u_n^{h,d,m,y}-u_n^{h-1,d,m,y}|$. We observe that both quantities are small compared to the nominal solar generation values, and peak during sunrise and sunset. Further, the absolute change in uncertainty is, on average, smaller than the absolute uncertainty. This observation suggests that, if a realized day is above the mean solar generation of its corresponding reduced scenario at hour $h-1$, it will likely also be at hour $h$.

We now summarize the main managerial insights from our analysis and their implications:
\begin{itemize}
    \item Figures \ref{fig:unc1} suggests that uncertain parameters during different hours of the day, as well as during different reduced scenarios, should be modeled separately in order to obtain as accurate a model as possible. However, provided that the number of reduced scenarios is sufficiently large, we need not treat uncertain parameters differently depending on which month they correspond to.
    \item Figure \ref{fig:unc3} (left) suggests uncertainty sets motivated by the {\msom central limit theorem (CLT)} may perform well in practice.
    \item Figure \ref{fig:unc3} (right) implies that ``smoothness constraints'' which link uncertainty during consecutive hours and require that it varies slowly, may perform well in practice.
\end{itemize}
Motivated by these findings, we propose a data-driven uncertainty set in the next subsection.

\subsection{A Data-Driven Uncertainty Set} \label{ssec:RO-uncertainty-set}
We now propose a data-driven uncertainty set which immunizes the sample-robust model proposed in the previous section against intermittency in solar generation capacity. This immunization is significant, since the SAA model we derived in the previous section relies on day-ahead solar capacity forecasts being deterministically accurate, which does not occur in practice.

As discussed in Section \ref{ssec:RO-uncertainty-defn}, the (uncertain) solar generation at node $n$ and time $(h,d,m,y)$ is taken to be $v^{h,d,m,y}_n = \bar{v}^{h,d}+u^{h,d,m,y}_n$, where $u^{h,d,m,y}_n$ is a component of the vector of uncertain parameters $\boldsymbol u \in \mathcal{U}_u$. Moreover, as discussed in Section \ref{ssec:RO-uncertainty-analysis}, the practical behaviour of $\mathcal{U}_u$ can be characterized by three constraints. First, the absolute value of $u^{h,d,m,y}_n$ is never too large in any hour. Second, {\msom $u^{h,d,m,y}_n$} never changes by too much between an hour and the previous one. Third, {\msom $u^{h,d,m,y}_n$} obeys the CLT in aggregate. Combining the ideas from Sections \ref{ssec:RO-uncertainty-defn} and \ref{ssec:RO-uncertainty-analysis} therefore gives rise to the following uncertainty set:
\begin{equation} \label{eqn:uncertainty-set}
\begin{split}
    \mathcal{U}_u = \left\{ \right. \boldsymbol{u} \in \mathbb{R}^{|\mathcal{H}| \times |\mathcal{D}| \times |\mathcal{M}| \times |\mathcal{Y}| \times |\mathcal{N}|}: \quad  
    & |u_n^{h,d,m,y}| \leq U_{\text{MAX}}^{h,d}, \ \left| u^{h,d,m,y}_n - u^{h-1,d,m,y}_n \right|  \leq U_{\text{SV}}^{h,d} \quad \forall h, d,m,y,n,\\
    &  \left. \left| \sum_{n,h,m,y} u^{h,d,m,y}_n \right| \leq U_{\text{CLT}}^{d} \quad \forall d\right\}.
\end{split}
\end{equation}
This set is parameterized by the following three sets of parameters, which we now define and provide guidance on how to set their values:
\begin{itemize}
    \item $U_{\text{MAX}}^{h,d} \in \mathbb{R}_+$: The maximum allowable perturbation for hour $h$ of reduced scenario $d$. We estimate $U_{\text{MAX}}^{h,d}$ as a fraction $\gamma_{\max} \geq 0$ of the maximum absolute deviation from the cluster centroid across all data points that are assigned to the same cluster: $ U_{\text{MAX}}^{h,d} = \gamma_{\max} \cdot \max_{i: c(i)=d} \left| v_0^{h,i}-\bar{v}^{h,d} \right| .$
    
    \item $U_{\text{SV}}^{h,d} \in \mathbb{R}_+$: Imposes smoothness across time, i.e., the uncertainty in solar generation cannot be too different between hours $h-1$ and $h$ of reduced scenario $d$. We estimate $U_{\text{SV}}^{h,d}$ as a fraction $\gamma_{\text{C}} \geq 0$ of the maximum absolute difference in uncertainty across all data points that are assigned to the same cluster: $ U_{\text{SV}}^{h,d} = \gamma_{\text{C}} \cdot \max_{i: c(i)=d} \left| \left( v_0^{h,i}-\bar{v}^{h,d} \right) - \left( v_0^{h-1,i}-\bar{v}^{h-1,d} \right) \right| .$
    
    \item $U_{\text{CLT}}^{d} \in \mathbb{R}_+$: Controls the limiting behavior of the realized uncertainty.  For each reduced scenario $d$, and for fixed $n,m,y$, we focus on the aggregated daily uncertainty $\sum_{h} u_n^{h,d,m,y}$. The central limit theorem implies that the sum across $n,m,y$ of the realized aggregated daily uncertainties $\sum_{n,m,y} \left(\sum_{h} u_n^{h,d,m,y}\right)$ converges to a zero-mean normal distribution with standard deviation $\sigma^d$. 
    We take $\sigma^d$ as the corrected sample standard deviation and estimate $U_{\text{CLT}}^{d}$ as a fraction $\gamma_{\text{CLT}} \geq 0$ of the sample standard deviation: $ U_{\text{CLT}}^{d} = \gamma_{\text{CLT}} \cdot \sqrt{|\mathcal{N}|\cdot|\mathcal{Y}|\cdot|\mathcal{M}|\cdot|\mathcal{H}|} \cdot \sigma^d .$

\end{itemize}

\subsection{A Constraint-Aggregation Approach} \label{ssec:RO-constraint-aggregation}
In this section, we integrate the uncertainty set derived in section \ref{ssec:RO-uncertainty-set} within our SAA model. This is non-trivial: in our SAA model the uncertain solar capacities $v^{h,d,m,y}_n$ appears in constraints \eqref{eqn:kirchhoff} and \eqref{eqn:sell} individually. Therefore, as is well documented in the RO literature \citep[see, e.g.,][]{roos2020reducing}, directly imposing robustness constraint-wise is equivalent to simultaneously setting the capacity factors $v^{h,d,m,y}_n$ to their minimum value in each hour of the planning horizon and at each node in the network, while completely ignoring the smoothness and CLT constraints. To avoid this, we adopt a constraint aggregation approach which requires that the aggregated demand across all nodes and all time periods is feasible for all realizations of the uncertain parameter, and the total amount of energy sold across all nodes and all time periods satisfies legal requirements.


To simplify our notation for the rest of this section, we now concatenate all variables except the ones corresponding to solar investment decisions into the vector: \vspace{-3mm} $$\boldsymbol x := (\boldsymbol f, \boldsymbol r, \boldsymbol x_O, \boldsymbol x_N, \boldsymbol w).$$ Further, we denote the aggregated installed solar capacity in node $n$ and in year $y$ by \vspace{-3mm}$$ \bar{z}_n^y := \sum_{y'=1}^y \xi^{y-y'} z_n^{y'}.$$ By doing so, in both constraints of interest, the coefficients of $\boldsymbol x$ are deterministic, whereas the coefficients of $\boldsymbol{\bar z}$ involve the uncertain parameters. Let us first focus on the demand constraints \eqref{eqn:kirchhoff}. The resulting aggregated constraint, which needs to hold for all realizations of $\boldsymbol{v} \in  \mathcal{U}${\msom , all scenarios $d \in \mathcal{D}$, and all months $m \in \mathcal{M}$} is
\begin{equation} \label{eqn:unc-constraint-1}
\begin{split}\msom 
    \sum_{n,h,y} & \msom  \quad \left[
    \sum_{a \in \mathcal{I}(n)} \tau_a(f_{a}^{h,d,m,y})
    + \sum_{a \in \mathcal{O}(n)} \tau_a(-f_{a}^{h,d,m,y})
    + v^{h,d} \left( \sum_{y'=1}^y \xi^{y-y'} z_n^{y'} \right) \right. \\
    & \msom  \quad \left. + R\cdot r_{n}^{h,d,m,y} + x_{O,n}^{h,d,m,y} + x_{N,n}^{h,d,m,y} - w_{n}^{h,d,m,y} \right] \ \geq \sum_{n,h,y} d_{n}^{h,m,y} , \qquad \forall \boldsymbol{v} \in \mathcal{U}.  
\end{split}
\end{equation}
Moreover, it is not too hard to see that this condition holds if and only if $$\msom \boldsymbol{a_1}^\top \boldsymbol{x}
    + \min_{\boldsymbol{v} \in \mathcal{U}} \sum_{n,h,y} v^{h,y}_n \bar{z}_n^y \geq b_1, $$ where we denote by $\boldsymbol{a_1}$ all deterministic coefficients in Equation \eqref{eqn:kirchhoff} (i.e., all coefficients except $ \boldsymbol v$) and by $b_1 := \sum_{n,h,y} d_{n}^{h,m, y}$. We now derive this constraint's deterministic equivalent via the following lemma, which holds due to strong duality (proof deferred to {\msom Section} \ref{appx:proofs}):
{\msom     
\begin{lemma}\label{lemma:detequiv}
The robust constraint \eqref{eqn:unc-constraint-1} admits the following deterministic reformulation:
\begin{equation} \label{eqn:rc-1}
\begin{split}
    & \boldsymbol{a_1}^\top \boldsymbol{x} + \left( \sum_{n,h,y} \bar{v}^{h,d} \sigma_{1,n}^{h,d,m,y} + U_{\text{MAX}}^{h,d} \phi_{1,n}^{h,d,m,y} + U_{\text{MAX}}^{h,d} \chi_{1,n}^{h,d,m,y} + U_{\text{SV}}^{h,d} \psi_{1,n}^{h,d,m,y} + U_{\text{SV}}^{h,d} \omega_{1,n}^{h,d,m,y} \right)  + U_{\text{CLT}}^{d} \tau_1^d \geq b_1, \\
    & \sigma_{1,n}^{h,d,m,y} \leq \bar{z}_n^y , \quad \forall h, d,m,y,n, \\
    & -\sigma_{1,n}^{h,d,m,y} - \phi_{1,n}^{h,d,m,y} +\chi_{1,n}^{h,d,m,y} - \psi_{1,n}^{h,d,m,y} + \omega_{1,n}^{h,d,m,y} +\psi_{1,n}^{h+1,d,m,y} - \omega_{1,n}^{h+1,d,m,y} + \tau_1^d = 0, \quad \forall h, d,m,y,n, \\
    & \phi_{1,n}^{h,d,m,y} \leq 0,\ \chi_{1,n}^{h,d,m,y} \leq 0,\ \psi_{1,n}^{h,d,m,y} \leq 0,\ \omega_{1,n}^{h,d,m,y} \leq 0,\ \tau_1^d \leq 0 , \quad \forall h, d,m,y,n.
\end{split}
\end{equation}
\end{lemma}

We now focus on constraints \eqref{eqn:sell}. In this case, the resulting aggregated constraint, which needs to hold for all realizations of $\boldsymbol{v} \in  \mathcal{U}$, can be written as
\begin{equation} \label{eqn:unc-constraint-2}
\begin{split}\msom 
    \sum_{n,y} & \msom  \quad \left[
    \sum_{h} w_{n}^{h,d,m,y} - \beta v_{n}^{h,d,m,y} \bar{z}_n^{y} \right]  \leq \beta  \sum_{n,h,y} \max\left\{ 0, -d_n^{h,d,m,y} \right\} , \qquad \forall \boldsymbol{v} \in \mathcal{U}, \\
    \Leftrightarrow & \msom  \quad
    \boldsymbol{a_2}^\top \boldsymbol{x}
    + \max_{\boldsymbol{v} \in \mathcal{U}} \sum_{n,h,y} \beta  v^{h,d,m,y}_n \bar{z}_n^y \leq b_2,
\end{split}
\end{equation}
where we denote by $\boldsymbol{a_2}$ all deterministic coefficients in Equation \eqref{eqn:sell} (i.e., all coefficients except the ones that involve $ \boldsymbol v$) and by $b_2 := \beta  \sum_{n,h,y} \max\left\{ 0, -d_n^{h,m,y} \right\}$. Applying the same machinery as we did in Lemma \ref{lemma:detequiv} and denoting by $\boldsymbol{\sigma_2}, \boldsymbol{\tau_2}, \boldsymbol{\phi_2}, \boldsymbol{\chi_2}, \boldsymbol{\psi_2}, \boldsymbol{\omega_2}$ the corresponding dual variables, yields the deterministic reformulation:
\begin{equation} \label{eqn:rc-2}
\begin{split}
    & \boldsymbol{a_2}^\top \boldsymbol{x} + \left( \sum_{n,h,y} \bar{v}^{h,d} \sigma_{2,n}^{h,d,m,y} + U_{\text{MAX}}^{h,d} \phi_{2,n}^{h,d,m,y} + U_{\text{MAX}}^{h,d} \chi_{2,n}^{h,d,m,y} + U_{\text{SV}}^{h,d} \psi_{2,n}^{h,d,m,y} + U_{\text{SV}}^{h,d} \omega_{2,n}^{h,d,m,y} \right) + U_{\text{CLT}}^{d} \tau_2^d \leq b_2, \\
    & \sigma_{2,n}^{h,d,m,y} \leq \beta  \bar{z}_n^y, \quad \forall h, d,m,y,n, \\
    & -\sigma_{2,n}^{h,d,m,y} - \phi_{2,n}^{h,d,m,y} +\chi_{2,n}^{h,d,m,y} - \psi_{2,n}^{h,d,m,y} + \omega_{2,n}^{h,d,m,y} +\psi_{2,n}^{h+1,d,m,y} - \omega_{2,n}^{h+1,d,m,y} + \tau_2^d = 0, \quad \forall h, d,m,y,n, \\
    & \phi_{2,n}^{h,d,m,y} \leq 0,\ \chi_{2,n}^{h,d,m,y} \leq 0,\ \psi_{2,n}^{h,d,m,y} \leq 0,\ \omega_{2,n}^{h,d,m,y} \leq 0,\ \tau_2^d \leq 0 , \quad \forall h, d,m,y,n.
\end{split}
\end{equation}

Finally, we obtain our RO model by imposing the deterministic equivalent constraints \eqref{eqn:rc-1} and \eqref{eqn:rc-2} in the SAA model, and dropping their robust counterparts \eqref{eqn:kirchhoff}-\eqref{eqn:sell}. The resulting problem is notably a linear problem with $2 \cdot |\mathcal{D}| \cdot \left( 5 \cdot |\mathcal{N}| \cdot |\mathcal{H}| \cdot |\mathcal{M}| \cdot |\mathcal{Y}| + 1 \right)$ new variables and $2 \cdot |\mathcal{N}| \cdot |\mathcal{H}| \cdot |\mathcal{D}| \cdot |\mathcal{M}| \cdot |\mathcal{Y}| + 2 |\mathcal{D}| \cdot |\mathcal{M}|  $ new constraints, and can thus be solved at scale.
}

\subsection{Preventing Overfitting: Distributionally Robust Optimization to the Rescue}\label{ssec:RO-distributionally}
In this section, we immunize our model against overfitting the solar capacity factors provided by OCP, using ideas from DRO. This is an important practical step. Indeed, as mentioned in the introduction, sample average approximations of capacity expansion problems generally perform well in large sample settings, but overfit in the presence of small sample sizes. Moreover, Moroccan weather patterns may change due to changes in the climate induced by increasing levels of carbon in the atmosphere, making historical data unreliable. Therefore, overfitting could certainly occur in our problem setting, where we have access to capacity factors on an hourly basis for one year.

{\latest To prevent this, we take a DRO approach to SAA inspired by the works of \cite{van2021data, anderson2022improving} among others.} Specifically, for ease of notation, let us vectorize the strategic decisions (investment in batteries and solar panels across all sites and years) as $\bm z_s := (\bm b, \bm z)$ and the operational decisions (cost to rent lines, cost to procure and sell energy) as $\bm x := (\bm f, \bm x_O, \bm x_N, \bm w)$ that appear in the objective. In addition, we denote by $\bm c_z$ and $\bm c_x$ the corresponding vectors of objective coefficients for the investment decisions (e.g., time-discounted cost of purchasing batteries) and operational decisions (e.g., time-discounted cost of renting lines), respectively. We then rewrite our model as:\vspace{-2mm}
\begin{align}
    \min_{(\bm z_s, \bm x) \in \mathcal{Z}} & \quad \underbrace{\bm c_z^\top \bm z_s}_{\text{strategic decisions}}
        +\underbrace{\sum_{y,m,d} {\msom \text{D}^{m,y} \text{P}^{d,m,y}} (\bm c_x^{m,y})^\top \bm{x}^{d,m,y}
    }_{\text{operational decisions}},
\end{align}
where $\mathcal{Z}$ denotes the feasible set defined by constraints {\msom \eqref{eqn:budget}, \eqref{eqn:batteries}-\eqref{eqn:flow}, \eqref{eqn:rc-1} and \eqref{eqn:rc-2}}.

We now improve this model, by optimizing over the worst case probability measure $\bm{Q}$ within a KL-divergence of $\delta$ from the empirical measure $\bm{P}$; recall that the KL-divergence of two measures $\bm{Q}$, $\bm{P}$ defined on the same probability space is the asymmetric distance
$\text{D}_\text{KL}(\bm Q|| \bm P) :=\sum_{d \in \mathcal{D}} Q^d \log\left(\nicefrac{Q^d}{P^d}\right).$
As shown by \cite{van2021data}, this approach behaves optimally in terms of minimizing the out-of-sample disappointment over all variants of SAA. Moreover, as shown by \cite{anderson2022improving}, this approach outperforms SAA in the small sample size regime. Formally, we have the DRO problem:\vspace{-3mm}
\begin{align}\label{prob:dro}
    \min_{(\bm z_s, \bm x) \in \mathcal{Z}} & \quad \bm c_z^\top \bm z_s
        + \sum_{y,m} {\msom \text{D}^{m,y}} \underset{\substack{\bm Q^{m,y} \in \mathcal{Q}(\bm P^{m,y},\delta) }}{\max} \sum_d Q^{d,m,y} (\bm c_x^{m,y})^\top \bm{x}^{d,m,y}
\end{align}
where $\mathcal{Q}(\bm P,\delta):=\{\bm{Q} \in \mathbb{R}^{|\mathcal{D}|}_+: \bm{e}^\top \bm{Q}=1, \sum_{d \in \mathcal{D}} Q^d \log(Q^d/P^d) \leq \delta\}$ denotes the set of all probability measures within a KL-divergence of $\delta$ of the empirical one. We now reformulate the DRO problem as a deterministic one, via the following lemma which leverages strong conic duality and is essentially due to \citet[Theorem 1]{kocuk2020conic}:
\begin{lemma}\label{lemma:dro}
Problem \eqref{prob:dro} admits the following deterministic reformulation:\vspace{-2mm}
\begin{align}\label{prob:dro_det}
    \min_{\substack{
    (\bm z_s, \bm x) \in \mathcal{Z}, \\
    \bm \alpha \in \mathbb{R}^{|\mathcal{M}|\cdot|\mathcal{Y}|}, \\
    \bm \beta \in \mathbb{R}^{|\mathcal{M}|\cdot|\mathcal{Y}|} \geq 0, \\
    \bm{\gamma}, \bm{\zeta} \in \mathbb{R}^{|\mathcal{D}|\cdot|\mathcal{M}|\cdot|\mathcal{Y}|}}} 
    & \quad \bm c_z^\top \bm z_s 
    + \sum_{y,m} \left( {\msom \text{D}^{m,y}} \alpha^{m,y} + \delta \beta^{m,y} + \sum_d \text{P}^{d,m,y} \gamma^{d,m,y} \right) \nonumber\\
    \text{\rm s.t.} & \quad \alpha^{m,y} - \zeta^{d,m,y} \geq (\bm c_x^{m,y})^\top \bm{x}^{d,m,y} , \ (-\beta^{d,m}, \zeta^{d,m,y}, \gamma^{d,m,y}) \in \mathcal{K}_{\text{exp}}^\star,
\end{align}
 $$\text{\rm where:} \quad \mathcal{K}_{\text{exp}}^\star:=\mathrm{cl}\left(\left\{(u,v,w) \in \mathbb{R}^3: -u \exp(v/u) \leq \exp(1)w, u<0\right\}\right)$$ denotes the dual cone to the exponential cone and $\mathrm{cl}$ denotes the closure of a set \citep[see][for a derivation of the cone]{serrano2015algorithms}. 
\end{lemma}
We remark that we could replace the KL divergence with another $\phi$-divergence or distance measure without significantly impacting the tractability of the deterministic reformulation; see \citet{rahimian2019distributionally} for a review of other DRO formulations. Lemma \ref{lemma:dro} reveals that we can convert our SAA problem into a DRO one while only increasing the number of decision variables by $2\cdot|\mathcal{M}|\cdot|\mathcal{Y}|+2\cdot|\mathcal{D}|\cdot|\mathcal{M}|\cdot|\mathcal{Y}|$ and the number of (dual) exponential cone constraints by $|\mathcal{D}|\cdot|\mathcal{M}|\cdot|\mathcal{Y}|$. Thus, our DRO formulation can be tractably solved by state-of-the-art conic solvers such as \verb|Mosek|.

\section{Numerical Results} \label{sec:exper}
In this section, we describe how the deterministic and robust optimization methodologies proposed in Sections \ref{sec:SAA}—\ref{sec:RO} {\msom can be} implemented in practice. We first explore the relationship between OCP’s investment level and long-run operational costs and carbon emissions (Section \ref{ssec:exper-invest-vs-oper}), and establish that investing {\latest $10$ billion MAD (resp. $20$ billion MAD)} reduces OCP’s carbon emissions by {\latest $70\%$ (resp. $95\%$)} in a profitable fashion. Next, we study the investment policy prescribed by the model with a {\latest $20$ billion MAD budget in detail} (Section \ref{ssec:exper-inspect}), and demonstrate that the model both anticipates solar generation before it occurs and captures notions of load shifting. {\latest Finally, we summarize our findings in Section \ref{ssec:exper-summary} and discuss their impact on OCP's operations in Section \ref{ssec:exper-impact}.}

{\msom  From a managerial perspective, the experiments described in this section are very similar to those conducted by OCP, using our model, to size their final investment in solar panels and batteries. Indeed, they provide important managerial insights into the number of solar panels and batteries that should be installed at each year of the planning horizon and at each site under a given investment budget, and how OCP's system operates on a day-to-day basis.

}
{\msom \paragraph{Problem Data:}
All experiments were run using historical solar capacity data collected by OCP's engineers, energy demand data calibrated according to OCP's sales forecasts, and manufacturer data on estimated battery and solar panel cost/efficiency. All other data, including the problem data described in Table \ref{tab:multiperioddata}, were set in collaboration with OCP's engineers, by iteratively running the models described in the previous two sections, jointly examining their output, discussing whether any aspects of the output did not make sense to OCP's engineers, and calibrating values such as the discount factor accordingly. To preserve OCP's privacy, both this problem data and our code are withheld from our paper and its online supplement.

\paragraph{Hyperparameter Tuning and Impact of RO:} {\latest Using a standard machine learning paradigm, we set all RO and DRO hyperparameters by splitting our data into a training set, a validation set, and a testing set and selecting the hyperparameters that performed best (in terms of minimizing both the cost of investment plus operations and the CO$_2$ emissions) on the validation set under an investment budget of $20$ billion MAD. In particular, combining insights from our cross-validation procedure (which we detail in Section \ref{appx:exper-tuning}) with discussions with the OCP team, we set the RO and DRO hyperparameters to $(\gamma_{\max},\gamma_C,\gamma_{\text{CLT}},\delta)=(0.5,0.5,0.5,0.01)$. As shown in Figure \ref{fig:robust_validation_heatmat}, under an investment budget of $20$ billion MAD, using RO and DRO improves the performance of our models on the validation set by around $16\%$ compared to a pure SAA model, which emphasizes the benefits of robustness in a numerically striking manner.
}

After tuning the hyperparameters used in the model, we also performed several robustness checks on the model in Section \ref{appx:exper-sensitivity}. Namely, we investigate the stability of OCP's operational cost, the total estimated carbon emissions reduction, and the amount invested in solar panels and batteries as we vary the number of scenarios. By pairing scenarios together, we also considered partially relaxing the constraint that load shifting cannot occur between days. We observe empirically that our optimal solutions are stable in both cases, which confirms the stability of our approach.
}


\subsection{Trading Off Investment Costs Against Operational Costs and CO$_2$ Emissions} \label{ssec:exper-invest-vs-oper}
In this section, we investigate the potential benefits of our cross-validated model for OCP in terms of both cost savings and emission reductions. Specifically, we explore the relationship between OCP's investment budget in the cross-validated model, and its long-run costs and carbon emissions.

\begin{figure}[!ht] 
\begin{adjustbox}{minipage=\linewidth,scale=1}
    \centering
    \begin{subfigure}[b]{0.48\linewidth}
    
        \centering
        \includegraphics[width=\linewidth]{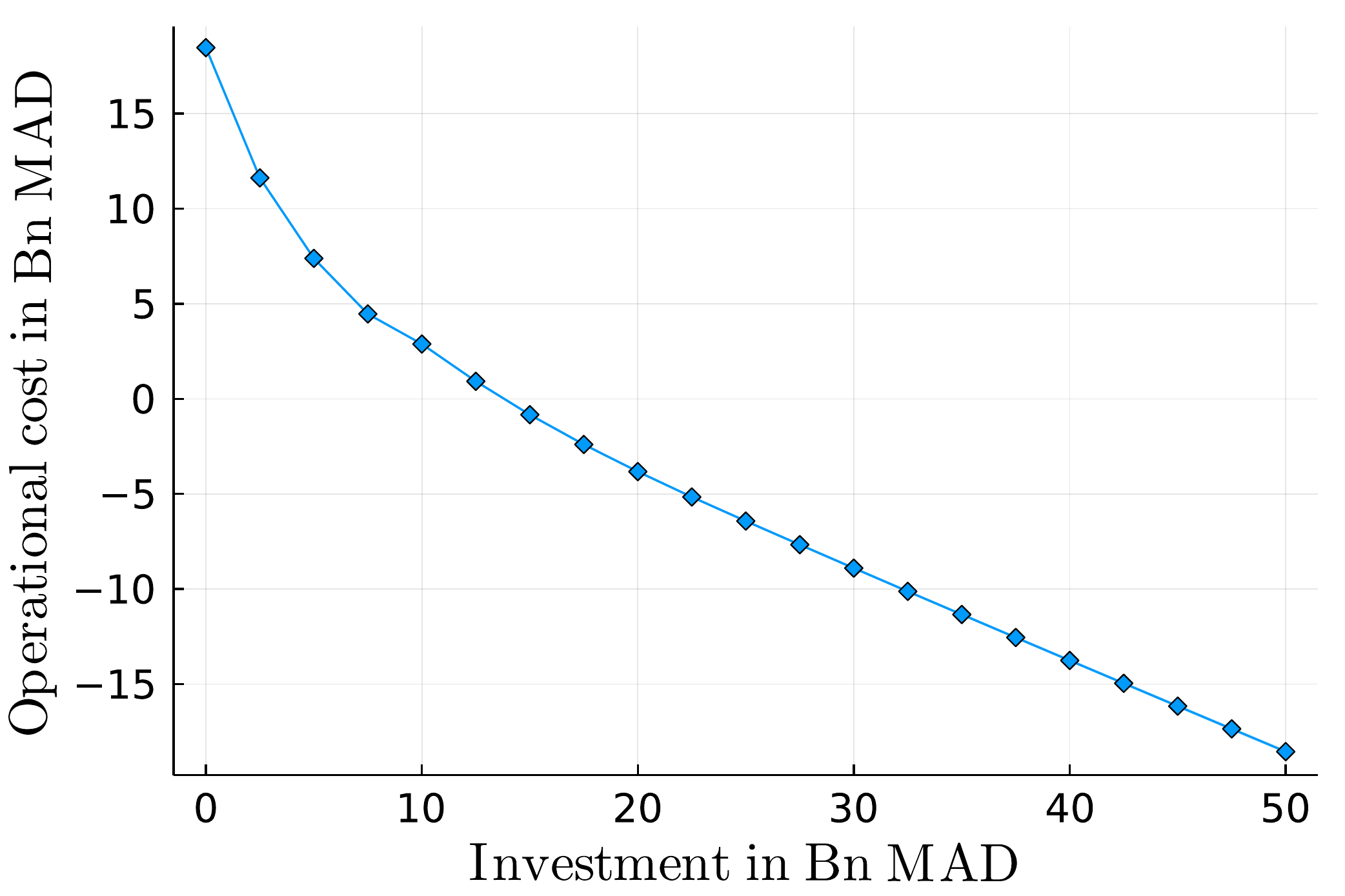}
        \label{fig:budget_oper_holdout}
    
    \end{subfigure}
    \hfill
    \begin{subfigure}[b]{0.48\linewidth}
    
        \centering
        \includegraphics[width=\linewidth]{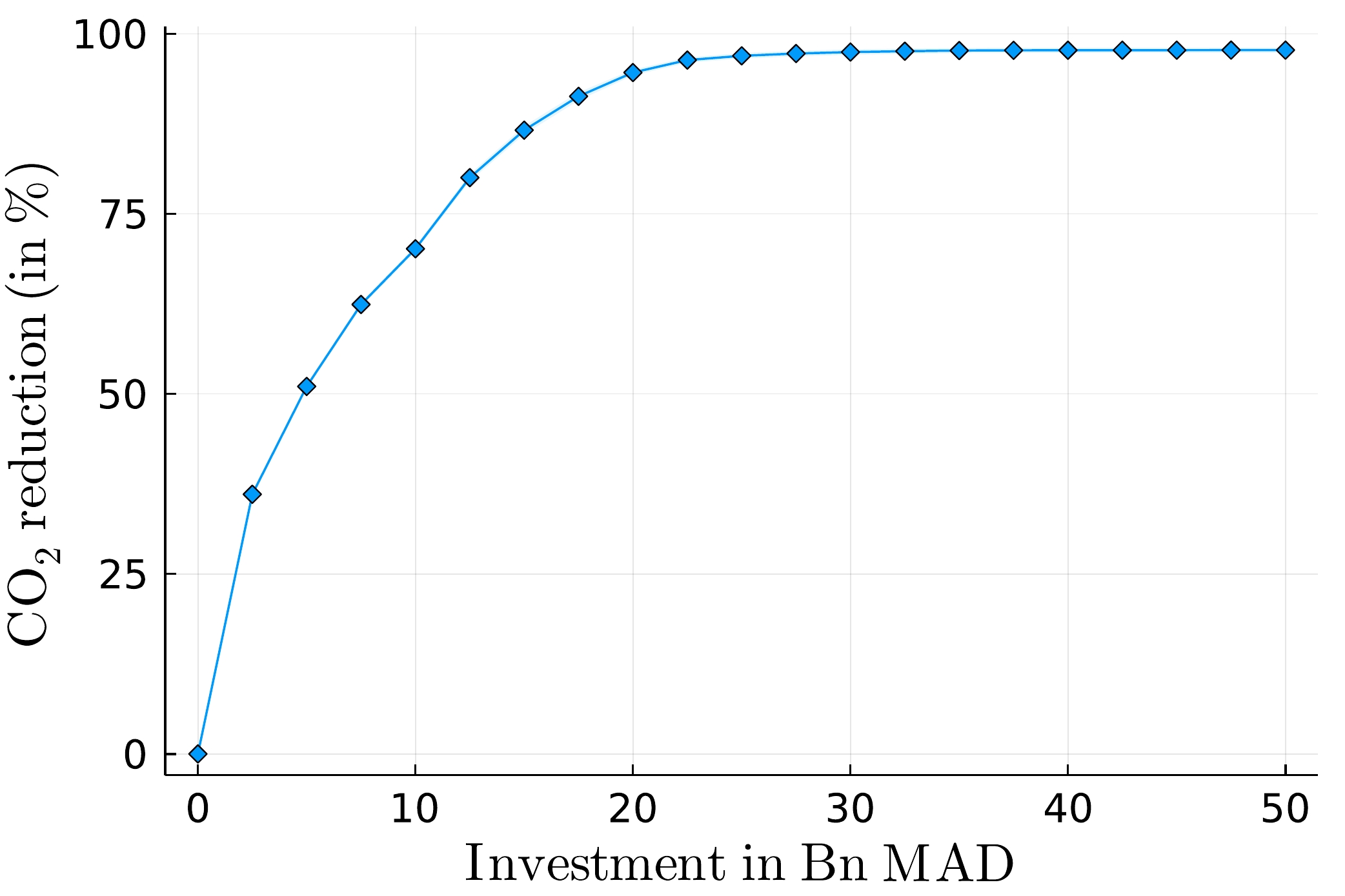}
        \label{fig:budget_co2_holdout}
    
    \end{subfigure}
\end{adjustbox}
\caption{Trade-off between overall investment cost and operational costs (left), CO$_2$ emissions (right). Increasing investment decreases OCP's cost of operations (left) and OCP's carbon emissions (right).}
\label{fig:tradeooff}
\end{figure}

We solve our cross-validated model with an overall investment budget of {\msom $B$ for each $B \in \{0,2.5, 5, \ldots 50\}$ billion MAD}. For each budget $B$, we compute OCP's anticipated operational costs (i.e., the cost of procuring energy from ONEE/NAREVA, plus cost of renting lines, minus profit from energy sold) and expected CO$_2$ emissions reduction over the $20$-year time horizon with $10$ reduced scenarios. We compare our cross-validated model against the no-investment baseline which satisfies OCP's energy needs via purchasing electricity from the grid locally at each site, in terms of both improvement in operational costs and emissions reduction, in Figure \ref{fig:tradeooff}. {\msom Further, we depict the net present value of the project under different investment budgets in Section \ref{appx:exper-NPV} (Figure \ref{fig:tradeooff-npv}). The baseline's operational costs on the testing set is $18$ billion MAD. }

{\latest As reflected\footnote{\msom Note that we reported a different (higher) return on investment for a given investment budget in an earlier version of this manuscript; this discrepancy arises because we previously used a less significant discount factor and different demand data; OCP has since updated their production plan and their expected energy demand in each year of the planning horizon.} in Figure \ref{fig:tradeooff}, increasing the investment decreases OCP's cost of operations (which, with an investment of $15$ billion MAD becomes negative owing to profit made from energy sold) and OCP's carbon emissions. The net present value of the project, excluding the salvage value of solar panels and batteries at the end of the planning horizon, increases with the value of the investment budget until it peaks at around $7.5$ billion MAD, and subsequently decreases with the investment size (see Figure \ref{fig:tradeooff-npv} in Section \ref{appx:exper-NPV}). Therefore, there exists a trade-off between fully decarbonizing OCP and maximizing OCP's profit from the initiative: investing around $10$ billion MAD returns a substantial net present value (NPV) of $5$ billion MAD and reduces OCP's carbon emissions by $70\%$; investing $20$ billion MAD returns a lower NPV of $2$ billion MAD and reduces OCP's carbon emissions by $95\%$; investing $30$ billion MAD reduces OCP's carbon emissions by $97\%$ but returns a NPV of negative $2.5$ billion MAD. } Therefore, OCP ultimately needs to choose an investment level that balances its attitude toward risk against both its interest in lowering global CO$_2$ emissions and the profit that it could make from this initiative. 

{\msom We remark that the phenomenon that completely decarbonizing a production system is substantially more expensive than partially decarbonizing it has been observed in other contexts. Indeed, \cite{downward2020using} recently investigated strategies for fully decarbonizing New Zealand's energy market, and found that doing so was prohibitively expensive compared to reducing carbon emissions to very low levels. This suggests that either installing alternatives to solar panels or making use of carbon capture technologies might be a more efficient way to decarbonize OCP completely. To address this issue, we are investigating installing wind generation as follow-up work.}

\subsection{The Model in Action}  \label{ssec:exper-inspect} 
\begin{figure}[!ht] 
\begin{adjustbox}{minipage=\linewidth,scale=1}
    \centering
    \begin{subfigure}[b]{0.48\linewidth}
    
        \centering
        \includegraphics[width=\linewidth]{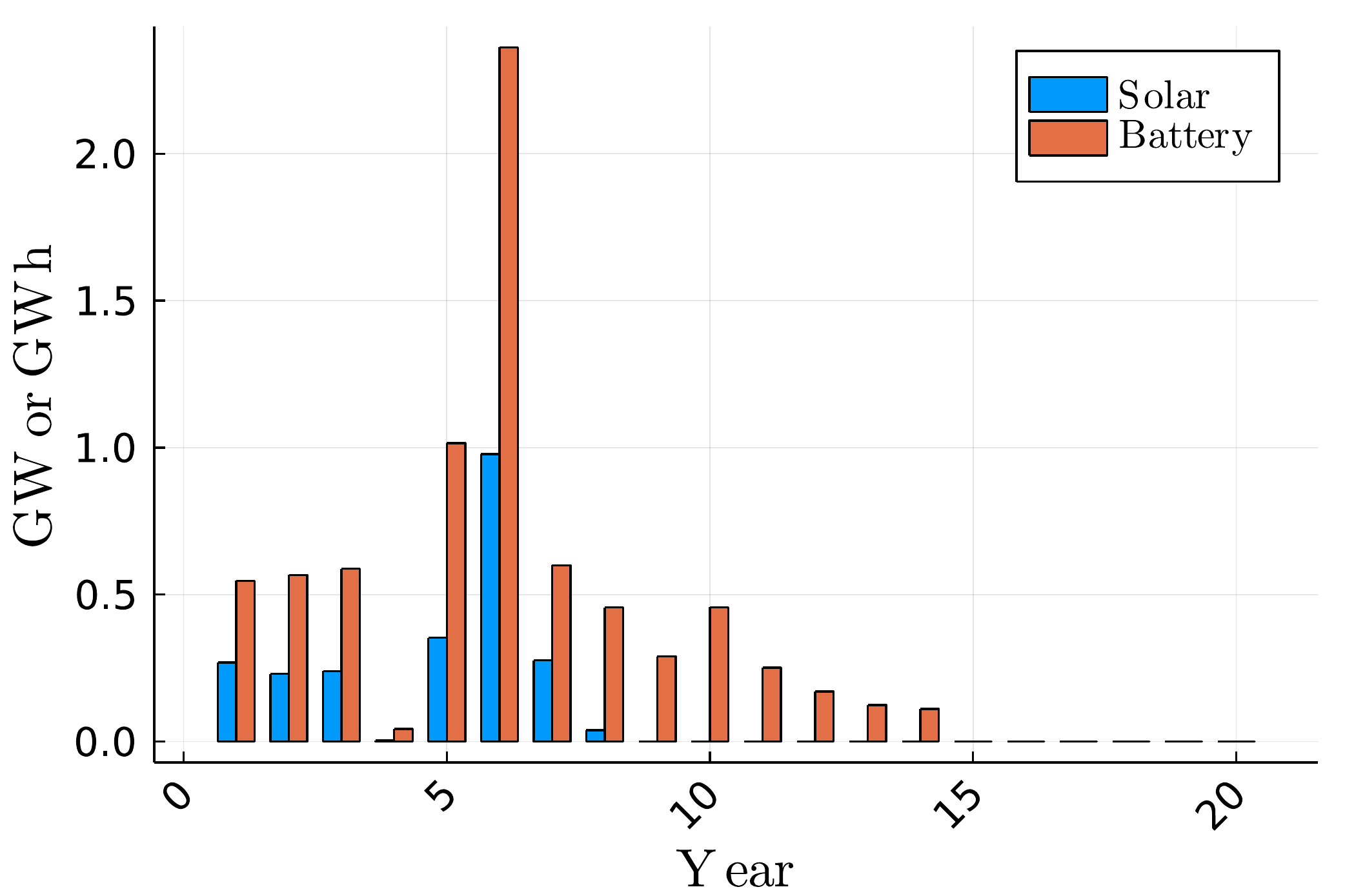}
        \label{fig:budget_investment_by_year}
    
    \end{subfigure}
    \hfill
    \begin{subfigure}[b]{0.48\linewidth}
    
        \centering
        \includegraphics[width=\linewidth]{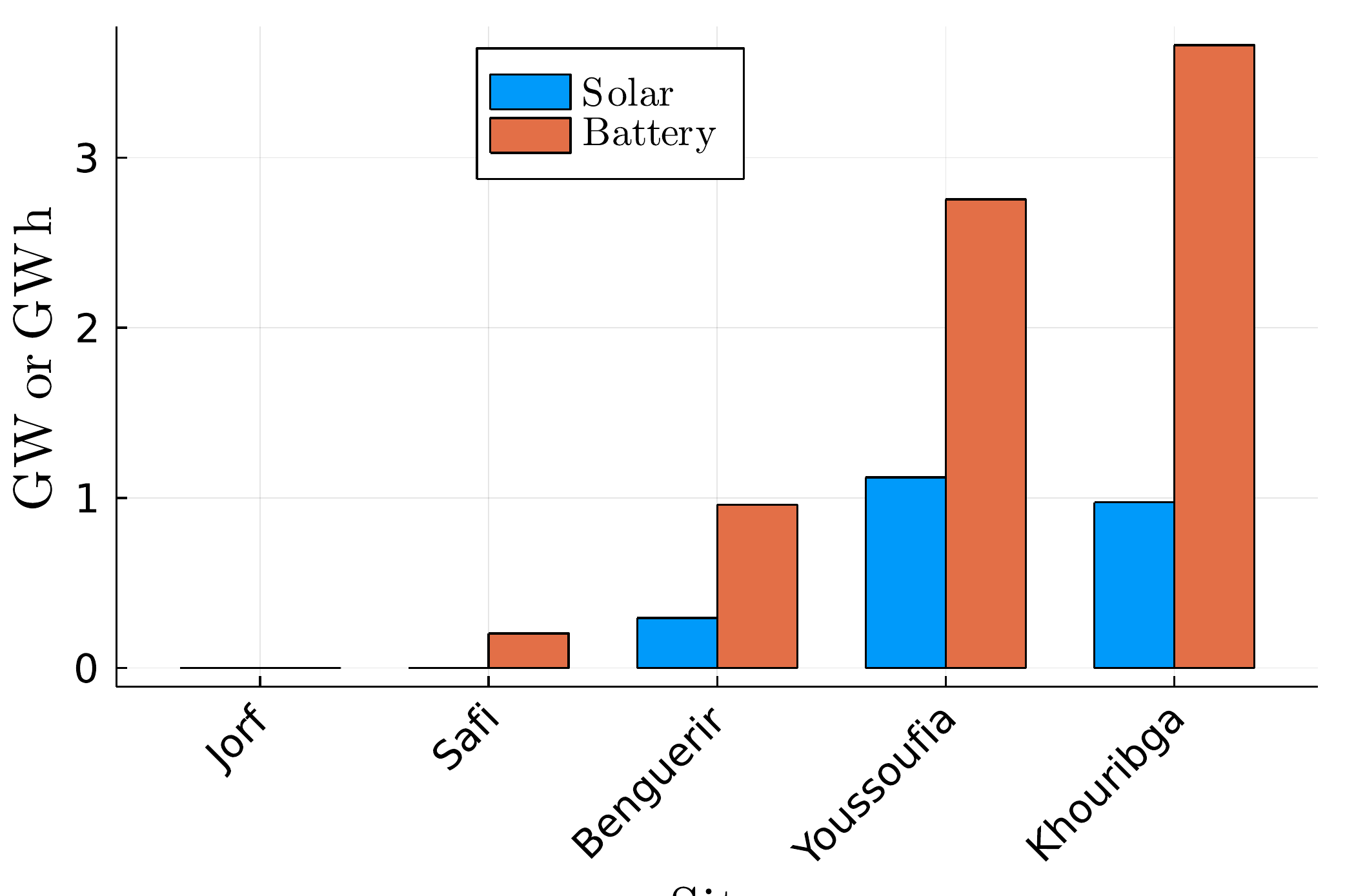}
        \label{fig:budget_investment_breakdown_by_site}
    
    \end{subfigure}
    \hfill
    \begin{subfigure}[b]{0.48\linewidth}
    
        \centering
        \includegraphics[width=\linewidth]{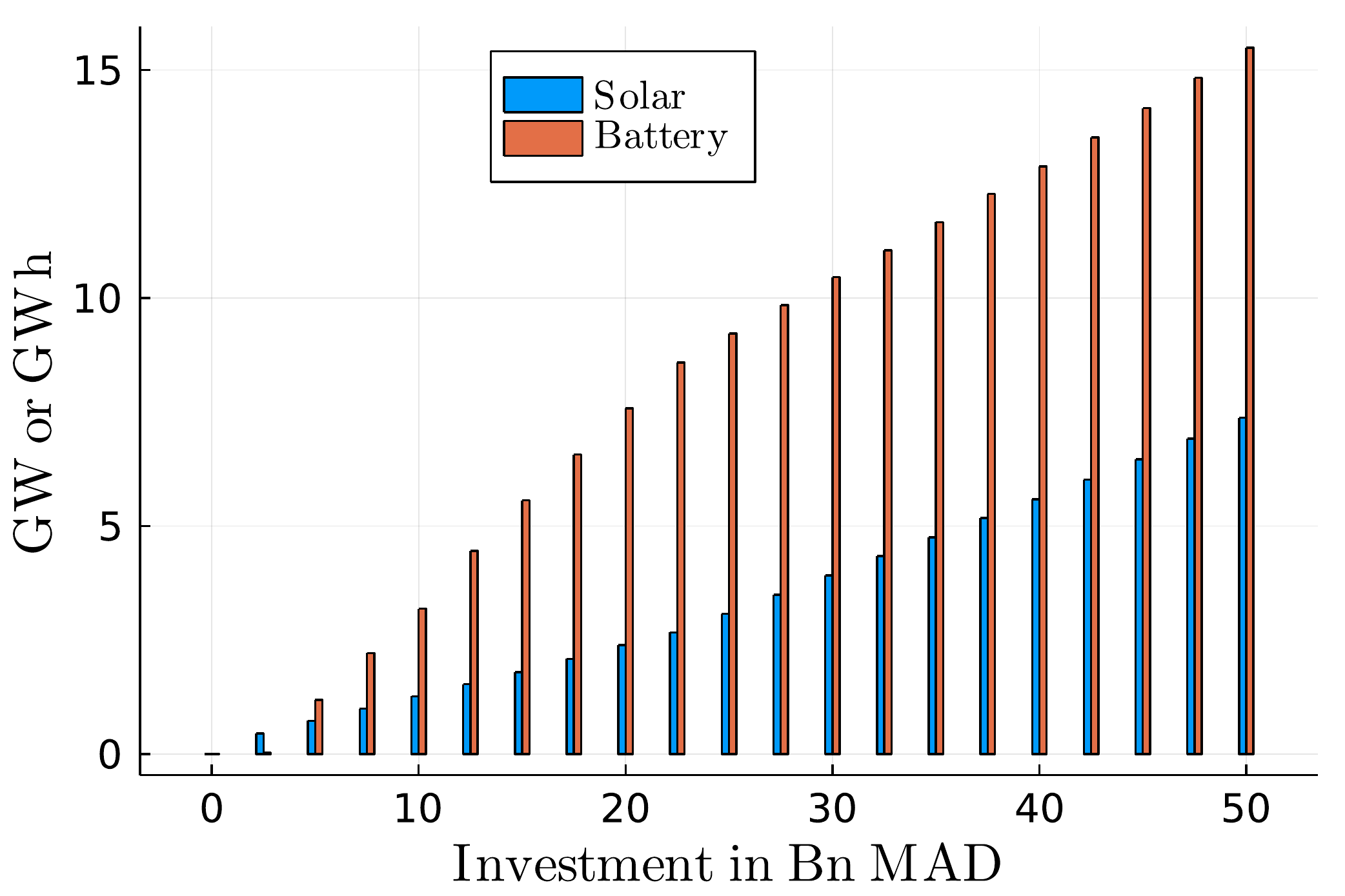}
        \label{fig:budget_investment_breakdown_by_budget}
    
    \end{subfigure}
\end{adjustbox}
\caption{Investment policy prescribed by the model over the entire horizon (aggregated across sites). We report the investment by year (left) and overall investment for increasing budget (right).}
\label{fig:budget_investment}
\end{figure}

{\msom 
In this section, we investigate the policy prescribed by our robustified and cross-validated model with a $20$-year planning horizon, $10$ reduced scenarios, and a $20$ billion MAD budget}, in order to obtain managerial insights into the structure of optimal solar capacity expansion strategies. 
{\msom 
\paragraph{Strategic decisions.} {\msom Figure \ref{fig:budget_investment} (top left) presents the total prescribed investment in solar generation and batteries by year of the planning horizon, and Figure \ref{fig:budget_investment} (top right) depicts the total investment at each site, aggregated over all years in the horizon. We observe that the model installs a significant amount of solar panels and batteries in years $1$--$5$ of the planning horizon, followed by a significant spike in year six of the planning horizon, and thereafter installs a smaller number of batteries and solar panels in the next nine years of the horizon. The spike at year six occurs because there is a cogeneration facility located at one of OCP's sites, scheduled to go offline in year six, and OCP anticipates significantly increasing its energy consumption in year six of the planning horizon. 

Figure \ref{fig:budget_investment} (bottom) depicts the prescribed investment in solar generation and batteries aggregated over the entire planning horizon for varying investment budget and reveals that the model prioritizes solar with a smaller budget, and includes more batteries as the investment size increases. This reflects two ideas. First, batteries degrade over time, so installing more as existing ones degrade can be profitable. Second, the marginal utility of installing solar panels is lower when there are already enough solar panels present in OCP's system to power OCP's operations during the day, because the energy they generate needs to be stored in batteries to power the system at night time, and this requires purchasing batteries, which can be expensive.

}

}


\begin{figure}[!ht] 
\begin{adjustbox}{minipage=\linewidth,scale=1}
    \centering
    \begin{subfigure}[b]{0.48\linewidth}
    
        \centering
        \includegraphics[width=\linewidth]{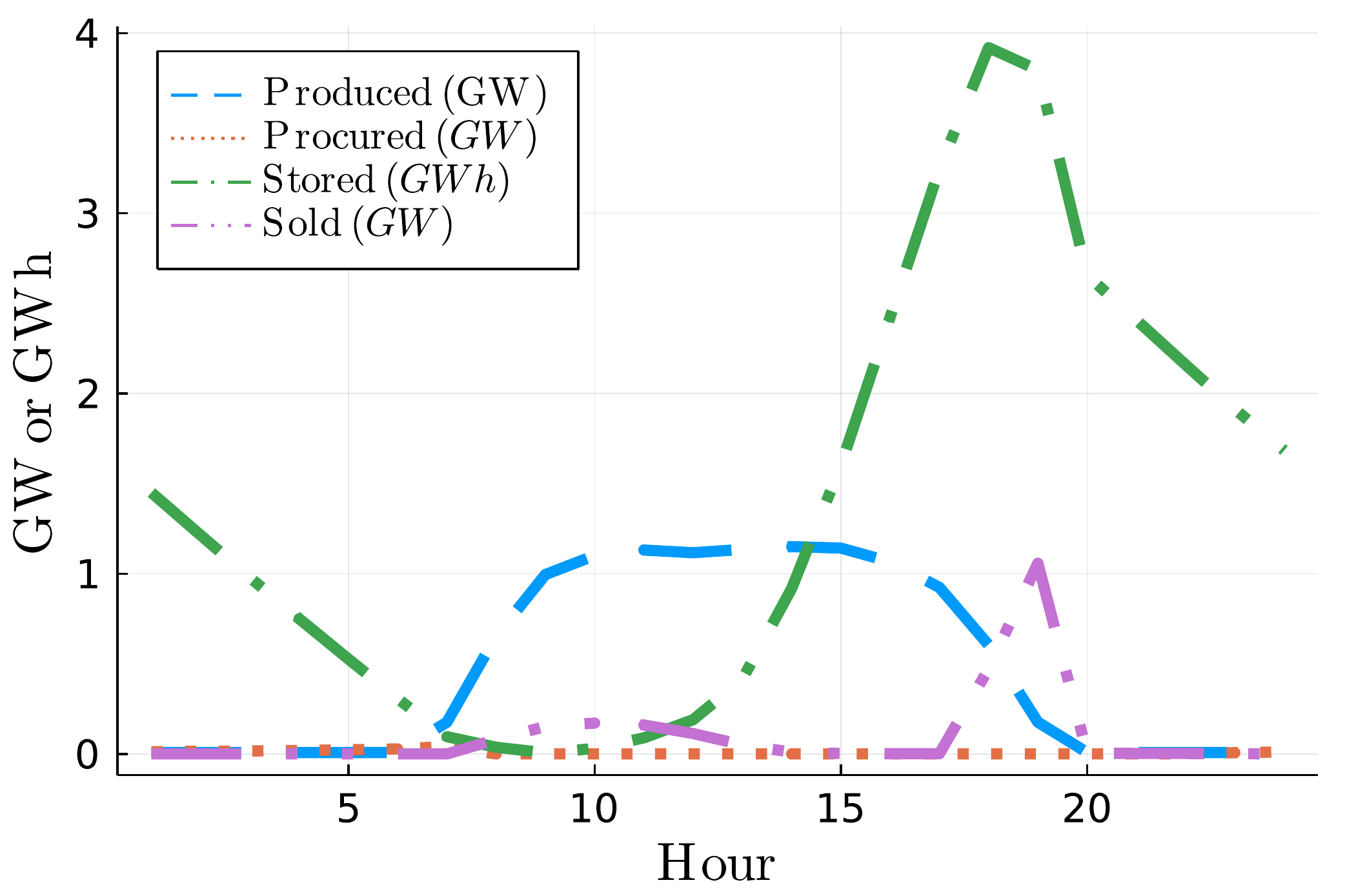}
        \label{fig:budget_energy_flow_1}
    
    \end{subfigure}
    \hfill
    \begin{subfigure}[b]{0.48\linewidth}
    
        \centering
        \includegraphics[width=\linewidth]{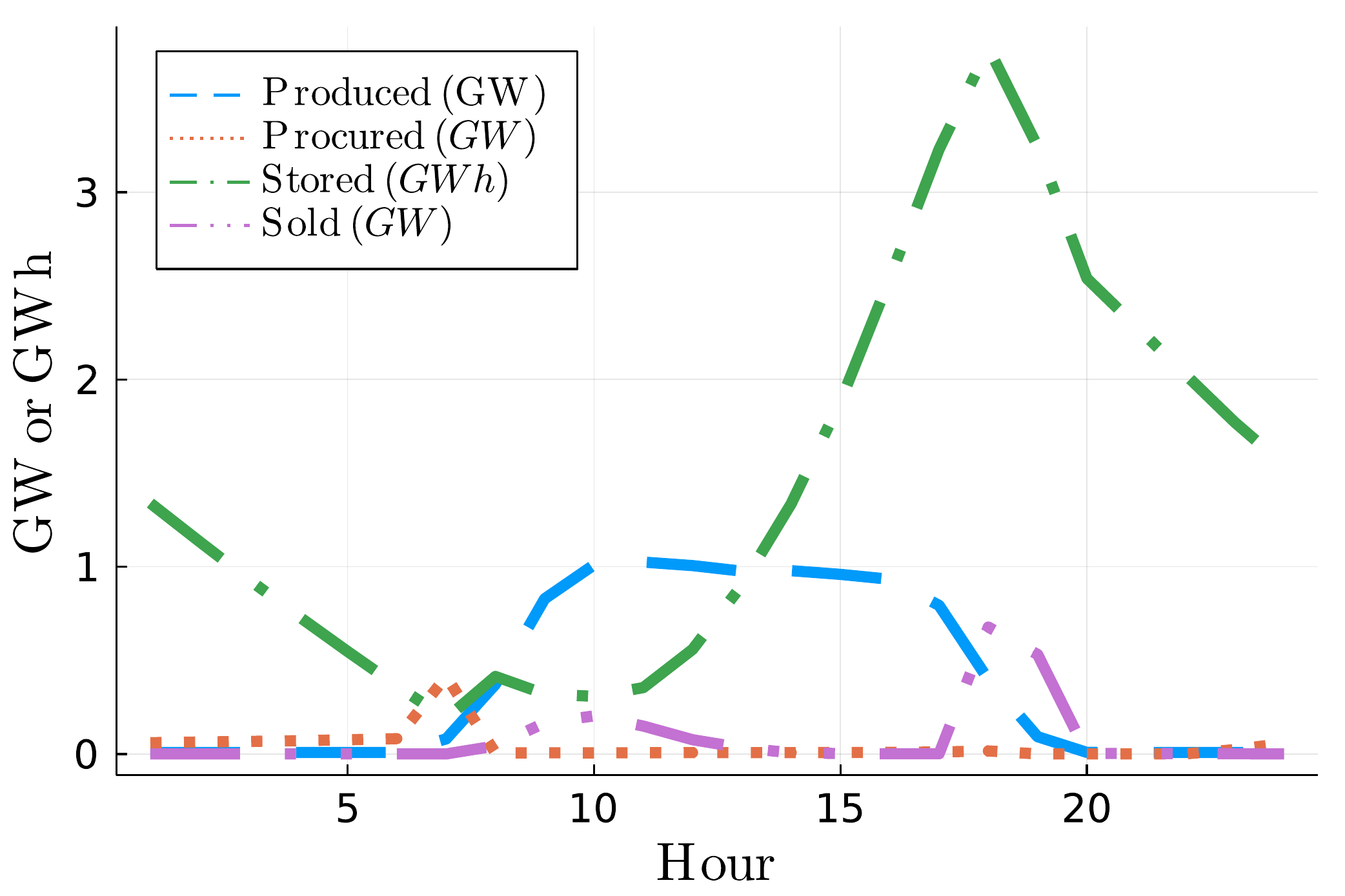}
        \label{fig:budget_energy_flow_2}
    
    \end{subfigure}
\end{adjustbox}
\caption{\msom Operational policy prescribed by the model for the two most common month-reduced scenario pairs (aggregated across sites and years). We report a high-generation scenario (left) and a low-generation scenario (right).}
\label{fig:budget_energy_flow}
\end{figure}

{\msom 
\paragraph{Operational Decisions.} {\latest Figure \ref{fig:budget_energy_flow} depicts the energy flow in OCP's network in two common reduced scenarios (as per our discussion in Section \ref{ssec:RO-uncertainty-analysis}, Figure \ref{fig:unc1}),} which respectively represent a sunnier, higher-generation day (Figure \ref{fig:budget_energy_flow} (left)) and a cloudier, lower-generation day (Figure \ref{fig:budget_energy_flow} (right)), aggregated across all sites and years in the planning horizon. The installed solar panels produce a lot of energy during both days, which is used to meet demand, refill OCP's batteries, and even occasionally sold to the grid. Moreover, the battery storage levels (reported in green) are kept sufficiently high by the solar energy generated during the day that OCP never purchases energy during the peak or shoulder periods of the day when the wholesale suppliers charge a higher unit price. Interestingly, in both scenarios, OCP purchases some electricity from the grid in the late evening (when electricity prices are low). This suggests that the model anticipates solar generation before it occurs and captures notions of load shifting, or arbitraging prices across time.  
}

\subsection{Summary of Main Findings}\label{ssec:exper-summary}
We now summarize the main managerial insights from our numerical results:
\begin{itemize}
    \item {\latest The use of robust and distributionally robust optimization can reduce the cost of decarbonization by $16\%$} or more (as demonstrated by our discussion in Section \ref{appx:exper-tuning}), which is significant in the context of a {\msom multi} billion-dollar investment in a renewable initiative by OCP. {\msom Considering} recent very significant investments in renewable energy to combat climate change, e.g., $330$ billion USD in new spending by the United States Congress to combat climate change via the Inflation Reduction Act \citep{biden2022inflation}, this suggests that robust optimization and other techniques from the Operations Management literature may have a significant role to play in future efforts to decarbonize.
    \item Figure \ref{fig:tradeooff} in Section \ref{ssec:exper-invest-vs-oper} demonstrates that a {\msom $10$ billion MAD} investment in solar panels and batteries reduces OCP's carbon emissions by approximately {\msom $70\%$} and is profitable over a $20$ year planning horizon. Moreover, a {\msom $20$ billion MAD investment in solar panels and batteries reduces OCP's carbon emissions by approximately {\latest $95\%$} and is also profitable}. This suggests that in addition to a global carbon arbitrage (where decreasing global $\text{CO}_2$ emissions increases global economic output) \citep[c.f.][]{adrian2022great}, there is—at least sometimes—a local carbon arbitrage, where decarbonizing is profitable for an energy-dependent manufacturer. {\msom However, achieving a $100\%$ reduction in carbon emissions is not profitable for OCP according to our model, which suggests that further regulations and/or subsidies which incentivize installing renewables may be needed to achieve the amount of decarbonization required to fulfill the Paris climate treaty.}
    \item Figure \ref{fig:budget_investment} in Section \ref{ssec:exper-inspect} demonstrates that when operating over a finite time horizon, investing in renewable energy generation earlier in the time horizon results in more cost savings for OCP. In a more general setting, these findings suggest that (a) other large manufacturers could consider decarbonization as a viable and potentially cheaper alternative to relying on increasingly scarce carbon-based energy sources, and (b) if they decarbonize, then investing in renewable energy sooner may be more profitable. This is particularly true given the recent surge in global energy prices since the Russian invasion of Ukraine \citep[c.f.][]{liadze2022economic}. 
\end{itemize}

\subsection{\msom Project Implementation and Impact of Research on OCP's Operations}\label{ssec:exper-impact}
OCP is currently using our approach as the basis of an approximately twenty billion MAD (approx. two billion USD) investment in solar panels and batteries, which our model estimates will simultaneously reduce OCP's carbon emissions which arise from energy consumption across the mining sites and factories considered by our model by {\latest $95\%$ and lower the time-adjusted operational costs of these assets by more than the investment cost over a $20$-year planning horizon. This investment is a significant subset of their recently announced $130$ billion MAD green initiative, which will fully decarbonize their operations by $2040$.}

From an impact perspective, our model has also influenced OCP's investment strategy by illustrating the value of batteries in decarbonizing their production process. Indeed, before our collaboration, OCP was contemplating a pure solar strategy. Still, our model demonstrates that decarbonizing by installing solar panels and batteries is more efficient than via solar panels alone, particularly after a saturation point where solar panels fully power OCP's system during the day.

In follow-up work, OCP is also considering expanding our model to consider other sites in their production pipeline, and incorporating wind generation to increase the cost efficiency of near $100\%$ decarbonization. OCP is also considering using renewable energy and batteries to power a new desalination plant that will help secure Morocco's water supply. 

\section{Conclusion}

{\msom 
We propose a comprehensive methodology developed in collaboration with OCP, one of the world's largest producers of phosphate and phosphate-based products, for partially decarbonizing its production pipeline using solar panels and batteries. Our methodology prescribes both the long-term strategic and daily tactical decisions that govern OCP's energy planning and is robust to the daily uncertainty in solar generation (using RO) and the uncertainty due to climate change (using DRO). }

\FloatBarrier

\if\anonymize0
\ACKNOWLEDGMENT{We are grateful to Yassine El Akel and Tarik Mortaji (OCP) for working with us on the data, implementation, and model outputs, the OCP Sustainability \& Green Industrial Development team for discussions, Kailyn Bryk for computational work, and Alexandre Jacquillat, J\'{o}nas Oddur J\'{o}nasson, Jean Pauphilet, Wolfram Wiesemann, and the associate editor and two referees for valuable comments which improved the manuscript.}
\fi


{\setlength{\bibsep}{0pt plus 0.0ex}
\scriptsize
\bibliographystyle{informs2014} 
\bibliography{thebib.bib} 
}

%
%
%
\ECSwitch
\ECHead{Supplementary Material}

\section{Omitted Proofs}\label{appx:proofs}
\subsection{Proof of Lemma \ref{lemma:detequiv}}

\proof{Proof. }

To reformulate constraint \eqref{eqn:unc-constraint-1} into its robust counterpart, we first expand the inner minimization problem:
\begin{equation} \label{eqn:unc-primal-1}
\begin{aligned}
    \min \quad & \sum_{n,h,d,m,y} \bar{z}_n^y  v^{h,d,m,y}_n  &\\ 
    \text{s.t.} \quad & v_n^{h,d,m,y} = \bar{v}^{h,d}+u^{h,d,m,y}_n \quad \forall h, d,m,y,n \quad & [\boldsymbol{\sigma_1}]\\
    & - U_{\text{MAX}}^{h,d} \leq u_n^{h,d,m,y} \leq U_{\text{MAX}}^{h,d} \quad \forall h, d,m,y,n \quad & [\boldsymbol{\phi_1}, \boldsymbol{\chi_1}]\\
    & - U_{\text{SV}}^{h,d} \leq u_n^{h,d,m,y} -  u_n^{h-1,d,m,y} \leq U_{\text{SV}}^{h,d} \quad \forall h, d,m,y,n \quad & [\boldsymbol{\psi_1},\boldsymbol{\omega_1}]\\
    & \sum_{n,h,m,y} u^{h,d,m,y}_n \leq U_{\text{CLT}}^{d} \quad \forall d \quad & [\boldsymbol{\tau_1}]\\
    & v_n^{h,d,m,y} \geq 0 \quad \forall  h, d,m,y,n &
\end{aligned}
\end{equation}
where, for simplicity, we take $u^{0,d}=u^{24,d},\forall d$, and denote the dual variable associated with each constraint in square brackets. Problem \eqref{eqn:unc-primal-1}'s dual is:
\begin{equation} \label{eqn:unc-dual-1}
\begin{split}
    \max \quad & \sum_{d} \left( \sum_{n,h,m,y} \bar{v}^{h,d} \sigma_{1,n}^{h,d,m,y} + U_{\text{MAX}}^{h,d} \phi_{1,n}^{h,d,m,y} + U_{\text{MAX}}^{h,d} \chi_{1,n}^{h,d,m,y} + U_{\text{SV}}^{h,d} \psi_{1,n}^{h,d,m,y} + U_{\text{SV}}^{h,d} \omega_{1,n}^{h,d,m,y} \right) + U_{\text{CLT}}^{d} \tau_1^d  \\ 
    \text{s.t.} \quad & \sigma_{1,n}^{h,d,m,y} \leq \bar{z}_n^y \quad \forall h, d,m,y,n \\
    & -\sigma_{1,n}^{h,d,m,y} - \phi_{1,n}^{h,d,m,y} +\chi_{1,n}^{h,d,m,y} - \psi_{1,n}^{h,d,m,y} + \omega_{1,n}^{h,d,m,y} +\psi_{1,n}^{h+1,d,m,y} - \omega_{1,n}^{h+1,d,m,y} + \tau_1^d = 0 \quad \forall h, d,m,y,n \\
    & \phi_{1,n}^{h,d,m,y} \geq 0,\ \chi_{1,n}^{h,d,m,y} \geq 0,\ \psi_{1,n}^{h,d,m,y} \geq 0,\ \omega_{1,n}^{h,d,m,y} \geq 0,\ \tau_1^d \geq 0 \quad \forall h, d,m,y,n
\end{split}
\end{equation}
where we take $\psi_{1,n}^{25,d,m,y}=\psi_{1,n}^{1,d,m,y},\ \omega_{1,n}^{25,d,m,y}=\omega_{1,n}^{1,d,m,y}$. Linear optimization duality then guarantees that we can replace the minimization problem (Problem \eqref{eqn:unc-primal-1}) in constraint \eqref{eqn:unc-constraint-1} with its dual (Problem \eqref{eqn:unc-dual-1}), and this yields the claimed result. \quad \Halmos
\endproof



\FloatBarrier
\section{\msom Operationalizing the SAA Model}\label{appx:SAA-operationalizing}
We have the following real-time model:
\begin{align}
    \min \quad & 
    \sum_{i=t}^T \left[ \sum_{a} c_{r,a}^{i} \vert f_{a}^{i}\vert 
    + \sum_{n} \left(
     p_{O,n}^{i} x_{O,n}^{i} 
     + p_{N,n}^{i} x_{N,n}^{i}
     - p_{w,n}^{i} w_{n}^{i} \right)
    \right] & \nonumber\\ 
    \text{s.t.} \quad & 
    \sum_{a \in \mathcal{I}(n)} \tau_a(f_{a}^{i})+ \sum_{a \in \mathcal{O}(n)} \tau_a(-f_{a}^{i}) + x_{O,n}^{i}+x_{N,n}^{i}-w_{n}^{i} + R\cdot r_{n}^{i} \nonumber\\
    & \geq d_{n}^{i} - v_n^{i} \left( \sum_{y'=1}^y \xi^{y-y'} \tilde{z}_n^{y'} \right) ,\quad & \forall n,\ t \leq i \leq T ,\nonumber\\
    & \sum_{i=t}^T w_n^i \leq \beta  \sum_{i=1}^{t-1} \left[  v_n^{i} \left( \sum_{y'=1}^y \xi^{y-y'} \tilde{z}_n^{y'} \right) + \max\left\{ 0, -d_n^{i} \right\} \right] - \sum_{i=1}^{t-1} \tilde{w}_{n}^{i} , \quad & \forall n ,\nonumber\\ 
    & s^{t}_n=\psi \tilde{s}_n^{t-1}-\tilde{r}_{n}^{t-1},\ s^{i+1}_n=\psi s_n^{i-1}-r_{n}^{i-1} , \quad & \forall n,\ t < i < T, \nonumber\\
    &  s^{T+1}_n=\psi s_n^{T}-r_{n}^{T},  \quad & \nonumber\\
    & s^{i}_n \leq \sum_{y'=1}^y \nu^{y-y'} \tilde{b}^{y'}_n , \quad & \forall n,\ t \leq i \leq T, \nonumber\\
    & x_{O,n}^{i} \leq G_o^i, \ x_{N,n}^{i} \leq G_n^i  ,& \forall n,\ t \leq i \leq T , \nonumber\\
    & \vert f_{a}^{i}\vert \leq K_a, \ s_{n}^{i}, x_{O,n}^{i}, x_{N,n}^{i}, w_{n}^{i} \geq 0, \quad & \forall a,\ n,\ t \leq i \leq T,\label{eqn:operational}
\end{align}
where $s^{T+1}_n$ is the amount of energy stored by OCP at the end of the day, which for simplicity we set to be equal to the amount of energy stored at the start of the day in order that the optimizer does not completely drain the batteries in time period $T$. Although we do not implement this here, we could also consider more complex approaches which address end-of-horizon effects more accurately, e.g., extending the planning horizon over multiple days or rewarding energy storage levels at the end of the day via a salvage function \citep[see][]{grinold1983model}. 

\section{\msom Net Present Value of Decarbonizing OCP} \label{appx:exper-NPV}
\begin{figure}[!ht] 
    \centering
    \includegraphics[width=0.48\linewidth]{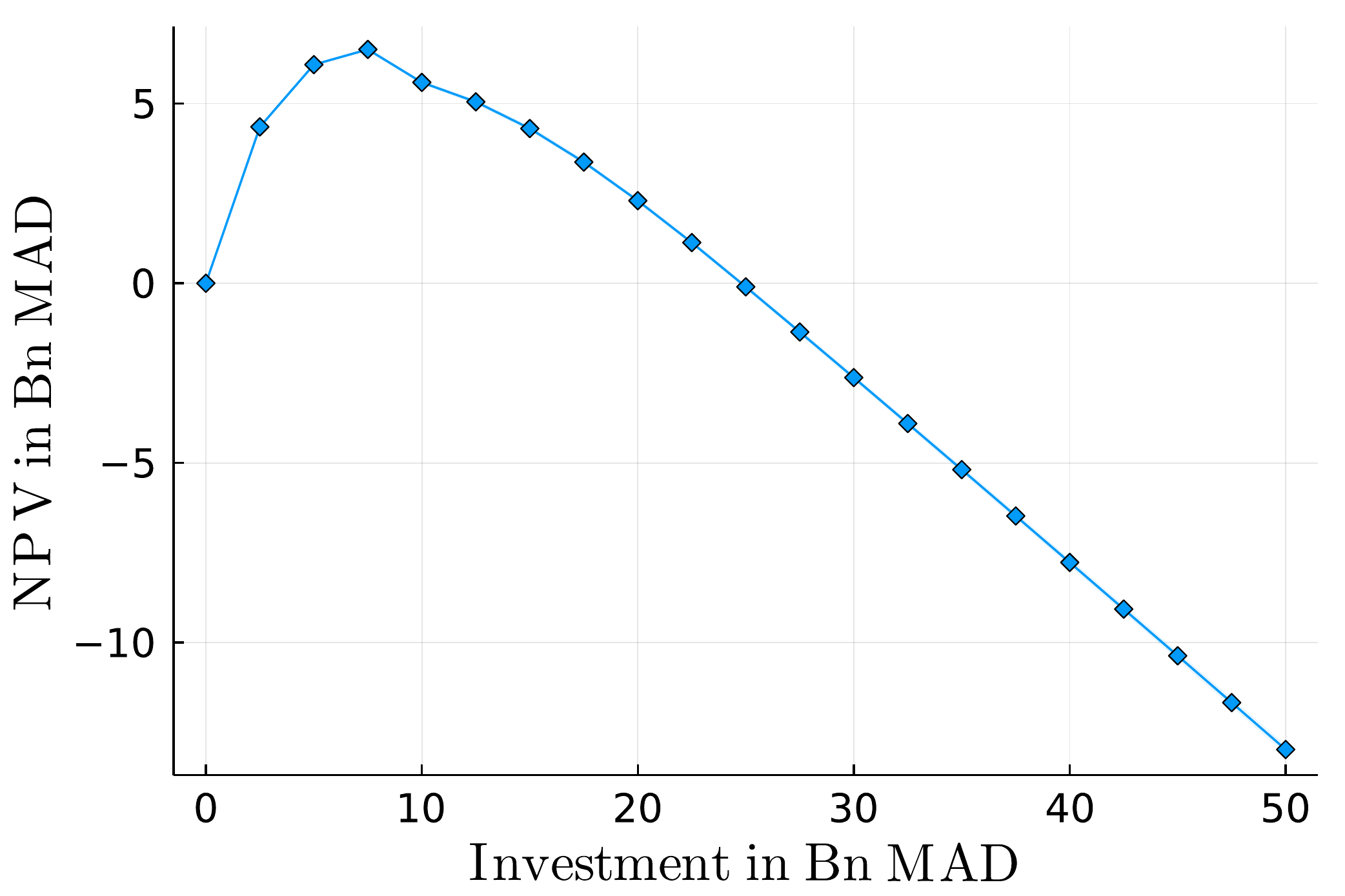}
    \caption{Trade-off between overall investment cost and the project's net present value over the planning horizon. Increasing the investment budget increases the project's value up to $7.5$ billion MAD, and subsequently decreases the project's value.}
\label{fig:tradeooff-npv}
\end{figure}

In this section, we complement Figure \ref{fig:tradeooff} by studying the trade-off between the size of OCP's investment and the project's net present value (NPV). {\latest Denoting by $\text{cost}_z(\boldsymbol z^t)$ the cost of the strategic decisions at time $t$ and $\text{cost}_x(\boldsymbol x^t; \boldsymbol z^1, \dots, \boldsymbol z^t)$ the cost of the operational decisions at time $t$ given the strategic decisions made until that time, we compute the NPV as follows:
$$
\text{NPV}(\boldsymbol z, \boldsymbol x) 
= \sum_t (\rho)^t \left( 
\text{cost}_x(\boldsymbol x^t; \boldsymbol 0, \dots, \boldsymbol 0) 
- \text{cost}_x(\boldsymbol x^t; \boldsymbol z^1, \dots, \boldsymbol z^t)
- \text{cost}_z(\boldsymbol z^t)
\right).
$$
Note that we do not include the salvage value of the solar panels and batteries at the end of the planning horizon when computing the NPV. We observe in Figure \ref{fig:tradeooff-npv} that the NPV of the project is approximately maximized at an investment of 7.5 billion MAD, although an investment of 10 billion MAD yields nearly the same return. }

{\latest
\section{\msom Cross-Validation for the RO and DRO Hyperparameters}
\label{appx:exper-tuning}
In this section, we propose a methodology for cross-validating the RO and DRO hyperparameters and demonstrate that properly selecting these parameters leads to out-of-sample cost savings for OCP. We also investigate the cross-validated model's sensitivity to distributional shifts between the training data and the testing data, e.g., induced by climate change. 

\subsection{Cross-Validation Methodology}  \label{appx:exper-tuning-methodology}
Following a standard machine learning paradigm, we randomly split the year of site-wise solar generation data supplied by OCP into training, validation, and testing sets on a month-by-month basis (with a training/validation/test split of 20 days/4 days/4 days). Using the training data, we estimate reduced scenarios using the methodology laid out in Section \ref{ssec:scenario-reduction}, and solve Problem \eqref{prob:dro_det} to optimality for each value of the hyperparameters \vspace{-3mm}$$(\gamma_{\max}, \gamma_C, \gamma_{\text{CLT}}, \delta) \in \left\{0, 1\right\}^3 \times \left\{0, 0.001, 0.01, 0.1, 1.0\right\}.$$ To select the optimal hyperparameters, we evaluate the performance of each model on the validation set, and select the combination of parameters that performs best in terms of minimizing both the cost of investment plus operations and the CO$_2$ emissions. To reduce the dimensionality of the problem, we set the number of reduced scenarios to $10$ and use an investment budget of $20$ billion MAD. Moreover, to reduce the sensitivity of our approach to noise, we repeat the entire validation process ten times and take the hyperparameters that perform best on average.


\subsection{Effect of Robustness on the Investment} 
We first study the effect of introducing robustness in the model on the amount invested in solar panels (Figure \ref{fig:RObudgets} (left)) and in batteries (Figure \ref{fig:RObudgets} (right)), with a constant investment budget of {$20$ billion MAD}. We make the following observations: 
{ \latest
\begin{itemize}
    \item Increasing the amount of daily uncertainty in the problem by increasing the RO uncertainty set size leads to a larger investment in solar panels and a smaller investment in batteries. This is perhaps surprising, since batteries offer OCP a recourse action after uncertainty in the system is revealed, while solar panels do not. However, it can be explained by three substructures in OCP’s problem. First,  roughly  speaking,  batteries  allow  load  shifting  to  take  place  and  reduce  the  overall  cost accrued in hours of the day where the sun does not shine, while solar panels reduce the cost accrued in hours when the sun does shine. Second, in OCP’s problem, the marginal price of electricity is highest in peak periods when the sun does shine. Third, due to the nature of the uncertainty sets designed in Section \ref{sec:RO}, increasing the size of the RO uncertainty allows nature to increase the relative cost of the most expensive hours of the day. Correspondingly, when the overall investment budget is held constant, increasing the amount of robustness increases the amount invested in solar panels and decreases the amount invested in batteries.

    \item Increasing the amount of long-term uncertainty in the problem by increasing the DRO ambiguity set size (controlled by hyperparameter $\delta$) results in a smaller investment in both solar panels and batteries, meaning that the full investment budget of 20 billion MAD is not utilized in our robust models. This can be explained by the fact that if there are diminishing returns in investing more capital in solar panels and batteries, then being more robust increases the model's risk aversion and decreases its appetite for investment. As reflected in Figure \ref{fig:budget_investment}, this issue of diminishing marginal returns certainly occurs with an investment budget of $20$ billion MAD.
\end{itemize}



\begin{figure}[!ht] 
\begin{adjustbox}{minipage=\linewidth,scale=1}
    \centering
    \begin{subfigure}[b]{0.48\linewidth}
    
        \centering
        \includegraphics[width=\linewidth]{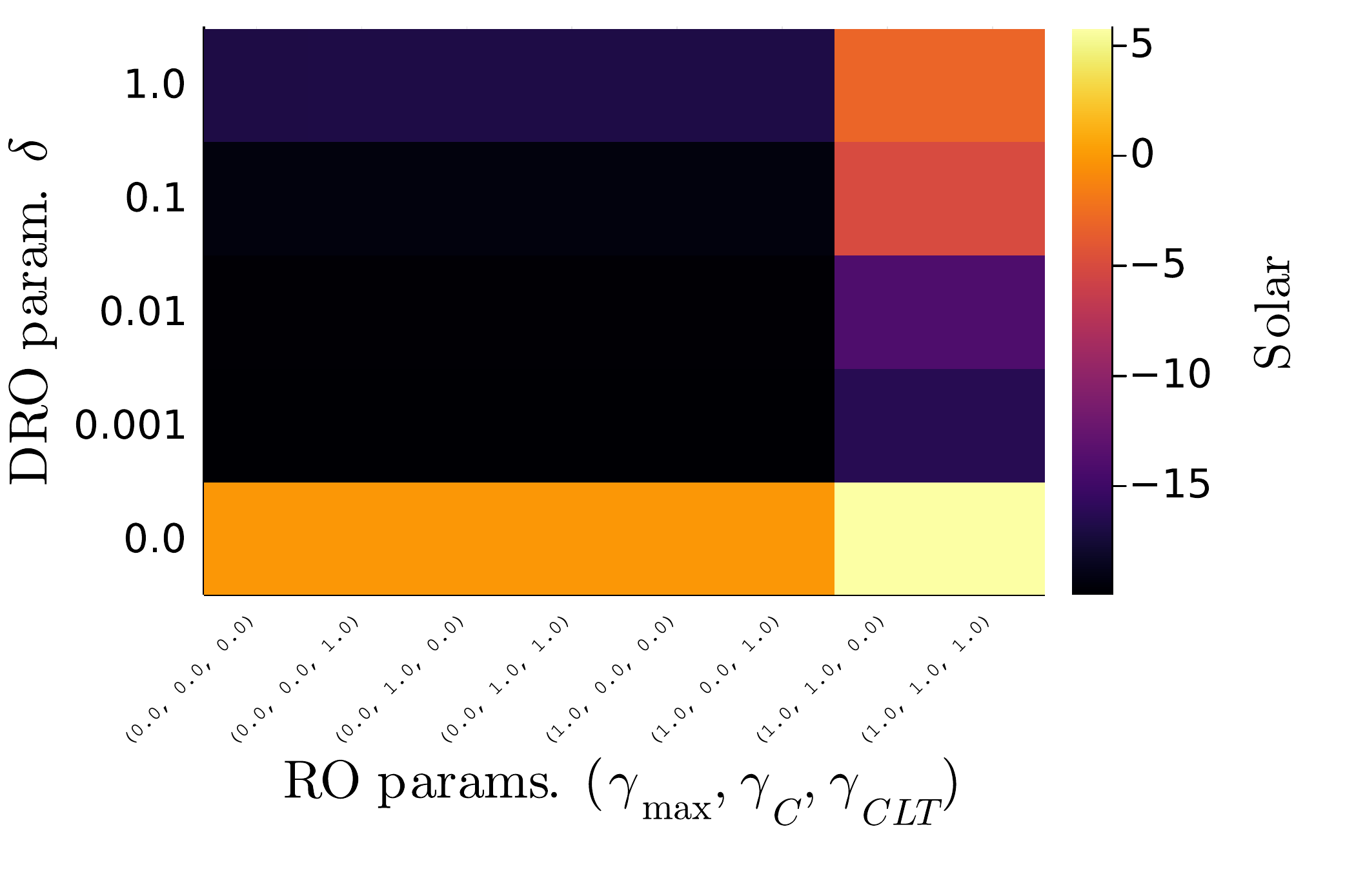}
        \label{fig:robust_battery_heatmat}
    
    \end{subfigure}
    \hfill
    \begin{subfigure}[b]{0.48\linewidth}
    
        \centering
        \includegraphics[width=\linewidth]{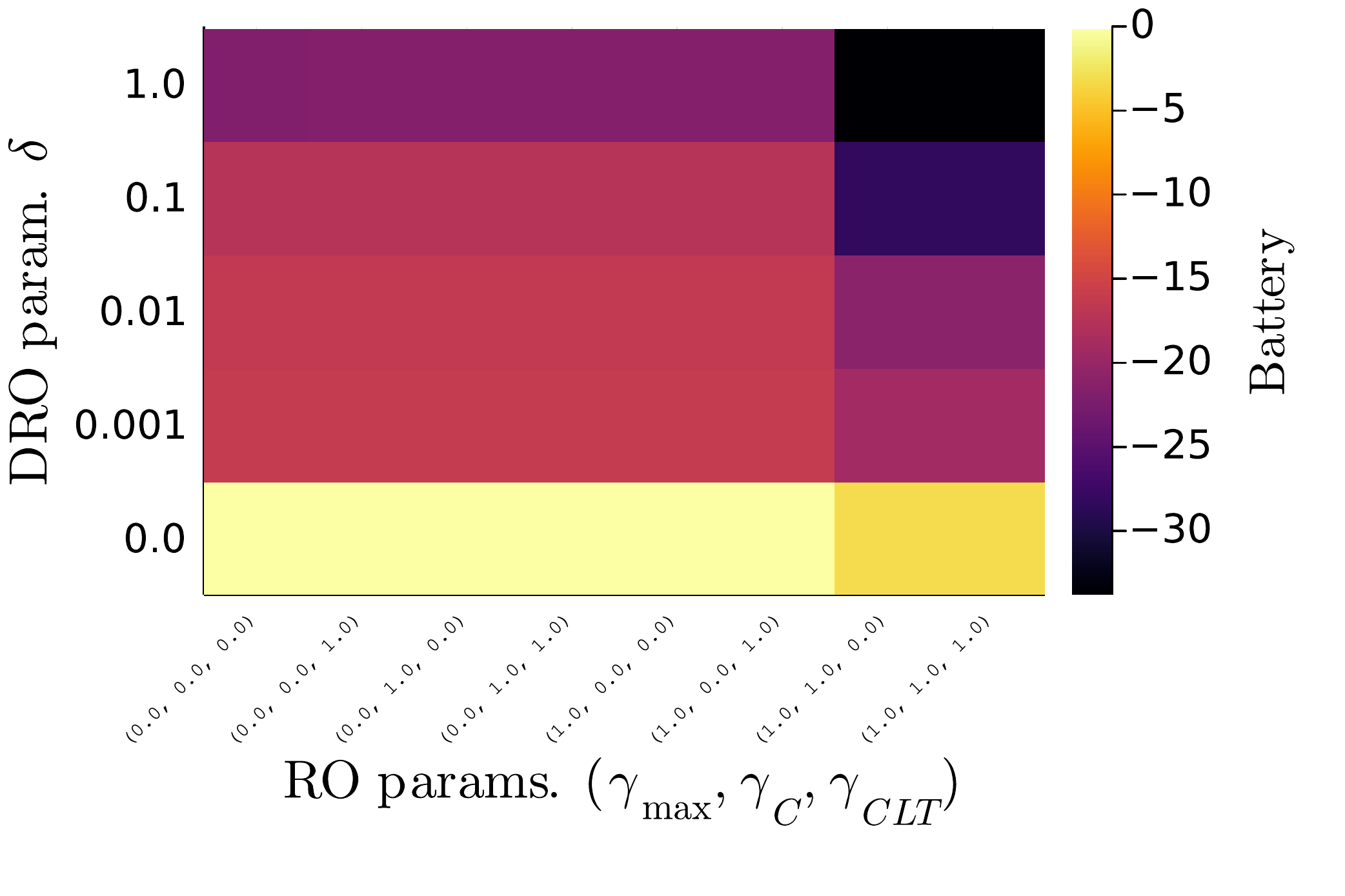}
        \label{fig:robust_solar_heatmat}
    
    \end{subfigure}
\end{adjustbox}
\caption{Effect of RO and DRO uncertainty budgets on investment in solar panels and batteries with a {$20$ billion MAD }investment budget.}
\label{fig:RObudgets}
\end{figure}

\subsection{Effect of Robustness on Operational Costs and CO$_2$ Emissions} 
We now investigate the average relative improvement in operational costs and CO$_2$ emissions reduction predicted by the validation set, compared to the SAA model, for a selected set of values of the RO and DRO hyperparameters. 
As depicted in Figure \ref{fig:robust_validation_heatmat}, the RO- and DRO-guarded models provide an improvement of $16.3\%$ in terms of operational costs and $3.5\%$ in terms of CO$_2$ emissions reduction over the SAA model. Nonetheless, being overly conservative harms the model's performance (as can be seen from the top and rightmost blocks of Figure \ref{fig:robust_validation_heatmat}). We present the percentage improvement in operational costs on the validation set using different robust hyperparameter values in Table \ref{tab:validation} (we independently repeat the cross-validation process/split the data ten times, so we report both the mean and standard deviation of the results). On the validation dataset, the combination of hyperparameters that performs best are $(\gamma_{\max},\gamma_C,\gamma_{\text{CLT}},\delta)=(1.0,0.0,1.0,0.001)$ and $(\gamma_{\max},\gamma_C,\gamma_{\text{CLT}},\delta)=(1.0,0.0,1.0,1.0)$ in terms of minimizing, respectively, the cost of investment plus operations and the CO$_2$ emissions.

\begin{figure}[!ht] 
\begin{adjustbox}{minipage=\linewidth,scale=1}
    \centering
    \begin{subfigure}[b]{0.48\linewidth}
    
        \centering
        \includegraphics[width=\linewidth]{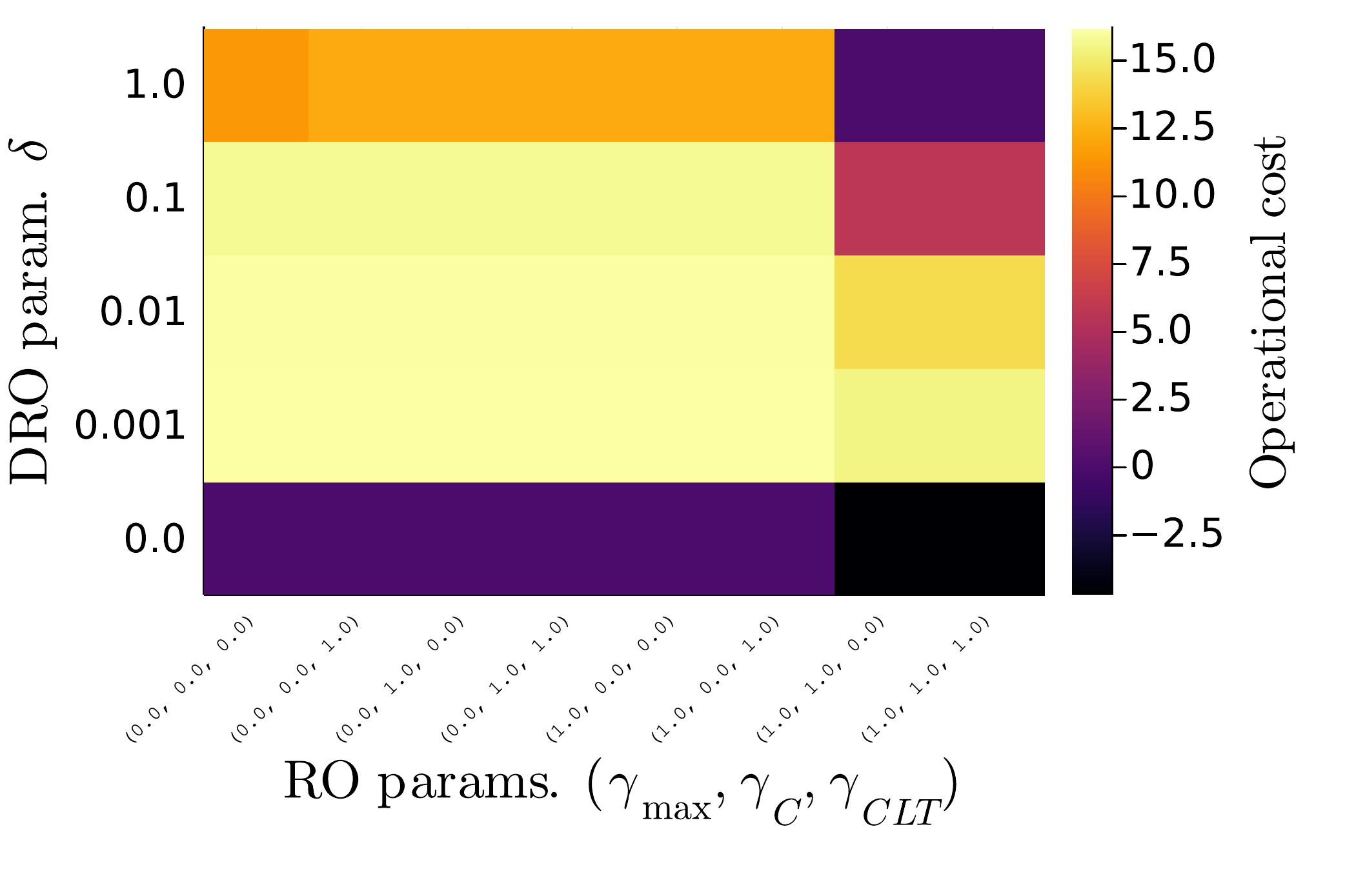}
    
    \end{subfigure}
    \hfill
    \begin{subfigure}[b]{0.48\linewidth}
    
        \centering
        \includegraphics[width=\linewidth]{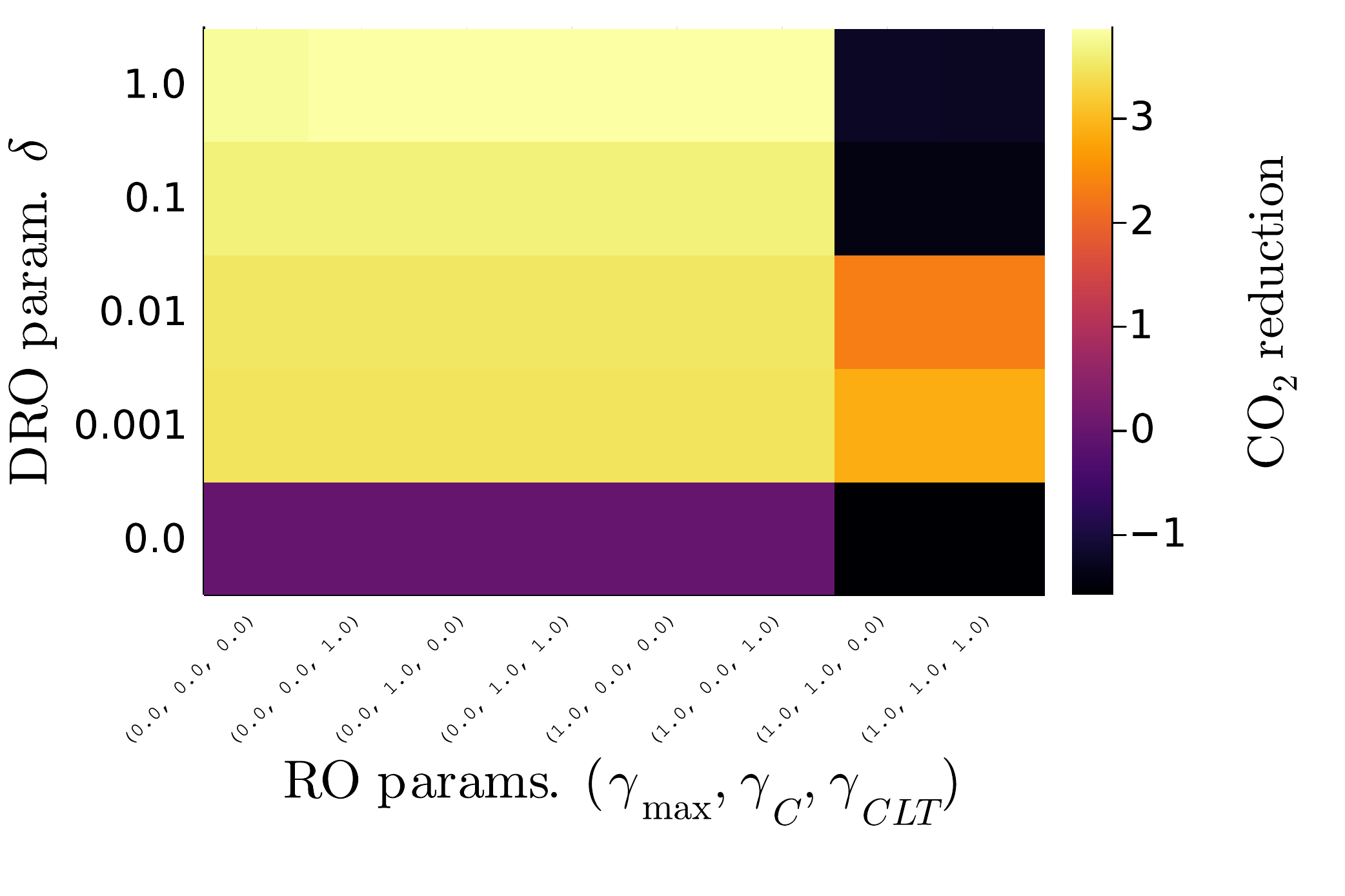}
    
    \end{subfigure}
    \end{adjustbox}
    \caption{{Cross-validating the RO and DRO budgets improves OCP's operational costs (by up to $16.3\%$) and reduces CO$_2$ emissions (by up to $3.5\%$) using validation data. }}
    \label{fig:robust_validation_heatmat}
\end{figure}

\begin{table}[!ht]\footnotesize
\centering
\caption{Mean and standard deviation of $\%$ improvement in operational costs (top half) and in CO$_2$ reduction (bottom half) compared to the SAA model using validation data. \label{tab:validation}}
\begin{tabular}{rrrrrr}
\toprule
$(\gamma_{\max},\gamma_C,\gamma_{CLT})$ & $\delta=$0.0 & $\delta=$0.001 & $\delta=$0.01 & $\delta=$0.1 & $\delta=$1.0 \\ \midrule
(0.0, 0.0, 0.0) & 0.0 (0.0) & 16.26 (1.0) & 16.19 (0.99) & 15.91 (0.97) & 11.53 (2.62) \\
(0.0, 0.0, 1.0) & 0.0 (0.0)   & 16.26 (1.0) & 16.19 (0.99) & 15.91 (0.97) & 12.32 (1.96) \\
(0.0, 1.0, 0.0) & 0.0 (0.0)   & 16.26 (1.0) & 16.19 (0.99) & 15.91 (0.97) & 12.32 (1.96) \\
(0.0, 1.0, 1.0) & 0.0 (0.0)   & 16.26 (1.0) & 16.19 (0.99) & 15.91 (0.97) & 12.32 (1.96) \\
(1.0, 0.0, 0.0) & 0.0 (0.0)   & 16.26 (1.0) & 16.19 (0.99) & 15.91 (0.97) & 12.32 (1.96) \\
(1.0, 0.0, 1.0) & 0.0 (0.0)   & \textbf{16.26 (1.0)} & 16.19 (0.99) & 15.91 (0.97) & 12.32 (1.96) \\
(1.0, 1.0, 0.0) & -4.69 (2.02) & 15.54 (1.1) & 14.43 (1.61) & 5.78 (5.89)  & 0.02 (7.19) \\
(1.0, 1.0, 1.0) & -4.69 (2.02) & 15.54 (1.1) & 14.43 (1.61) & 5.78 (5.89)  & 0.02 (7.19) \\ \midrule
(0.0, 0.0, 0.0) & 0.0 (0.0)  & 3.09 (0.45) & 3.12 (0.45) & 3.25 (0.46)  & 3.46 (0.53)          \\
(0.0, 0.0, 1.0) & 0.0 (0.0)  & 3.09 (0.45) & 3.12 (0.45) & 3.25 (0.46)  & 3.51 (0.52)          \\
(0.0, 1.0, 0.0) & 0.0 (0.0)  & 3.09 (0.45) & 3.12 (0.45) & 3.25 (0.46)  & 3.51 (0.52)          \\
(0.0, 1.0, 1.0) & 0.0 (0.0)  & 3.09 (0.45) & 3.12 (0.45) & 3.25 (0.46)  & 3.51 (0.52)          \\
(1.0, 0.0, 0.0) & 0.0 (0.0)  & 3.09 (0.45) & 3.12 (0.45) & 3.25 (0.46)  & 3.51 (0.52)          \\
(1.0, 0.0, 1.0) & 0.0 (0.0)  & 3.09 (0.45) & 3.12 (0.45) & 3.25 (0.46)  & \textbf{3.51 (0.52)} \\
(1.0, 1.0, 0.0) & -1.74 (0.71) & 2.55 (0.53) & 2.06 (0.8)  & -1.55 (2.99) & -1.36 (2.25)         \\
(1.0, 1.0, 1.0) & -1.74 (0.71) & 2.55 (0.52) & 2.06 (0.8)  & -1.55 (2.99) & -1.36 (2.25)        
 \\ \bottomrule
\end{tabular}
\end{table}

\subsection{Model's Out-Of-Sample Performance}
Next, we verify the model's out-of-sample performance using holdout data. In particular, in Table \ref{tab:test}, we utilize the testing set that we reserved as part of our cross-validation methodology (described in Section \ref{appx:exper-tuning-methodology}). In Table \ref{tab:perturbed}, we generate a ``perturbed testing set'' by perturbing each entry $v$ in our hourly solar generation testing data according to $\tilde{v} = v \cdot \left[ 1 + \mathcal{N}(\mu, (\nicefrac{\mu}{10})^2) \right],$ where $\mathcal{N}$ denotes the normal distribution and $\mu$ is drawn uniformly at random from $[-0.25, 0.25].$ We again present the average relative improvement in operational costs and CO$_2$ emissions reduction compared to the SAA model, for the same selected set of values of the RO and DRO hyperparameters. In both cases, the out-of-sample results are consistent with the estimates obtained using validation data. 

    
    
    
    

\begin{table}[!ht]\footnotesize
\centering
\caption{Mean and standard deviation of $\%$ improvement in operational costs (top half) and in CO$_2$ reduction (bottom half) compared to the SAA model using testing data. \label{tab:test}}
\begin{tabular}{rrrrrr}
\toprule
$(\gamma_{\max},\gamma_C,\gamma_{CLT})$ & $\delta=$0.0 & $\delta=$0.001 & $\delta=$0.01 & $\delta=$0.1 & $\delta=$1.0 \\ \midrule
(0.0, 0.0, 0.0) & 0.0 (0.0) & 16.51 (0.87) & 16.45 (0.88) & 16.18 (0.98) & 11.81 (2.46) \\
(0.0, 0.0, 1.0) & 0.0 (0.0) & 16.51 (0.87) & 16.45 (0.88) & 16.18 (0.98) & 12.63 (2.05) \\
(0.0, 1.0, 0.0) & 0.0 (0.0) & 16.51 (0.87) & 16.45 (0.88) & 16.18 (0.98) & 12.63 (2.05) \\
(0.0, 1.0, 1.0) & 0.0 (0.0) & 16.51 (0.87) & 16.45 (0.88) & 16.18 (0.98) & 12.63 (2.05) \\
(1.0, 0.0, 0.0) & 0.0 (0.0) & 16.51 (0.87) & 16.45 (0.88) & 16.18 (0.98) & 12.63 (2.05) \\
(1.0, 0.0, 1.0) & 0.0 (0.0) & \textbf{16.51 (0.87)} & 16.45 (0.88) & 16.18 (0.98) & 12.63 (2.05) \\
(1.0, 1.0, 0.0) & -4.75 (2.1)  & 15.78 (1.16) & 14.63 (1.6)  & 5.83 (5.97)  & 0.11 (7.09)  \\
(1.0, 1.0, 1.0) & -4.75 (2.1)  & 15.78 (1.16) & 14.63 (1.6)  & 5.83 (5.97)  & 0.11 (7.09)  \\ \midrule
(0.0, 0.0, 0.0) & 0.0 (0.0) & 3.47 (0.27) & 3.51 (0.28) & 3.63 (0.3)   & 3.8 (0.43)          \\
(0.0, 0.0, 1.0) & 0.0 (0.0) & 3.47 (0.27) & 3.51 (0.28) & 3.63 (0.3)   & 3.86 (0.4)          \\
(0.0, 1.0, 0.0) & 0.0 (0.0) & 3.47 (0.27) & 3.51 (0.28) & 3.63 (0.3)   & 3.86 (0.4)          \\
(0.0, 1.0, 1.0) & 0.0 (0.0) & 3.47 (0.27) & 3.51 (0.28) & 3.63 (0.3)   & 3.86 (0.4)          \\
(1.0, 0.0, 0.0) & 0.0 (0.0) & 3.47 (0.27) & 3.51 (0.28) & 3.63 (0.3)   & 3.86 (0.4) \\
(1.0, 0.0, 1.0) & 0.0 (0.0) & 3.47 (0.27) & 3.51 (0.28) & 3.63 (0.3)   & \textbf{3.86 (0.4)}          \\
(1.0, 1.0, 0.0) & -1.58 (0.7)  & 2.87 (0.47) & 2.33 (0.82) & -1.41 (3.06) & -1.24 (2.31)        \\
(1.0, 1.0, 1.0) & -1.58 (0.7)  & 2.87 (0.47) & 2.33 (0.82) & -1.41 (3.06) & -1.24 (2.31)          
 \\ \bottomrule
\end{tabular}
\end{table}

\begin{table}[!ht]\footnotesize
\centering
\caption{Mean and standard deviation of $\%$ improvement in operational costs (top half) and in CO$_2$ reduction (bottom half) compared to the SAA model using perturbed testing data. \label{tab:perturbed}}
\begin{tabular}{rrrrrr}
\toprule
$(\gamma_{\max},\gamma_C,\gamma_{CLT})$ & $\delta=$0.0 & $\delta=$0.001 & $\delta=$0.01 & $\delta=$0.1 & $\delta=$1.0 \\ \midrule
(0.0, 0.0, 0.0) & 0.0 (0.0) & 16.92 (1.64) & 16.89 (1.67) & 16.73 (1.81) & 12.86 (3.14) \\
(0.0, 0.0, 1.0) & 0.0 (0.0) & 16.92 (1.64) & 16.89 (1.67) & 16.73 (1.81) & 13.67 (2.86) \\
(0.0, 1.0, 0.0) & 0.0 (0.0) & 16.92 (1.64) & 16.89 (1.67) & 16.73 (1.81) & 13.67 (2.86) \\
(0.0, 1.0, 1.0) & 0.0 (0.0) & 16.92 (1.64) & 16.89 (1.67) & 16.73 (1.81) & 13.67 (2.86) \\
(1.0, 0.0, 0.0) & 0.0 (0.0) & 16.92 (1.64) & 16.89 (1.67) & 16.73 (1.81) & 13.67 (2.86) \\
(1.0, 0.0, 1.0) & 0.0 (0.0) & \textbf{16.92 (1.64)} & 16.89 (1.67) & 16.73 (1.81) & 13.67 (2.86) \\
(1.0, 1.0, 0.0) & -4.59 (2.9) & 16.53 (1.85) & 15.56 (2.3)  & 7.1 (6.55)   & 1.84 (6.93)  \\
(1.0, 1.0, 1.0) & -4.59 (2.9) & 16.53 (1.85) & 15.56 (2.3)  & 7.1 (6.55)   & 1.84 (6.93)  \\ \midrule
(0.0, 0.0, 0.0) & 0.0 (0.0)  & 3.12 (0.29) & 3.16 (0.3) & 3.29 (0.33) & 3.56 (0.46)          \\
(0.0, 0.0, 1.0) & 0.0 (0.0)  & 3.12 (0.29) & 3.16 (0.3) & 3.29 (0.33) & 3.61 (0.43)          \\
(0.0, 1.0, 0.0) & 0.0 (0.0)  & 3.12 (0.29) & 3.16 (0.3) & 3.29 (0.33) & 3.61 (0.43)          \\
(0.0, 1.0, 1.0) & 0.0 (0.0)  & 3.12 (0.29) & 3.16 (0.3) & 3.29 (0.33) & 3.61 (0.43)          \\
(1.0, 0.0, 0.0) & 0.0 (0.0)  & 3.12 (0.29) & 3.16 (0.3) & 3.29 (0.33) & 3.62 (0.43)          \\
(1.0, 0.0, 1.0) & 0.0 (0.0)  & 3.12 (0.29) & 3.16 (0.3) & 3.29 (0.33) & \textbf{3.62 (0.43)} \\
(1.0, 1.0, 0.0) & -1.43 (0.63) & 2.66 (0.32) & 2.24 (0.6) & -1.2 (2.86) & -0.97 (2.15)         \\
(1.0, 1.0, 1.0) & -1.43 (0.63) & 2.66 (0.32) & 2.24 (0.6) & -1.2 (2.86) & -0.97 (2.15)         
 \\ \bottomrule
\end{tabular}
\end{table}

\subsection{The Impact of Robustness}
Combining insights from our cross-validation procedure with discussions with the OCP team, we set the RO and DRO hyperparameters as follows:
\begin{itemize}
    \item For RO, the following $(\gamma_{\max},\gamma_C,\gamma_{\text{CLT}})$ hyperparameter combinations perform equally well across all values of the DRO hyperparameter $\delta$: $(1,0,0),(0,1,0),(0,0,1),(0,1,1),(1,0,1)$. As long as we are not overly conservative by setting either both $\gamma_{\max}$ and $\gamma_C$ or all three RO hyperparameters to $1$, RO provides a small yet non-negligible improvement. In addition, we can obtain such improvement from all three types of RO constraints. To prevent conservatism (by excluding the worst-case realizations of the uncertainty in our solar generation data), while at the same time incorporating all three types of RO constraints, we set $(\gamma_{\max},\gamma_C,\gamma_{\text{CLT}}) = (0.5, 0.5, 0.5).$

    \item For DRO, $\delta=0.001$ performs best in terms of minimizing the cost of investment plus operations, whereas $\delta=1.0$ performs best in terms of minimizing the CO$_2$ emissions. Therefore, we set $\delta=0.01.$
\end{itemize}

Figure \ref{fig:robust_distributions} compares the distributions of out-of-sample
operational costs (left) and CO$_2$ emissions reductions (right)
between the SAA and the cross-validated robust model.
To obtain such distributions, we perform multiple (ten) training/validation/testing splits and, for each such split, we evaluate both models' performance using non-perturbed and perturbed testing data (for various levels of perturbation). Note that the randomness is due to both the data-splitting process and the random perturbations we introduce. We observe that the cross-validated robust model saves an additional $0.4$ billion MAD in operational costs and further reduces CO$_2$ emissions by $3.5\%$
compared to the SAA model.

\begin{figure}[!ht] 
\begin{adjustbox}{minipage=\linewidth,scale=1}
    \centering
    \begin{subfigure}[b]{0.48\linewidth}
    
        \centering
        \includegraphics[width=\linewidth]{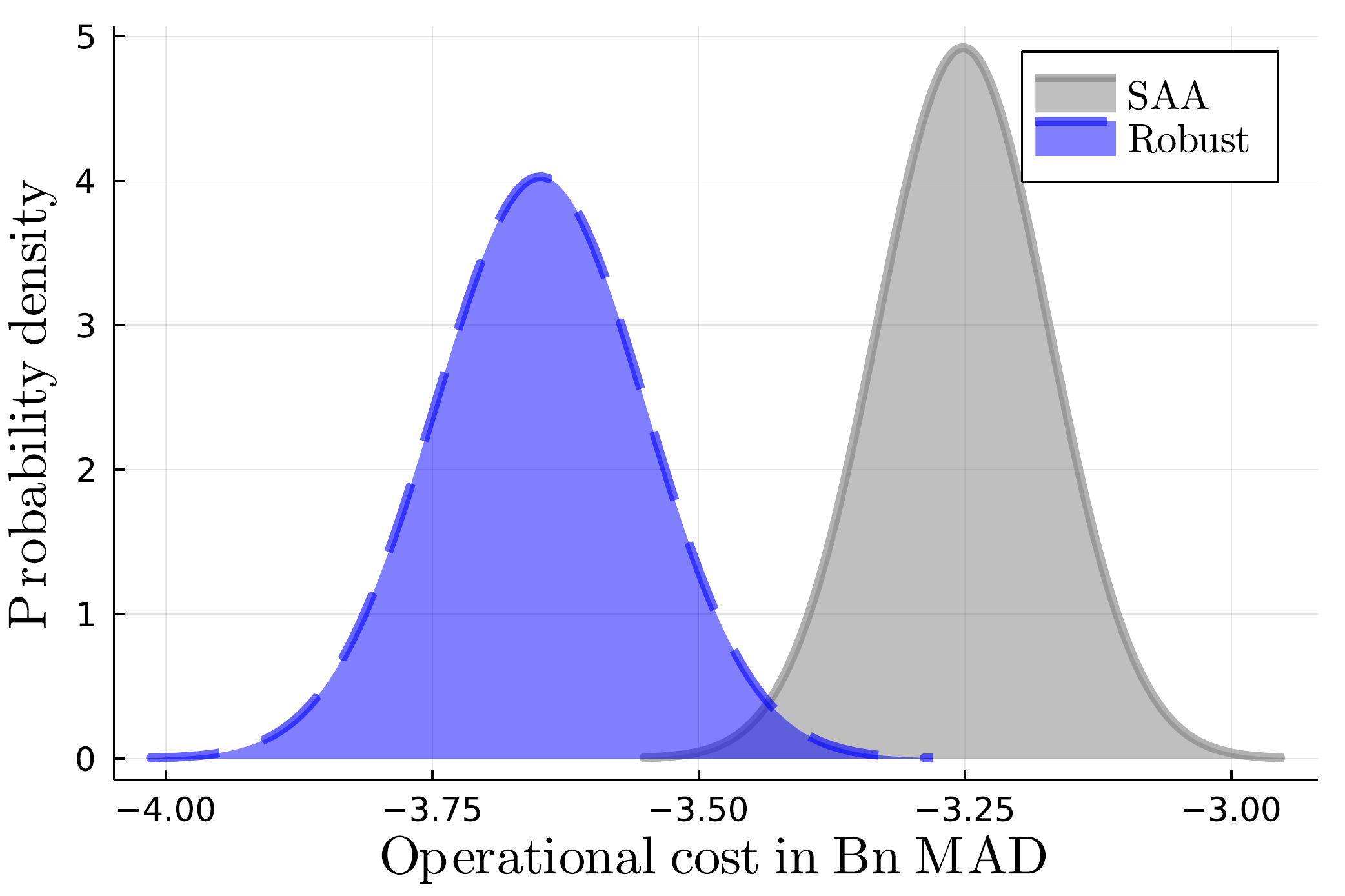}
    
    \end{subfigure}
    \hfill
    \begin{subfigure}[b]{0.48\linewidth}
    
        \centering
        \includegraphics[width=\linewidth]{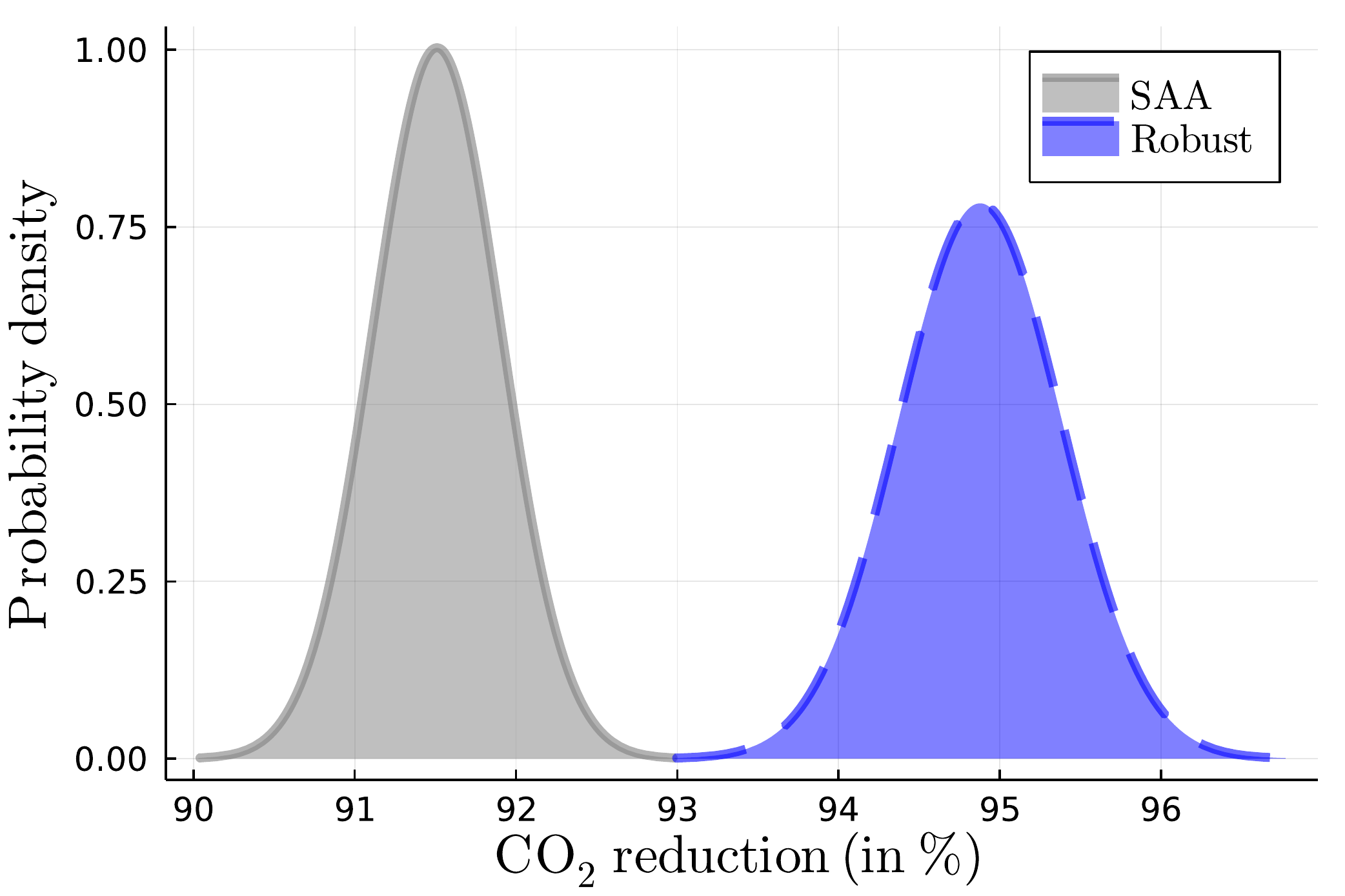}
    
    \end{subfigure}
    \end{adjustbox}
    \caption{{The cross-validated robust model leads to significant operational cost savings and CO$_2$ emissions reductions using testing data. }}
    \label{fig:robust_distributions}
\end{figure}

}

{\latest 

\section{Sensitivity Analysis with Respect to Reduced Scenarios}
\label{appx:exper-sensitivity}
In this section, we investigate the cross-validated model's sensitivity to the number of reduced scenarios and the assumption that scenarios can only be succeeded by scenarios of the same type. 

\subsection{Sensitivity to Number of Scenarios} 
We first explore the model's sensitivity to the number of reduced scenarios. We use a horizon of $20$ years, an investment budget of $20$ billion MAD, set the RO and DRO hyperparameters to $(\gamma_{\max},\gamma_C,\gamma_{\text{CLT}},\delta)=(0.5,0.5,0.5,0.01)$ as suggested by our cross-validation procedure, and vary the number of scenarios in $\{3,5,10,15,20,30,40,50\}$. Figure \ref{fig:n_typical_days_num_scenarios} depicts the relationship between the number of scenarios and the improvement in OCP's long-run operational costs compared to a naive baseline of not installing solar panels or batteries and satisfying OCP's energy needs via purchasing electricity from the grid locally at each site. \textcolor{blue}{We conclude that the improvement in operational cost (top left), the total reduction in carbon emissions (top right), and the investment policy prescribed by the model (bottom) do not vary significantly as we increase the number of scenarios. }

\begin{figure}[!ht] 
\begin{adjustbox}{minipage=\linewidth,scale=1}
    \centering
    \begin{subfigure}[b]{0.48\linewidth}
    
        \centering
        \includegraphics[width=\linewidth]{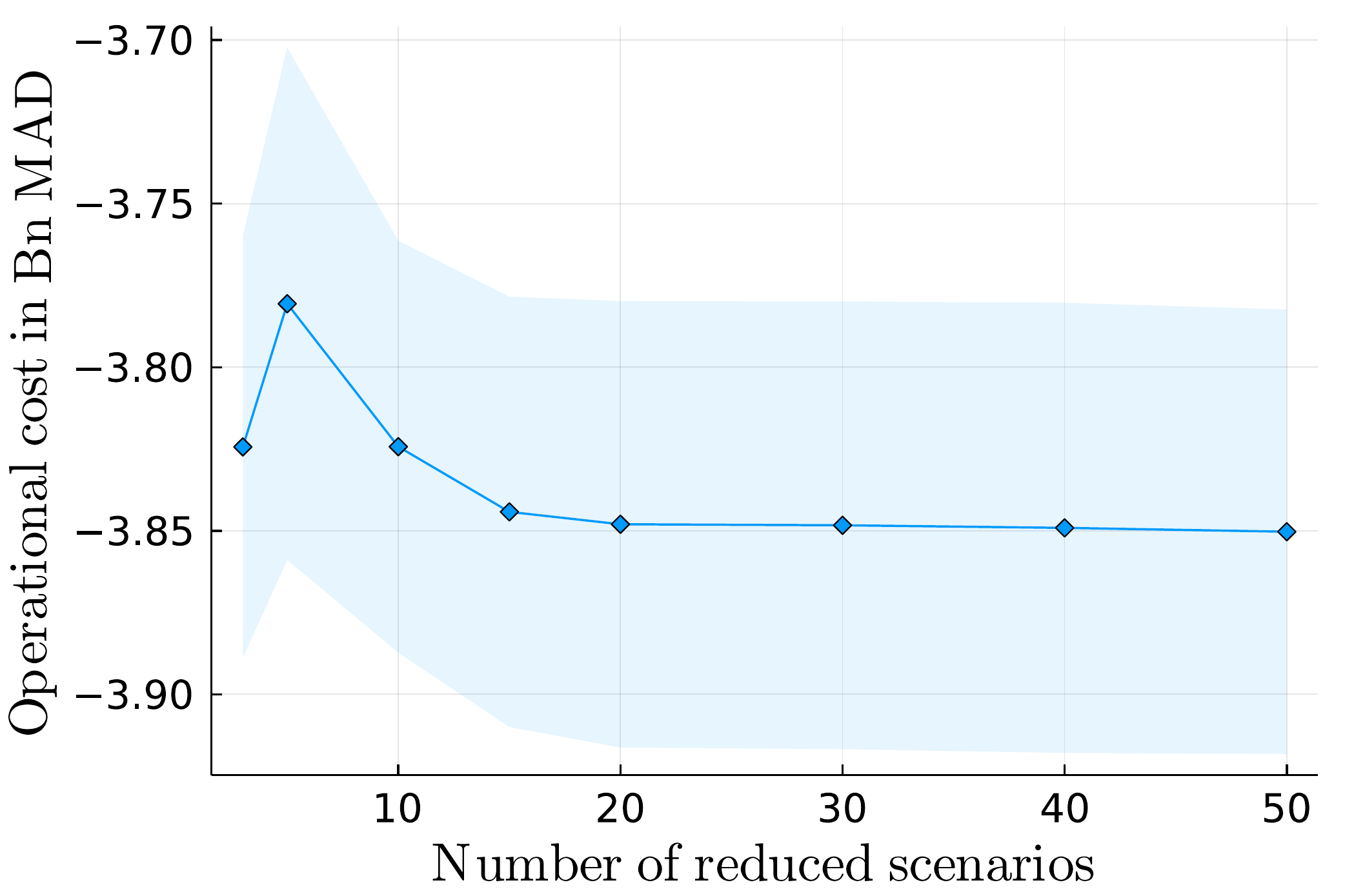}
        \label{fig:n_typical_days_operational}
    
    \end{subfigure}
    \hfill
    \begin{subfigure}[b]{0.48\linewidth}
    
        \centering
        \includegraphics[width=\linewidth]{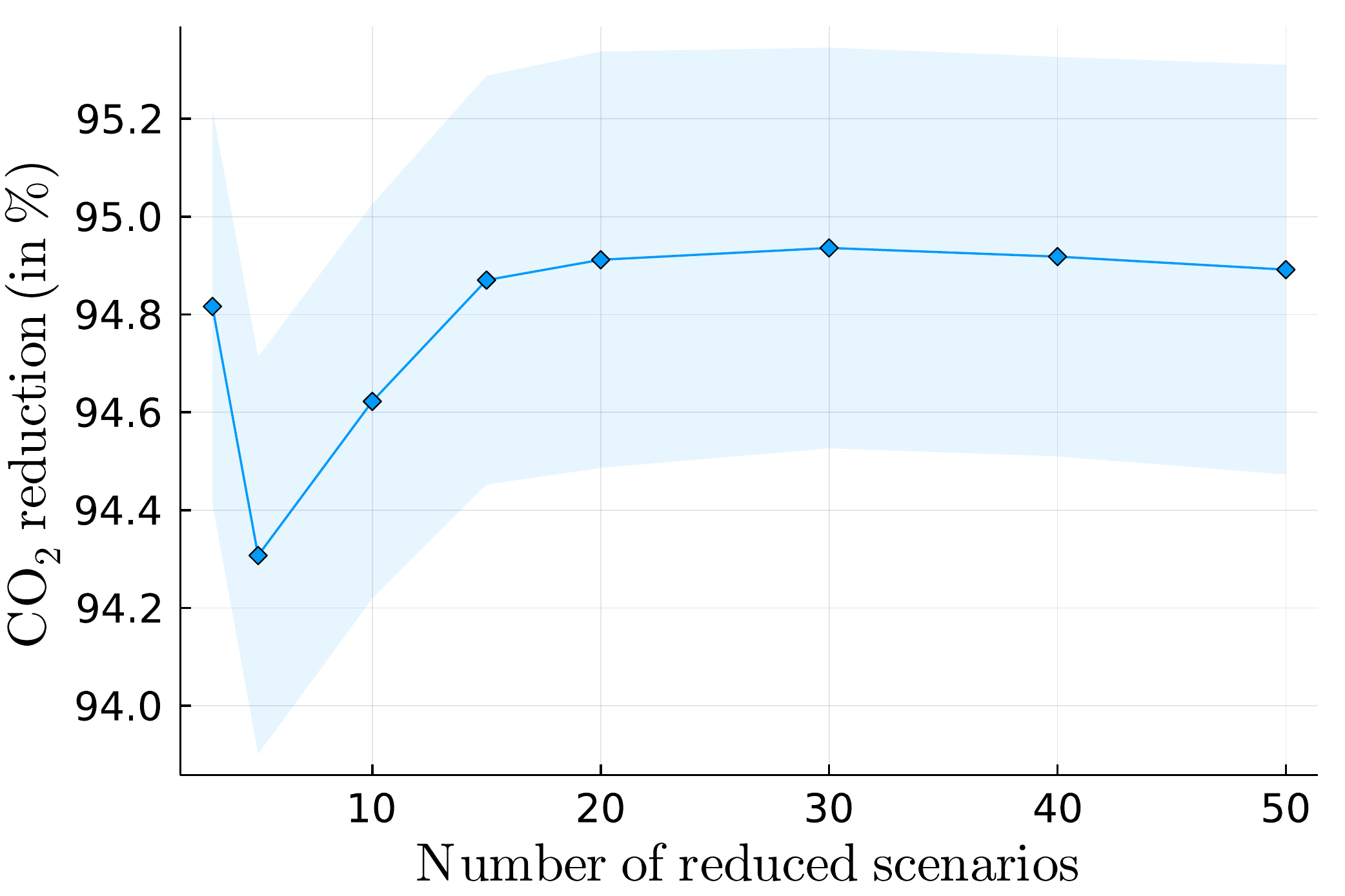}
        \label{fig:n_typical_days_co2}
    
    \end{subfigure}
    \hfill
    \begin{subfigure}[b]{0.48\linewidth}
    
        \centering
        \includegraphics[width=\linewidth]{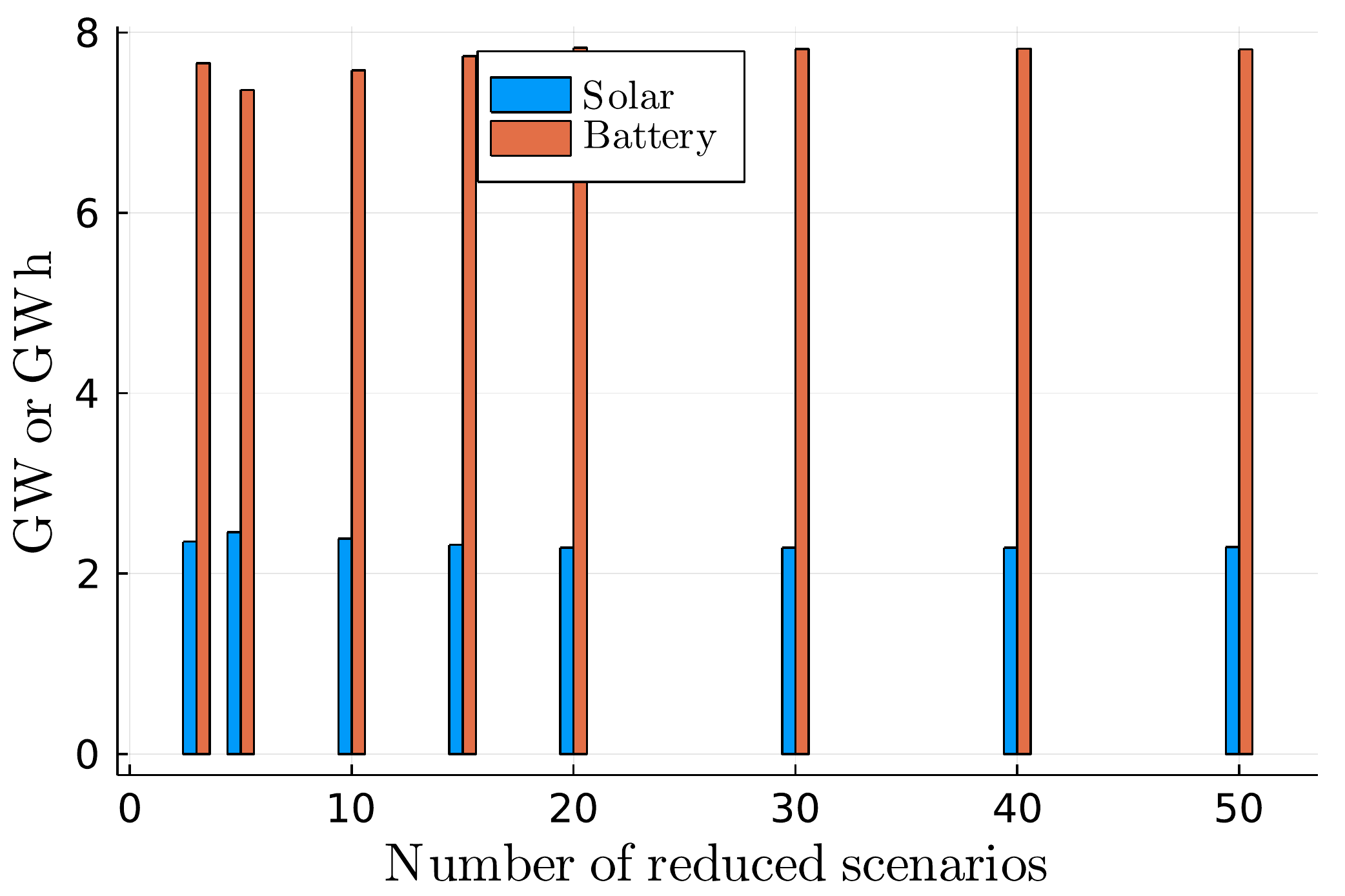}
        \label{fig:n_typical_days_investment}
    
    \end{subfigure}
\end{adjustbox}
\caption{\msom The model is insensitive to the number of scenarios. We report the operational cost (top-left), CO$_2$ emissions reduction (top-right), and investment policy (bottom) as we vary the number of scenarios, with a total investment budget of $20$ billion MAD}
\label{fig:n_typical_days_num_scenarios}
\end{figure}

    
    
    
    

\subsection{Sequences of Reduced Scenarios} 
We now relax the assumption that the amount of energy stored in batteries at the end of hour $24$ of each scenario is the amount of energy available in hour $1$ of the same scenario, and its implication that scenarios can only be succeeded by scenarios of the same type. Specifically, we consider sequences of scenarios of lengths two and three by repeatedly taking the outer product of our set of scenarios with itself. This results in $48$ hour-long ``extended scenarios'' (for a sequence length of two) and $72$ hour-long ``extended scenarios'' (for a sequence length of three), where any scenario can be succeeded by any other scenario and it is possible to shift load between the first and second day in each extended scenario. 

We use a horizon of $20$ years, an investment budget of $20$ billion MAD, and set the RO and DRO hyperparameters to $(\gamma_{\max},\gamma_C,\gamma_{\text{CLT}},\delta)=(0.5,0.5,0.5,0.01)$. In the left part of Figure \ref{fig:sequences}, we consider $3$ single-day scenarios and form sequences of $2$ and $3$ to obtain, respectively, $48$ and $72$ hour-long ``extended scenarios.'' In the right part of Figure \ref{fig:sequences}, we consider $5$ single-day scenarios and form sequences of $2$ to obtain $48$ hour-long ``extended scenarios.'' Our results suggest that using sequences of length $2$ or $3$ rather than single-day scenarios has a negligible impact on the  operational cost, the CO$_2$ emissions reduction, and the optimal investment. This is explained by Morocco's consistent weather patterns from month-to-month and by the decrease in the energy stored in batteries over time, which, in turn, lead to optimal operational policies that do not involve shifting a significant amount of load between days. To further illustrate this point, in Figure \ref{fig:operation-policy-sequences}, we present the operational policy prescribed by the model for two month-extended reduced scenario pairs and verify that the variation in load shifting between different scenarios is indeed negligible.

\begin{figure}[!ht] 
\begin{adjustbox}{minipage=\linewidth,scale=1}
    \centering
    \begin{subfigure}[b]{0.48\linewidth}
    
        \centering
        \includegraphics[width=\linewidth]{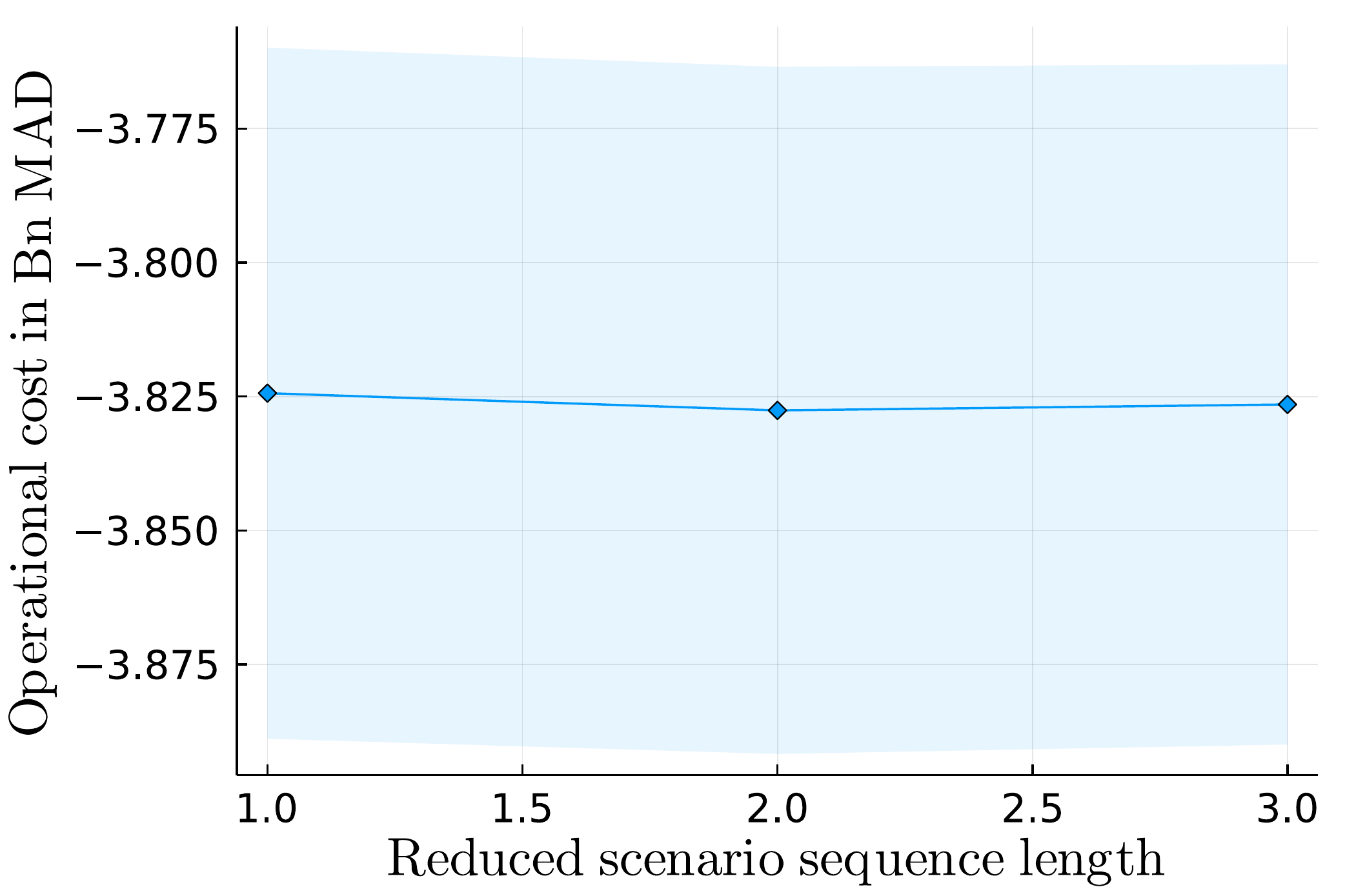}
        \label{fig:sequences_operational_ndays3}
    
    \end{subfigure}
    \hfill
    \begin{subfigure}[b]{0.48\linewidth}
    
        \centering
        \includegraphics[width=\linewidth]{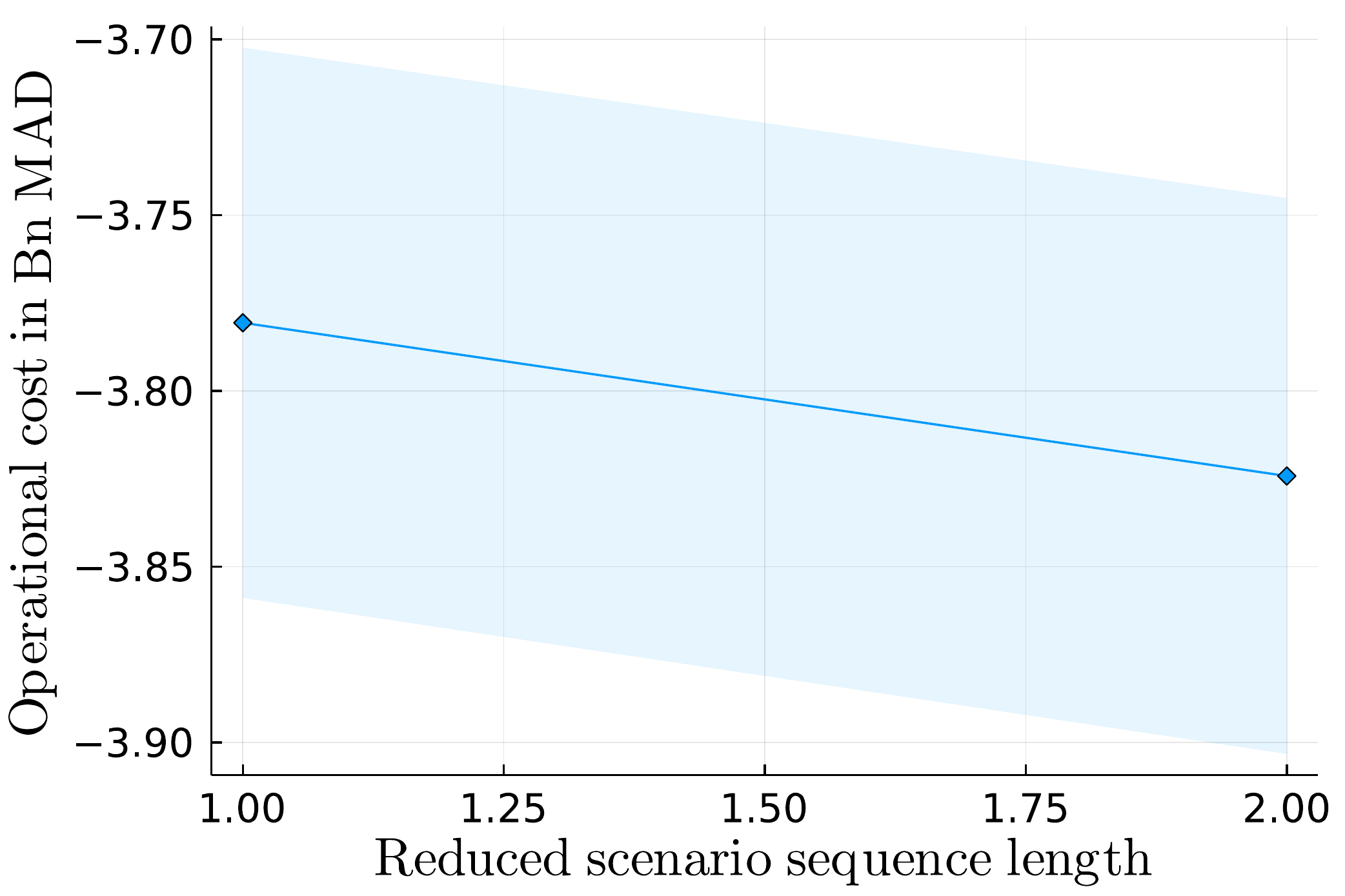}
        \label{fig:sequences_operational_ndays5}
    
    \end{subfigure}
    \hfill
    \begin{subfigure}[b]{0.48\linewidth}
    
        \centering
        \includegraphics[width=\linewidth]{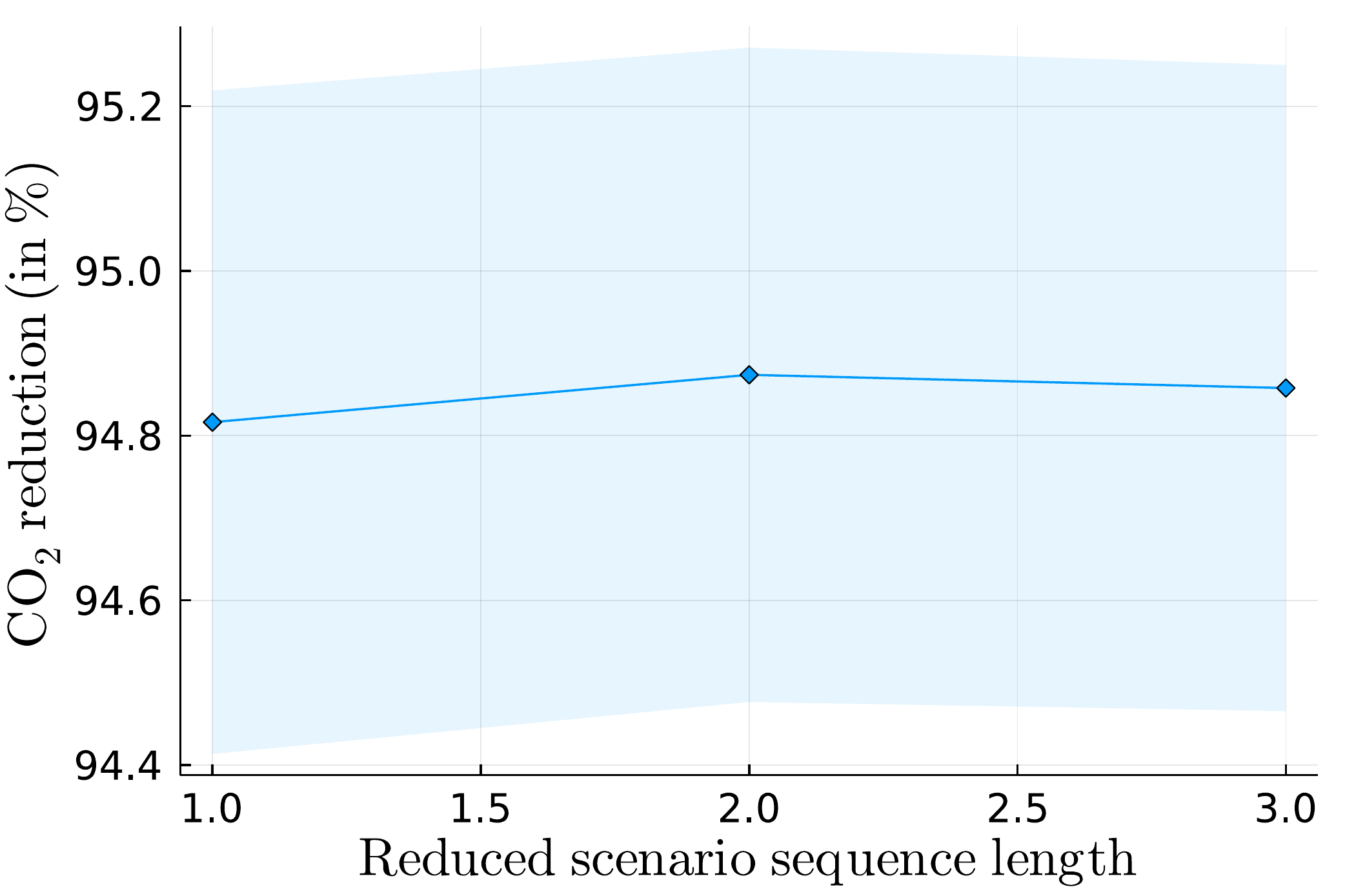}
        \label{fig:sequences_co2_ndays3}
    
    \end{subfigure}
    \hfill
    \begin{subfigure}[b]{0.48\linewidth}
    
        \centering
        \includegraphics[width=\linewidth]{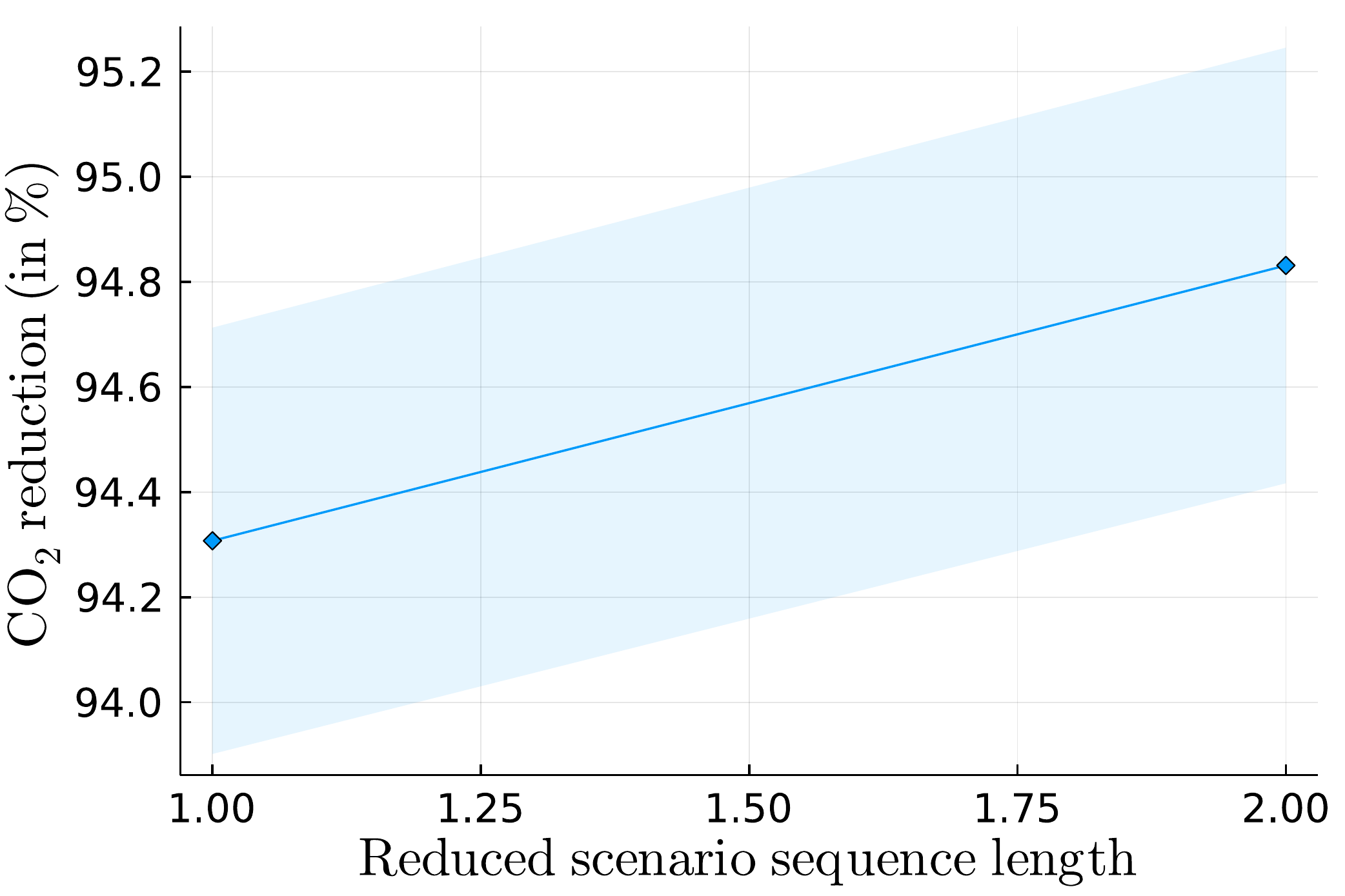}
        \label{fig:sequences_co2_ndays5}
    
    \end{subfigure}
    \hfill
    \begin{subfigure}[b]{0.48\linewidth}
    
        \centering
        \includegraphics[width=\linewidth]{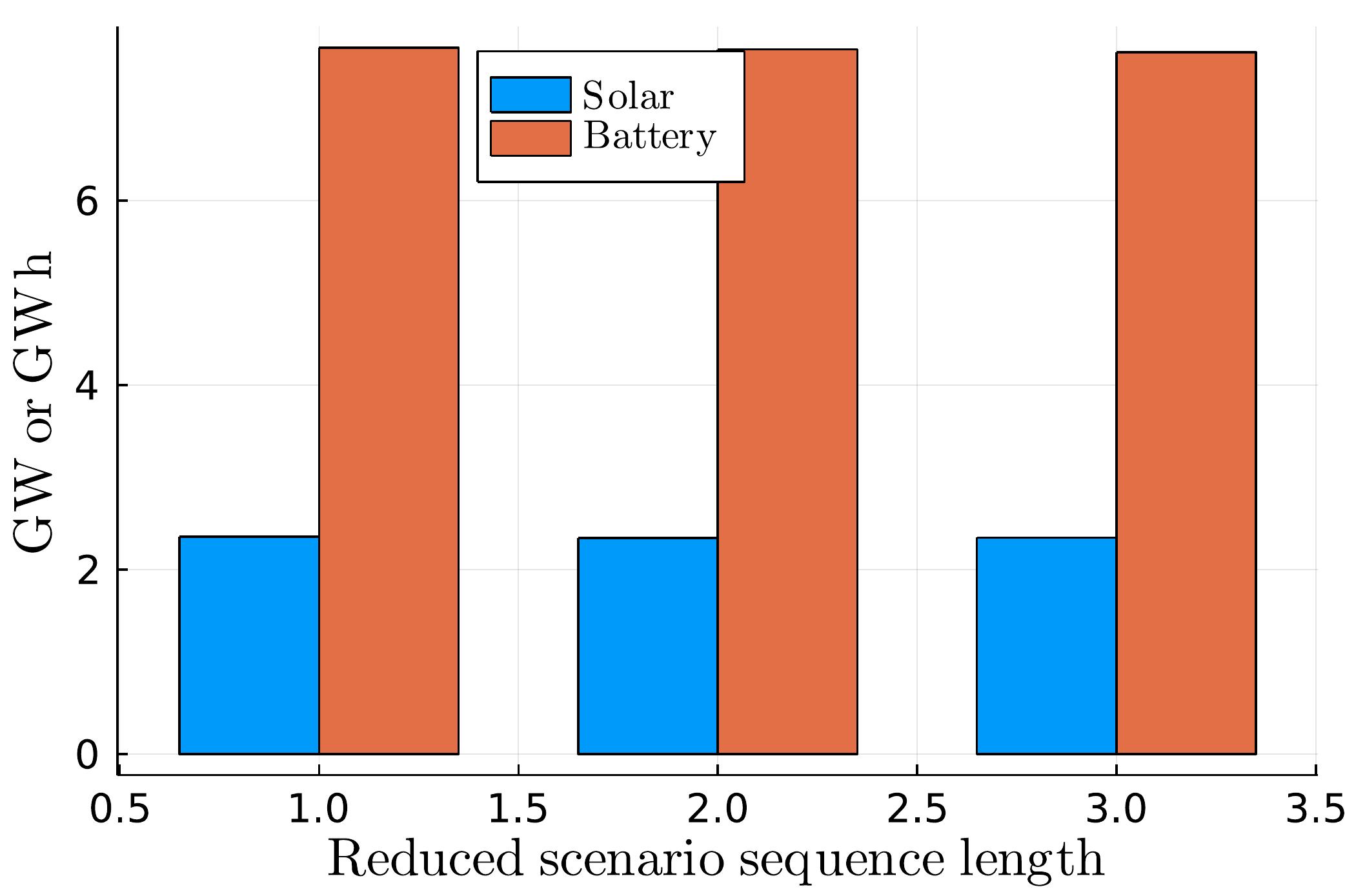}
        \label{fig:sequences_investment_ndays3}
    
    \end{subfigure}
    \hfill
    \begin{subfigure}[b]{0.48\linewidth}
    
        \centering
        \includegraphics[width=\linewidth]{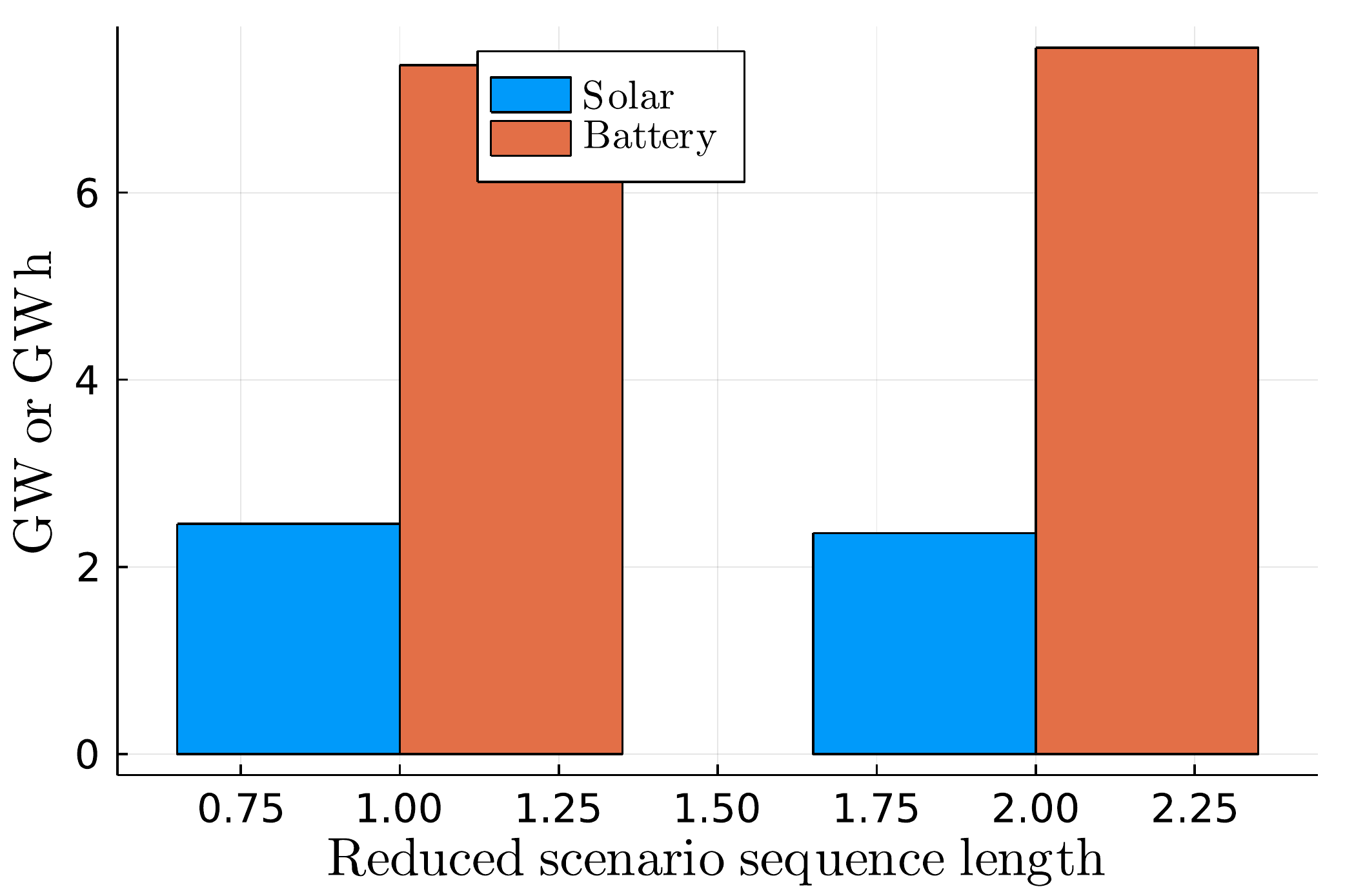}
        \label{fig:sequences_investment_ndays5}
    
    \end{subfigure}
\end{adjustbox}
\caption{\msom Coupling scenarios by using ``extended scenarios'' (formed as sequences of $2$ or $3$ of single-day scenarios) has a negligible impact on the optimal investment. We report the operational cost (with sequences of 3 and 5 reduced scenarios on the top-left and top-right, respectively), CO$_2$ emissions reduction (with sequences of 3 and 5 reduced scenarios on the middle-left and middle-right, respectively), and investment policy (with sequences of 3 and 5 reduced scenarios on the bottom-left and bottom-right, respectively) as we vary the sequence length.}
\label{fig:sequences}
\end{figure}

\begin{figure}[!ht] 
\begin{adjustbox}{minipage=\linewidth,scale=1}
    \centering
    \begin{subfigure}[b]{0.48\linewidth}
    
        \centering
        \includegraphics[width=\linewidth]{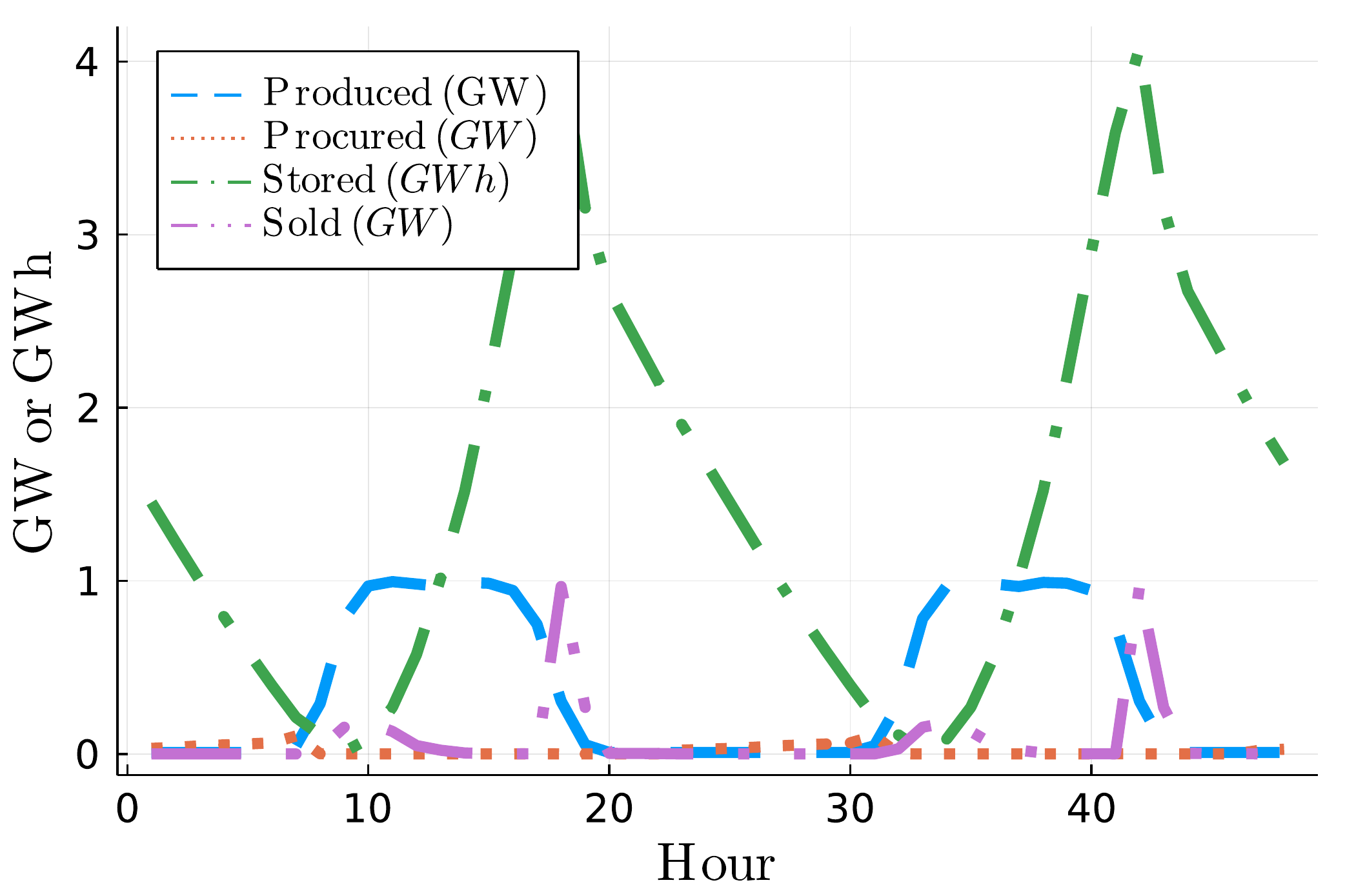}
    
    \end{subfigure}
    \hfill
    \begin{subfigure}[b]{0.48\linewidth}
    
        \centering
        \includegraphics[width=\linewidth]{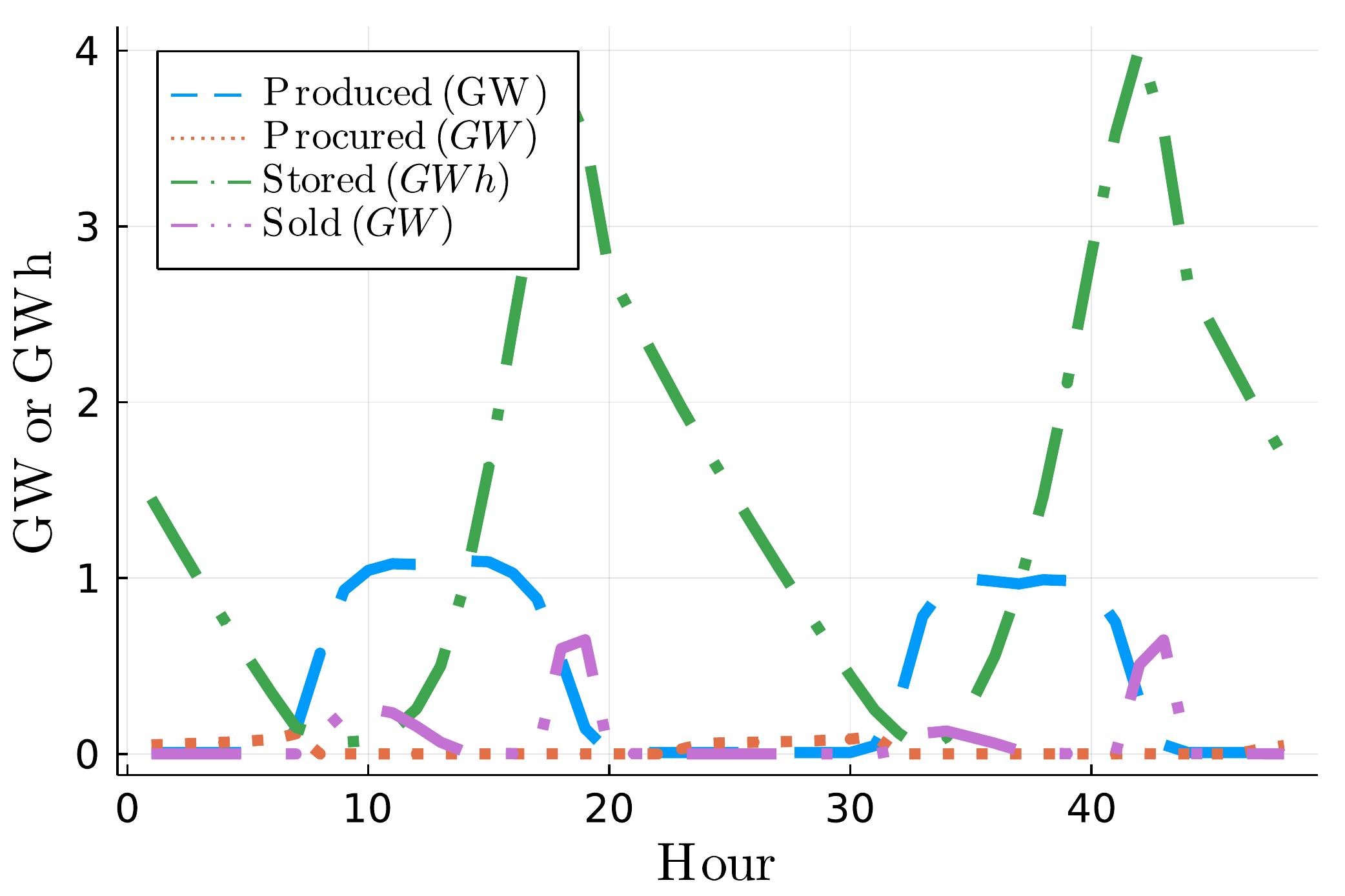}
    
    \end{subfigure}
    \end{adjustbox}
    \caption{{Operational policy prescribed by the model for two month-extended reduced scenario pairs (aggregated across sites and years). We report the policy for a high-generation scenario succeeded by the same scenario (left) and the policy for the same high-generation scenario succeeded by a lower-generation scenario (right).}}
    \label{fig:operation-policy-sequences}
\end{figure}

\section{Reduced Scenarios Visualization} \label{appx:scenarios-visualization}
In this section, we provide a visualization of the estimated reduced scenarios using the methodology we develop in Section \ref{ssec:scenario-reduction}. Figure \ref{fig:rs-vis} presents two aspects of our scenario reduction process:
\begin{itemize}
    \item The left side of the figure illustrates the first part of the process, where each scenario is estimated as a cluster centroid. Each circle in the diagram signifies a 72-dimensional solar capacity factor data point ($\boldsymbol v^i$) from the training set, encapsulating hourly solar capacity factors for a particular day across all three mining sites. These data points are summarized by their mean (x-axis) and standard deviation (y-axis) for the entire day and across all locations. Stars in the figure indicate reduced scenarios or centroids ($\boldsymbol{\bar{v}}^d$). The color coding illustrates the assignment of training points to these scenarios.

    \item The right side of the figure delves into the second part of our scenario reduction process, where we estimate the likelihood of each scenario occurring in different months. Every circle represents an estimated scenario, with the x-axis indicating the most common month of the scenario's occurrence and the y-axis displaying the scenario's mean solar capacity factor. The chart shows that the scenarios with the highest generation commonly occur during the summer.
\end{itemize}

\begin{figure}[h!]
\centering
\begin{subfigure}{.45\textwidth}
  \centering
  \includegraphics[width=\linewidth]{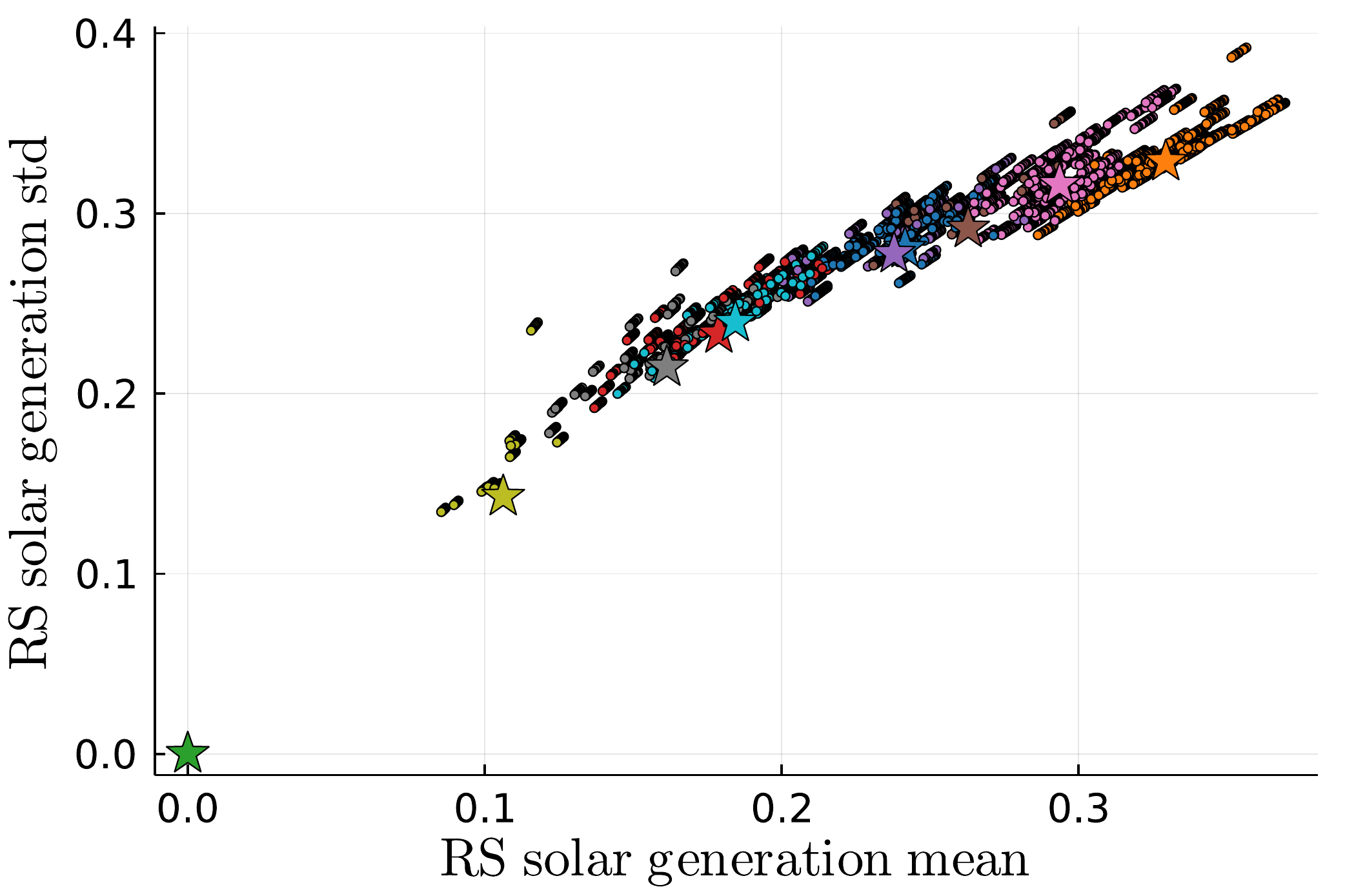}
\end{subfigure}%
\hfill
\begin{subfigure}{.45\textwidth}
  \centering
  \includegraphics[width=\linewidth]{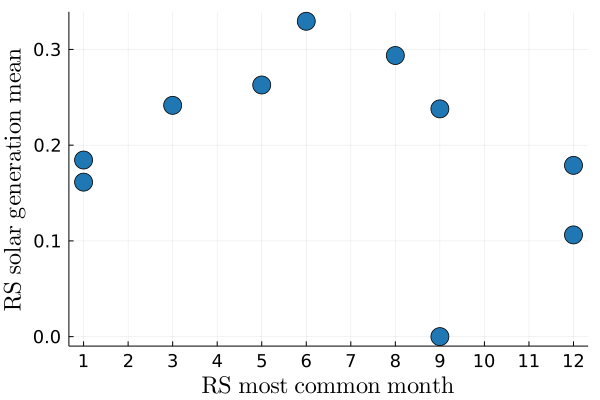}
\end{subfigure}
\vspace{10pt}
\caption{Our two-stage scenario reduction process of Section \ref{ssec:scenario-reduction}: On the left, circles represent 72-dimensional solar capacity factor data points from the training set summarized using their mean and standard deviation for the day across all sites; stars represent the estimated reduced scenarios.
On the right, we examine the likelihood of each scenario's occurrence across months, with the x-axis and y-axis representing the most common occurrence month and the mean solar capacity factor for each scenario, respectively.}
\label{fig:rs-vis}
\end{figure}
}

\end{document}